\documentclass[twoside,11pt]{article}

\usepackage{amsmath}
\usepackage{jmlr2e}
\usepackage{inputenc}
\usepackage[T1]{fontenc}
\usepackage{textcomp}
\usepackage[english]{babel}
\usepackage{graphicx}
\usepackage{enumerate}
\usepackage{dsfont}
\usepackage[normalem]{ulem} %
\usepackage[norelsize,ruled,vlined,commentsnumbered]{algorithm2e}
\usepackage{hyperref}

\def\P{\mathbb{P}}
\def\E{\mathbb{E}} %

\renewcommand{\d}{\mathrm{d}}

\newcommand{\hatN}{\hat{\mathbf{N}}_{M}} 	
\newcommand{\R}{\mathbb{R}}                   		
\newcommand{\C}{\mathbf{C}}            			%

\newcommand{\mbf}[1]{\mathbf #1}
\newcommand{\eq}{\begin{equation}}
\newcommand{\qe}{\end{equation}}

\newcommand{\hatl}{\hat{\mbf L}_{M}}

\newcommand{\B}{\mathbb B}

\newcommand{\X}{\mathcal X}	
\newcommand{\dS}{\mathcal S}				
\newcommand{\Y}{\mathcal Y}				
\renewcommand{\L}{\mathcal L^{D}}					
\newcommand{\Proj}{\mathfrak P}				

\newcommand{\Q}{\mathbf Q}				
\newcommand{\F}{\mathfrak F}				
\newcommand{\M}{\mathbf{M}_{M}} 			
\newcommand{\hatM}{\hat{\mathbf{M}}_{M}}	
\newcommand{\U}{\mathbf{U}}

\newcommand{\f}{\mathbf{f}}
\newcommand{\h}{\mathbf{h}}

\newcommand{\Rbf}{\mathbf{R}}
     
\newcommand{\one}{\mathds{1}}
\renewcommand{\O}{\mathbf{O}_{M}}			
\newcommand{\hatO}{\hat{\mathbf{O}}_{M}}	

\renewcommand{\B}{\mathbf{B}}

\newcommand{\D}[1]{\mathfrak{Diag}[{#1}]}        
\newcommand{\hatPj}{\hat{\mathbf P}_{M}}		

\newcommand{\pen}{{\rm pen}}

\jmlrheading{}{2015}{}{07/15}{??/??}{Yohann De Castro, \'Elisabeth Gassiat and Claire Lacour}

\ShortHeadings{Estimation of nonparametric HMMs}{De Castro, Gassiat and Lacour}
\firstpageno{1}

\begin{document}

\title{Minimax adaptive estimation of nonparametric hidden Markov models}

\author{\name Yohann De Castro \email yohann.decastro@math.u-psud.fr \\
       \addr Laboratoire de Math\'ematiques d'Orsay, Univ. Paris-Sud, CNRS, Universit\'e Paris-Saclay,\\ 91405 Orsay, France.
       \AND
       \name \'Elisabeth Gassiat \email elisabeth.gassiat@math.u-psud.fr\\
             \addr Laboratoire de Math\'ematiques d'Orsay, Univ. Paris-Sud, CNRS, Universit\'e Paris-Saclay, \\91405 Orsay, France.       \AND
       \name Claire Lacour \email claire.lacour@math.u-psud.fr\\
             \addr Laboratoire de Math\'ematiques d'Orsay, Univ. Paris-Sud, CNRS, Universit\'e Paris-Saclay,\\ 91405 Orsay, France.
     }

\editor{}

\maketitle

\begin{abstract}
We consider stationary hidden Markov models with finite state space and nonparametric modeling of the emission distributions. 
It has remained unknown until very recently that such models are identifiable. In this paper, we propose a new penalized least-squares estimator for the emission distributions
which is statistically optimal and practically tractable.  We prove a non asymptotic oracle inequality for our nonparametric estimator of the emission distributions. A consequence is that  this new estimator is rate minimax adaptive up to a logarithmic term.
Our methodology is based on projections of the emission distributions onto nested subspaces of increasing complexity. The popular spectral estimators are unable to achieve the optimal rate but may be used as initial points in our procedure. 
Simulations are given that show the improvement obtained when applying the least-squares minimization consecutively to the spectral estimation.
\end{abstract}


\section{Introduction}

\subsection{Context and motivations}
Finite state space hidden Markov models (HMMs for short)  are widely used to model data evolving in time  and coming from heterogeneous populations.
They seem to be reliable tools to model practical situations in a variety of applications such as economics, genomics, signal processing and image analysis, ecology, environment,  speech recognition, to name but a few.  
From a statistical view point, finite state space HMMs are stochastic processes $(X_j,Y_j)_{j\geq 1}$ where $(X_j)_{j\geq 1} $ is a Markov chain with finite state space    and conditionally on $(X_j)_{j\geq 1}$ the $Y_j$'s are independent with a distribution depending only on $X_j$. 
The observations are $Y_{1:N} = (Y_1,\ldots,Y_N)$ and the associated states $X_{1:N} = (X_1,\ldots,X_N)$ are unobserved. 
The parameters of the model are the initial distribution, the transition matrix of the hidden  chain, and the emission distributions of the observations, that is the probability distributions of the $Y_{j}$'s conditionally to $X_{j}=x$ for all possible $x$'s. In this paper we shall consider stationary ergodic HMMs so that the initial distribution is the stationary distribution of the (ergodic) hidden Markov chain.

Until very recently, asymptotic performances of estimators were  proved 
only in the parametric setting (that is, with finitely many unknown parameters). Though, nonparametric methods for HMMs have been considered in applied papers, but with no theoretical guarantees,
see for instance \cite{CC00} for voice activity detection, \cite{LWM03} for climate state identification, \cite{Lef03} for automatic speech recognition, \cite{SC09} for facial expression recognition, \cite{VBM13} for methylation comparison of proteins, \cite{MR2797735}  for copy number variants identification in DNA analysis. 

The preliminary obstacle to obtain theoretical results on general finite state space nonparametric HMMs was to understand when such models are indeed  identifiable. 
Marginal distributions of finitely many observations are finite mixtures of products of the emission distributions. It is clear that identifiability can not be obtained based on the marginal distribution of only one observation. It is needed, and it is enough,  to consider the marginal distribution of at least three consecutive observations to get identifiability, see  \cite{gassiat2013finite}, following \cite{allman2009identifiability} and \cite{hsu2012spectral}.  

\subsection{Contribution}

The aim of our paper is to propose a 
new approach to estimate nonparametric HMMs with a statistically optimal and practically tractable method. We obtain this way
nonparametric estimators of the emission distributions that achieve the minimax rate of estimation in an adaptive setting. 

Our perspective is based on estimating the projections of the emission laws onto nested subspaces of increasing complexity. Our analysis encompasses any family of nested subspaces of Hilbert spaces and works with a large variety of models. In this framework one could think to use the spectral estimators as proposed by \cite{hsu2012spectral,anandkumar2012method}
 in the parametric framework, by extending them to the nonparametric framework. But a careful analysis of the tradeoff between sampling size and approximation complexity shows that they do not lead to rate optimal estimators of the emission densities, see \cite{YCES} for a formal statement and proof. This can be easily understood. Indeed,  the spectral estimators of the emission densities are computed as functions of the empirical estimator of the marginal distribution of three consecutive observations on $\Y^3$ (where $\Y$ is the observation space), for which, roughly speaking, when  $\Y$ is a subset of $\R$, the optimal rate is $N^{-s/(2s+3)}$, $N$ being the number of observations and $s$ the smoothness of the emission densities. Thus the rate obtained this way for the emission densities is also $N^{-s/(2s+3)}$. But since those emission densities describe one dimensional random variables on $\Y$, one could hope to be able to obtain the sharper rate $N^{-s/(2s+1)}$. This is the rate we obtain, up   to a $\log N$ term, with our new method. Let us explain how it works.

Using  
the HMM modeling, and using sieves for the emission densities on $\Y$, we propose a penalized least squares estimator  in the model selection framework. 
We prove an oracle inequality  for the $L_{2}$-risk of the estimator of the density of $(Y_{1},Y_{2},Y_{3})$, see Theorem \ref{prop:oracle}. Since the complexity of the model is that given by the sieves for the emission densities, this leads, up to a $\log N$ term, to the adaptive minimax rate computed as for the density of only one observation $Y_1$ though we estimate the density of $(Y_{1},Y_{2},Y_{3})$. Roughly speaking, when the observations are one dimensional, that is when  $\Y$ is a subset of $\R$, the obtained rate for the density of $(Y_{1},Y_{2},Y_{3})$ is of order $N^{-s/(2s+1)}$ up to a $\log N$ term, $N$ being the number of observations and $s$ the smoothness of the emission densities. 

The key point is then to be able to go back to the emission densities. This is the cornerstone of our main result. We prove in Theorem \ref{LEMMA:INEQ2} that, under the assumption {\bf [HD]} defined in Section \ref{sec:MAIN}, the quadratic risk for the density of $(Y_{1},Y_{2},Y_{3})$ is lower bounded by some positive constant multiplied by the quadratic risk for the emission densities. This technical assumption is generically satisfied in the sense that it holds for all possible emission densities for which the $L_2$-norms and Hilbert dot products do not lie on a particular algebraic surface with coefficients depending on the transition matrix of the hidden chain. Moreover, we prove that,
 when the number of hidden states equals two, this assumption is always verified when the two emission densities are distinct, see Lemma \ref{LEMMA:INEQ3}.
 
 Our methodology requires 
 that we have a preliminary estimator of the transition matrix. 
To get such an estimator, it is possible to use spectral methods. 
Thus our approach is the following. First, get a preliminary estimator of the initial distribution and the transition matrix of the hidden chain. Second, apply penalized least squares estimation on the density of three consecutive observations, using HMM modeling, model selection on the emission densities, and initial distribution and stationary matrix of the hidden chain set at the estimated value. This gives emission density estimators which have minimax adaptive rate, as our main result states, see Theorem \ref{theo:adapt}. A simplified version of this theorem can be given as follows.
\begin{theorem} Assume $(Y_{j})_{j\geq 1}$ is a hidden Markov model on $\mathbb{R}$, with latent Markov chain
$(X_{j})_{j\geq 1}$ with $K$ possible values and true transition matrix $\Q^{\star}$. Denote $f^{\star}_k$ the density of
$Y_n$ given $X_n=k$, for $k=1,\ldots,K$.
 Assume the true transition matrix $\Q^{\star}$ is full rank and the true emission densities $f^{\star}_k$, $k=1,\ldots,K$ are linearly independent, with smoothness $s$.
 Assume that {\bf [HD]} holds true.
Then, up to label switching, for $N$ the number of observations large enough,  the estimators $\hat\Q, \hat{f}_{k}$, $k=1,\ldots,K$ built in Section \ref{sec:lsmethod} and \ref{sec:spectral} satisfy
\begin{align*}
\E\left[ \|\Q^{\star}-\hat\Q\|^2\right]&=O\Big(\frac{\log N}{N}\Big),\\
\E\left[\|f^{\star}_{k}-\hat{f}_{k}\|_{2}^{2}\right]&=
O\Big(\big[\frac{\log N}{N}\big]^{\frac s{2s+1}}\Big),\;k=1,\ldots,K.
\end{align*}
\end{theorem}

\noindent
Moreover, since the family of sieves we consider is that given by finite dimensional spaces described by an orthonormal basis, we are able to use the spectral estimators of the coefficients of the densities as initial points in the least squares minimization. This is important since, in the HMM framework, least squares minimization does not have an explicit solution and may lead to several local minima. However, since the spectral estimates are proved to be consistent, we may be confident that their use as initial point is enough. Simulations indeed confirm this point.

To conclude we claim that our results support a powerful new approach to estimate,  for the first time, nonparametric HMMs with a statistically optimal and practically tractable method.

\subsection{Related works}
The papers  \cite{allman2009identifiability}, \cite{hsu2012spectral} and \cite{anandkumar2012method} paved the way to obtain identifiability  under reasonable assumptions. In \cite{anandkumar2012method} the authors point out  a structural link between multivariate mixtures with conditionally independent observations and finite state space HMMs.  In \cite{hsu2012spectral} the authors propose a spectral method to estimate all parameters for finite state space HMMs (with finitely many observations), under the assumption that the transition matrix of the hidden chain is non singular, and that the (finitely valued) emission distributions are linearly independent. Extension to emission distributions on any space, under the  linear independence assumptions (and keeping the assumption of non singularity of the transition matrix), allowed to prove the general identifiability result for finite state space HMMs, see \cite{gassiat2013finite}, where also model selection likelihood methods and nonparametric kernel methods are proposed to get nonparametric estimators. Let us  notice also \cite{ver13} that proves theoretical consistency of the posterior in nonparametric Bayesian methods for finite state space HMMs with adequate assumptions. Later, \cite{alexandrovich2014nonparametric} obtained identifiability when the emission distributions are all distinct (not necessarily linearly independent) and still when the transition matrix of the hidden chain is full rank.
In the nonparametric multivariate mixture model,  \cite{song2013nonparametric} prove that any linear functional of the emission distributions may be estimated with parametric rate of convergence in the context of reproducing kernel Hilbert spaces.
The latter uses spectral methods, not the same but similar to the ones proposed in \cite{hsu2012spectral} and \cite{anandkumar2012method}. 

Recent papers that contain theoretical results on different kinds of nonparametric HMMs are \cite{gassiat:rousseau:13}, where the emitted distributions are translated versions of each other, and \cite{TDSLC} in which the authors consider regression models with hidden regressor variables that can be Markovian on a continuous state space.

\subsection{Outline of the paper}
In Section \ref{sec:notations}, we set {the} notations, the model we shall study, and the assumptions we shall consider. We then state an identifiability lemma (see Lemma \ref{identifiabilityL2}) that will be useful for our estimation method. 
In Sections \ref{sec:lsmethod} and  \ref{subsec:main1} we give our main results. We explain the penalized least-squares estimation method in Section \ref{sec:lsmethod}, and we prove in Section \ref{subsec:main1} that, 
when the transition matrix is irreducible and aperiodic, when the emission distributions are linearly independent and the penalty is adequately chosen, then, under a technical assumption, the penalized least squares estimator is asymptotically minimax adaptive up to a $\log N$ term, see Theorem~\ref{theo:adapt} and Corollary~\ref{cor:adapt}. For this, we first prove an oracle inequality for the estimation of the density of $(Y_{1},Y_{2},Y_{3})$, see Theorem \ref{prop:oracle}, then we prove the key result relating the risk of the density of $(Y_{1},Y_{2},Y_{3})$ to that of the emission densities, see Theorem \ref{LEMMA:INEQ2}. The latter holds under a technical assumption which we prove to be always verified in case $K=2$, see Lemma
\ref{LEMMA:INEQ3}.
Finally, we need  the performances of the spectral estimator of the transition matrix and of the stationary distribution which are given in Section \ref{sec:spectral}, see Theorem \ref{thm:spectral}, proved in \cite{YCES}. We finally present simulations in Section \ref{sec:simul} to illustrate our theoretical results. Those simulations show in particular the  improvement obtained when applying the least-squares minimization consecutively to the spectral estimation. Detailed proofs are given in Section \ref{sec:proofs}.

\section{Notations and assumptions}
\label{sec:notations}

\subsection{Nonparametric hidden Markov model}
Let $K$, $D$ be positive integers and let $\L$ be the Lebesgue measure on $\R^{D}$. Denote by $\X$ the set $\{1,\ldots,K\}$ of hidden states, $\Y\subset\mathbb R^{D}$ the observation space, and $\Delta_{K}$ the space of probability measures on $\X$ identified to the $(K-1)$-dimensional simplex. Let $(X_n)_{n\geq1}$ be a Markov chain on $\X$ with $K\times K$ transition matrix $\Q^{\star}$ and initial distribution $\pi^{\star}\in\Delta_{K}$.  Let $(Y_{n})_{n\geq1}$ be a sequence of observed random variables on $\Y$. Assume that, conditional on $(X_n)_{n\geq1}$, the observations $(Y_{n})_{n\geq1}$ are independent and, for all $n\in\mathbb N$, the distribution of $Y_{n}$ depends only on $X_{n}$. 
Denote by $\mu_{k}^{\star}$ the conditional law of $Y_{n}$ conditional on $\{X_{n}=k\}$, and assume that $\mu_{k}^{\star}$ has density $f_{k}^{\star}$ with respect to the measure $\L$ on $\Y$:
\[
\forall k\in\X\,,\quad \d\mu_{k}^{\star}=f_{k}^{\star}\d\L\,.
\]
Denote by $\mathfrak F^{\star}:=\{f_{1}^{\star},\ldots,f_{K}^{\star}\}$ the set of emission densities with respect to the Lebesgue measure. Then, for any integer $n$, the distribution of $(Y_{1},\ldots,Y_{n})$ has density with respect to $(\L)^{\otimes n}$
\[
\sum_{k_1,\ldots,k_n=1}^K \pi^{\star}(k_1) \Q^{\star}(k_1,k_2)\ldots\Q^{\star}(k_{n-1},k_n)f^{\star}_{k_1}(y_1)\ldots f^{\star}_{k_n}(y_n).
\]
We shall denote $g^{\star}$ 
 the density of $(Y_1,Y_2,Y_3)$.
 
 In this paper we shall address two observations schemes. 
We shall  consider $N$ i.i.d. samples $(Y_{1}^{(s)},Y_{2}^{(s)},Y_{3}^{(s)})_{s=1}^{N}$ of three consecutive observations {\bf (Scenario A)} or consecutive observations of the same chain {\bf (Scenario B)}:
\[
\forall s\in\{1,\ldots,N\},\quad (Y_{1}^{(s)},Y_{2}^{(s)},Y_{3}^{(s)}):=(Y_{s},Y_{s+1},Y_{s+2})\,.
\]
 
\subsection{Projections of the population joint laws} 
\label{sec:ObservableSpace}
Denote by $(\mathbf{L}^{2}(\Y,\L),\lVert \cdotp \lVert_{2})$ the Hilbert  space of square integrable functions on $\Y$ with respect to the Lebesgue measure $\L$ equipped with the usual inner product $\langle \cdotp,\cdotp\rangle$ on $\mathbf{L}^{2}(\Y,\L)$. Assume $\mathfrak F^{\star}\subset\mathbf{L}^{2}(\Y,\L)$.

Let $(M_{r})_{r\geq1}$ be an increasing sequence of integers, and let $(\Proj_{M_{r}})_{r\geq 1}$ be a sequence of nested subspaces with dimension $M_{r}$ such that their union is dense in $\mathbf{L}^{2}(\Y,\L)$. Let $\Phi_{M_{r}}:=\{\varphi_{1},\ldots,\varphi_{M_{r}}\}$ be an orthonormal basis of $\Proj_{M_{r}}$. Recall that for all $f\in\mathbf{L}^{2}(\Y,\L)$,
\eq
\label{eq:ConvergenceL2}
\lim_{r\to\infty}\sum_{m=1}^{M_{r}}\langle f,\varphi_{m}\rangle\varphi_{m}= f\,,
\qe
in $\mathbf{L}^{2}(\Y,\L)$. Note that changing $M_r$ may change all functions $\varphi_{m}$, $1\leq m \leq M_r$ in the basis $\Phi_{M_r}$, which we shall not indicate in the notation for sake of readability. Also, we drop the dependence on $r$ and write $M$ instead of $M_{r}$. Define the projection of the emission laws onto $\Proj_{M}$ by
\[
\forall k\in\X,\quad  f^{\star}_{M,k}:=\sum_{m=1}^{M}\langle f_{k}^{\star},\varphi_{m}\rangle\varphi_{m}\,.
\]
We shall write $\f^{\star}_{M}:=( f^{\star}_{M,1},\ldots, f^{\star}_{M,K})$ and $\f^{\star}:=( f^{\star}_{1},\ldots, f^{\star}_{K})$ throughout this paper.

\begin{remark}
One can consider the following standard examples:
\begin{enumerate}
\item[\bf (Spline)] The space of piecewise polynomials of degree bounded by $d_r$ based on the regular partition with $p_{r}^{D}$ regular pieces on $\Y=[0,1]^{D}$. 
It holds that $M_{r}=(d_{r}+1)^{D}p_{r}^{D}$. 
\item[\bf (Trig.)] The space of real trigonometric polynomials on $\Y=[0,1]^{D}$ with degree less than $r$. 
It holds that $M_{r}=(2r+1)^{D}$.
\item[\bf (Wav.)] A wavelet basis $\Phi_{M_{r}}$ of scale $r$ on $\Y=[0,1]^{D}$,  see
\cite{MR1228209}. 
In this case, it holds that $M_{r}=2^{(r+1)D}$.
\end{enumerate}
\end{remark}

\subsection{Assumptions}

We shall use the following assumptions on the hidden chain. 
\begin{enumerate}[{\bf [H1]}]
\item
{\it
The transition matrix $\Q^{\star}$ has full rank,
}
\item
{\it
The Markov chain $(X_n)_{n\geq1}$ is irreducible and aperiodic,
}
\item
\label{ass:PiMin}
{\it
The initial distribution $\pi^{\star}=(\pi^{\star}_{1},\ldots,\pi^{\star}_{K})$ is the stationary distribution}.
\end{enumerate}
Notice that under {\bf [H1]}, {\bf [H2]} and {\bf [H3]}, one has for all $k\in\X$, $\pi^{\star}_{k}\geq\pi^{\star}_{\mathrm{min}}>0$. We shall use the following assumption on the emission densities.
\begin{enumerate}[{\bf [H4]}]
\item
{\it
The family of emission densities $\F^{\star}:=\{f^{\star}_{1},\ldots,f^{\star}_{K}\}$ is linearly independent.
}
\end{enumerate}

\noindent
Those assumptions  appear in spectral methods, see for instance \cite{hsu2012spectral,anandkumar2012method}, and in identifiability issues, see for instance \cite{allman2009identifiability,gassiat2013finite}.

\subsection{Identifiability lemma}
For any $\f=(f_{1},\ldots,f_{K})\in (\mathbf{L}^{2}(\Y,\L))^{K}$ and any  transition matrix  $\Q$, denote by $g^{\Q,\f}:\Y^{3}\rightarrow \R$ the function given by
\begin{equation}
\label{eq:gQf}
g^{\Q,\f}\left(y_{1},y_{2},y_{3}\right)=\sum_{k_1,k_{2},k_3=1}^K \pi(k_1) \Q(k_1,k_2)\Q(k_{2},k_3)f_{k_1}(y_1) f_{k_2}(y_2) f_{k_3}(y_3),
\end{equation}
where $\pi$ is the stationary distribution of $\Q$. When $\Q=\Q^{\star}$ and $\f=\f^{\star}$, we get $g^{\Q^{\star},\f^{\star}}=g^{\star}$. When $f_{1},\ldots,f_{K}$ are probability densities on $\Y$, $g^{\Q,\f}$ is the probability distribution of three consecutive observations of a stationary HMM. We now state a lemma that gathers all what we need about identifiability.

For any transition matrix  $\Q$, let $T_{\Q}$ be the set of permutations $\tau$  such that for all $i$ and $j$, $\Q (\tau(i),\tau(j))=\Q(i,j)$.
The permutations in $T_{\Q}$ describe how the states of the Markov chain may be permuted without changing the distribution of the whole chain: for any $\tau$ in $T_{\Q}$, $(\tau (X_{n}))_{n\geq 1}$ has the same distribution as $(X_{n})_{n\geq 1}$. Since the hidden chain is not observed, if the emission distributions are permuted using $\tau$, we get the same HMM. In other words, if $\f^{\tau}=(f_{\tau(1)},\ldots,f_{\tau(K)})$, then
$
g^{\Q,\f^{\tau}}=g^{\Q,\f}.
$
Since identifiability up to permutation of the hidden states is obtained from the marginal distribution of three consecutive observations, we get the following lemma whose detailed proof is given in Section \ref{sec:proof-iden}.
\begin{lemma}
\label{identifiabilityL2}
Assume that $\Q$ is a transition matrix for which {\bf [H1]} and {\bf [H2]} hold. Assume that {\bf [H4]} holds.
Then for any $\h\in(\mathbf{L}^{2}(\Y,\L))^{K}$, 
\[
g^{\Q,\f^{\star}+\h}=g^{\Q,\f^{\star}}\Longleftrightarrow \exists \tau\in T_{\Q}\;{\text{such that}}\;h_{j}=f^{\star}_{\tau(j)}-f^{\star}_{j},\;j=1,\ldots,K.
\]
In particular, if $T_{\Q}$ reduces to the identity permutation,
\[
g^{\Q,\f^{\star}+\h}=g^{\Q,\f^{\star}}\Longleftrightarrow \h=(0,\ldots,0).
\]
\end{lemma}

\section{The penalized least-squares estimator}
\label{sec:lsmethod}

In this section we shall estimate the emission densities using the so-called penalized least squares method. Here, the least squares adjustment is made
on the density $g^\star$ of $(Y_1,Y_2,Y_3)$. Starting from the operator $\Gamma:\ t\mapsto\|t-g^{\star}\|_{2}^{2}-\|g^{\star}\|_{2}^{2}=\|t\|_{2}^2-2\int t g^{\star}$ which is minimal for the target $g^\star$, we introduce the corresponding empirical contrast $\gamma_{N}$. Namely, for any $t\in \mathbf{L}^{2}(\Y^{3},{\L}^{\otimes 3})$, set
\[
\gamma_{N}\left(t\right)=\|t\|_{2}^{2}-\frac{2}{N}\sum_{s=1}^{N}t\left(Z_{s}\right),
\]
with $Z_{s}:=(Y^{(s)}_{1},Y^{(s)}_{2},Y^{(s)}_{3})$ ({\bf{Scenario A}}) or $Z_{s}:=(Y_{s},Y_{s+1},Y_{s+2})$ ({\bf{Scenario B}}).
As $N$ tends to infinity, $\gamma_{N}\left(t\right)-\gamma_{N}\left(g^{\star}\right)$ converges almost surely to $\|t-g^{\star}\|_{2}^{2}$, thus the name least squares contrast function.
A natural estimator is then a function $t$ such that $\gamma_{N}\left(t\right)$ is minimal over a judicious approximation space which is a set of functions of form $g^{\Q,\f}$, $\Q$ a transition matrix and $\f\in{\mathcal F}^{K}$,
for $\mathcal F$ a subset of $\mathbf{L}^{2}(\Y,\L)$. We thus define a whole collection of estimates $\hat g_M$, each $M$ indexing an approximation subspace 
(also called model). 
Considering (\ref{eq:gQf}) we shall introduce a collection of
 model of functions by projection of possible $\f$'s on the subspaces $(\Proj_{M})_{M}$. 
Thus,
for any irreducible transition matrix $\Q$ with stationary distribution $\pi$, we define $\dS ({\Q},{M})$ as the set of functions $g^{\Q,\f}$
 such that $\f\in{\mathcal F}^{K}$ and, for each $k=1,\ldots,K$, 
 there exists $(a_{m ,k})_{1\leq m \leq M}\in\R^{M}$ such that
\[
 f_{k}=\sum_{m=1}^{M}a_{m, k}\varphi_{m}.
\]

\noindent
We now assume that we have in hand an estimator $\hat \Q$ of $\Q^{\star}$. 
For instance, one can use a spectral estimator, we recall such a construction in Section \ref{sec:spectral}. Then, $(\dS ({\hat \Q},{M}))_{M}$ 
is the collection of models we use for the least squares minimization.
For any $M$, define $\hat g_M$ as a minimizer of $\gamma_{N}(t)$ for $t\in \dS ({\hat \Q},{M})$.
Then $\hat g_M$ can be written as $\hat g_M=g^{\hat\Q,\hat{\f}_{M}}$ with $\hat{\f}_{M}\in{\mathcal F}^{K}$ and $\hat{f}_{M,k}=\sum_{m=1}^M\hat a_{m, k}\varphi_{m}$ ($k=1,\ldots,K$) for some $(\hat a_{m,k})_{1\leq m \leq M}\in\R^{M}$, $k=1,\ldots,K$.
It then remains to select the best model, that is to choose $M$ which minimizes
$\|\hat g_M-g^{\star}\|_{2}^{2}-\|g^{\star}\|_{2}^{2}$.  This quantity is close to $\gamma_{N}(\hat g_M)$, but we need to take into account the deviations of the process $\Gamma-\gamma_N$. Then we rather minimize
$\gamma_{N}(\hat g_M) + \pen(N,M)$ where $\pen(N,M)$ is a penalty term to be specified. 
Our final estimator will be a penalized least squares estimator. For this purpose we choose a penalty function $\pen(N,M)$ and 
we let 
\[
\hat M=\arg \min_{M=1,\ldots,N}\left\{\gamma_{N}(\hat g_M)+\pen(N,M)\right\}.
\]
Notice that, with $N$ observations, we consider $N$ subspaces as candidates for model selection. Then the estimator of $g^{\star}$ is $\hat{g}=\hat{g}_{\hat M}$, and the estimator of $\f^{\star}$ is $\hat{\f}:=\hat{\f}_{\hat M}$ so that $\hat{g}=g^{\hat \Q,\hat\f}$.\\

The least squares estimator does not have an explicit form such as in usual nonparametric estimation, so that one has to use numerical minimization algorithms. As initial point of the minimization algorithm, we shall use the spectral estimator, see Section \ref{sec:simul} for more details. Since the spectral estimator is consistent, see \cite{YCES}, the algorithm does not suffer from initialization problems.

\section{Adaptive estimation of the emission distributions}
\label{subsec:main1}

\subsection{Oracle inequality for the estimation of $g^{\star}$}
We now fix a subset ${\mathcal F}$ of $\mathbf{L}^{2}(\Y,\L)$, and we shall use the following assumption:
\begin{enumerate}[{\bf [HF]}]
\item
\label{F:K=2} 
{\it ${\mathcal F}$ is a  closed subset  of $\mathbf{L}^{2}(\Y,\L)$ such that:
 for any $f\in {\mathcal F}$, $\int f d\L =1$, $\|f\|_{2}\leq C_{\mathcal F,2}$ and $\|f\|_{\infty}\leq C_{\mathcal F,\infty}$ for some fixed positive $C_{\mathcal F,2}$ and $C_{\mathcal F,\infty}$.
}
\end{enumerate}

Our first main result is an oracle inequality  for the estimation of $g^{\star}$ which is stated below and proved in Section \ref{sec:proof-oracle}.
We denote by $\mathfrak S_{K}$ the set of permutations of $\{1,\ldots,K\}$. When $a$ is a vector, $\|a\|_{2}$ denotes its Euclidian norm, and when $A$ is a matrix, $\|A\|_{F}$ denotes its Frobenius norm.

\begin{theorem}
\label{prop:oracle}
Assume {\bf [H1]-[H4]} and {\bf [HF]}. Assume also $\f^{\star}\in{\mathcal F}^{K}$, and for all $M$,  $\f^{\star}_{M}\in{\mathcal F}^{K}$. Then, there exists positive constants $N_{0}$, $\rho^{\star}$
 and $A_{1}^{\star}$ 
 $($depending on $C_{\mathcal F,2}$ and $C_{\mathcal F,\infty}$ $(${\bf Scenario B}$)$ or on $\Q^{\star}$, $C_{\mathcal F,2}$ and $C_{\mathcal F,\infty}$  $(${\bf Scenario A}$))$ such that, if 
\[
\pen (N,M) \geq \rho^{\star}\frac{M\log N}{N} 
\]
then for all $x>0$, for all $N\geq N_{0}$, one has with probability $1-(e-1)^{-1}e^{-x}$, for any permutation $\tau\in \mathfrak S_{K}$, 
\begin{eqnarray*}
\|\hat{g}-g^{\star}\|_{2}^2
&\leq& 6 \inf_{M}\left\{\|g^{\star}-g^{\Q^{\star},{\f}^{\star}_{M}}\|_{2}^2+\pen(N,M)\right\}+ A^{\star}_{1}\frac{x}{N}\\
&&+18 C_{\mathcal F,2}^6\big(2\|\Q^\star-\mathbb P_{\tau}\hat\Q_N\mathbb P_{\tau}^{\top}\|_F^2+\|\pi^\star-\mathbb P_{\tau}\hat\pi\|_{2}^2\big).
\end{eqnarray*}
Here, $\mathbb P_{\tau}$ is the permutation matrix associated to $\tau$.
\end{theorem}

The important fact in this oracle inequality is that the minimal possible penalty is of order $M/N$ (up to logarithmic terms) and not $M^{3}/N$ as is usually the case when estimating a joint density of three random variables, so that we get a minimax rate adaptive estimator of the joint density $g^{\star}$. 

\subsection{Main result}
\label{sec:MAIN}

The problem is now to deduce from Theorem \ref{prop:oracle} a result on $\|f^{\star}_{k}-\hat{f}_{k}\|_{2}^{2}$, $k=1,\ldots,K$. This is the cornerstone of our work: we prove that, under a  technical assumption on the parameters of the unknown HMM, a direct lower bound links $\|\hat{g}-g^{\star}\|_{2}^2$ to $\sum_{k=1}^{K}\|f^{\star}_{1}-\hat{f}_{k}\|_{2}^{2}$, up to some positive constant. Let us now describe the assumption and comment on its genericity.

 For any $\f\in{\mathcal F}^{K}$, define $G(\f)$ the $K\times K$ matrix with coefficients
  $G(\f)_{i,j}=\langle f_{i},f_{j}\rangle$, $i,j=1,\ldots,K$. Notice that under the assumption {\bf [H4]},  $G(\f^{\star})$ is positive definite. Let $\Q$ be a transition matrix verifying {\bf [H1]-[H2]} and let $A_{Q}$ be the diagonal matrix having the stationary distribution $\pi$ of $\Q$ on the diagonal.
We shall now define a quadratic form with coefficients depending on  $\Q$ and $G(\f)$.
 If $U$ is a $K\times K$ matrix such that $U{\mathbf{1}}_{K}=0$, 
\begin{align*}
\mathcal D:=
 \sum_{i,j=1}^{K}\Big\{ \big( &\Q^{T}A_{Q} UG(\f) U^{T} A_{Q}\Q\big)_{i,j}\big( G(\f) \big)_{i,j}\big( \Q G(\f)\Q^{T}\big)_{i,j}\\
&+\big( \Q^{T}A_{Q} G(\f) A_{Q}\Q\big)_{i,j}\big(  UG(\f)U^{T} \big)_{i,j}\big( \Q G(\f)\Q^{T}\big)_{i,j}\\
&+\big( \Q^{T}A_{Q}G(\f)A_{Q}\Q\big)_{i,j}\big( G(\f) \big)_{i,j}\big( \Q UG(\f)U^{T} \Q^{T}\big)_{i,j}\Big\}\\
+2\sum_{i,j}\Big\{ \big( &\Q^{T}A_{Q} UG(\f) A_{Q}\Q\big)_{i,j}\big( UG(\f) \big)_{j,i}\big( \Q G(\f)\Q^{T}\big)_{i,j}\\
&+\big( \Q^{T}A_{Q} UG(\f) A_{Q}\Q\big)_{i,j}\big(  \Q UG(\f)\Q^{T} \big)_{j,i}\big( G(\f)\big)_{i,j}\\
&+ \big( UG(\f) \big)_{i,j}\big( \Q UG(\f) \Q^{T}\big)_{j,i}\big( \Q^{T}A_{Q}G(\f)A_{Q}\Q\big)_{i,j}\Big\}
\end{align*}
defines a semidefinite positive quadratic form  ${\mathcal D}$ in the coefficients $U_{i,j}$, $i=1,\ldots,K$, $j=1,\ldots,K-1$. The determinant  of this quadratic form is a polynomial in the coefficients of the matrices $\Q$, $A_{Q}$ and  $G(\f)$. Since the coefficients of $A_{Q}$ are rational functions of the coefficients of the matrix $\Q$,  this determinant is also a rational function of the coefficients of the matrices $\Q$ and  $G(\f)$. Define $H(\Q,G(\f))$  the numerator of the determinant.  Then $H(\Q,G(\f))$ is a polynomial in  the coefficients of the matrices $\Q$ and  $G(\f)$. 
Our assumption will  be:

\begin{center}
{\bf [HD]}\;\;\;$H(\Q^{\star},G(\f^{\star}))\neq 0$. 
\end{center}
Since $H$ is a polynomial function of $Q_{i,j}^{\star}$,  $i=1,\ldots,K$, $j=1,\ldots,K-1$, and $\langle f^{\star}_{i},f^{\star}_{j}\rangle$, $i,j=1,\ldots,K$, the assumption {\bf [HD]} is generically satisfied. We have been able to prove that  {\bf [HD]} always holds in the case $K=2$. We were only able to prove this result by direct computation, it is given in Section \ref{sec:bruteforce}.

\begin{lemma}
\label{LEMMA:INEQ3}
Assume $K=2$. Then for all $\Q^{\star}$ and $\f^{\star}$ such that {\bf [H1]-[H4]} hold,  one has $H(\Q^{\star},G(\f^{\star})) >0$.
\end{lemma}

\noindent
Notice now that, when {\bf [HD]} and {\bf [H1]-[H3]} hold, it is possible to define a compact neighborhood $\mathcal V$ of $\Q^{\star}$ such that, for all $\Q\in{\mathcal V}$,  $H(\Q,G(\f^{\star}))\neq 0$,  {\bf [H1]-[H3]} hold  for $\Q$ and $T_{\Q}\subset T_{\Q^{\star}}$. 
\\
For any $\h\in  \left(\mathbf{L}^{2}(\Y,\L)\right)^{K}$, define $\|\h \|_{\Q}^{2}:=\min_{\tau\in T_{\Q}}\{\sum_{k=1}^{K}\|h_{k} +f^{\star}_{k}-f^{\star}_{\tau(k)}\|_{2}^{2}\}$.
 Denote $\|\h \|_{2}^{2}:=
 \{\sum_{k=1}^{K}\|h_{k}
\|_{2}^{2}\}$.
We may now state the theorem which is the cornerstone of our main result.
 \begin{theorem}
\label{LEMMA:INEQ2}
Assume  {\bf [H1]-[H4]} and {\bf [HD]}.
Let $\mathcal K$ be a closed bounded  subset of  $\left(\mathbf{L}^{2}(\Y,\L)\right)^{K}$ such that if  $\h\in {\mathcal K}$, then $\int h_{i}d\L =0$, $i=1,\ldots,K$.  Let $\mathcal V$ be a compact neighborhood of $\Q^{\star}$ such that, for all $\Q\in{\mathcal V}$, $H(\Q,G(\f^{\star}))\neq 0$, {\bf [H1]-[H3]} holds for $\Q$ and $T_{\Q}\subset T_{\Q^{\star}}$. 
Then there exists a positive constant $c({\mathcal K}, {\mathcal V}, \F^{\star})$ such that
\[
\forall {\h}\in {\mathcal K}, \;\forall \Q\in{\mathcal V},\quad
\|g^{\Q,\f^{\star}+\h} -g^{\Q,\f^{\star}}\|_{2} \geq c({\mathcal K}, {\mathcal V}, \F^{\star}) \|\h \|_{\Q^{\star}}.
\]
\end{theorem}

\noindent
This theorem is proved in Section \ref{sec:proof-ineq2}.
\\

We are now ready to prove our main result on the penalized least squares estimator  of the emission densities. 
The following theorem gives an oracle inequality for the estimators of the emission distributions  provided the penalty is adequately chosen.
It is proved in Section~\ref{sec:proof-theoadapt}.

 \begin{theorem}[Adaptive estimation]
\label{theo:adapt}
Assume  {\bf [H1]-[H4]},   {\bf [HF]} and {\bf [HD]}. Assume also  that
for all $M$,  $\f^{\star}_{M}\in{\mathcal F}^{K}$.  
Let $\mathcal V$ be a compact neighborhood of $\Q^{\star}$ such that, for all $\Q\in{\mathcal V}$, $H(\Q,G(\f^{\star}))\neq 0$ and {\bf [H1]-[H3]} holds for $\Q$. 
Then, there exists a positive constant 
 $A^{\star}$
 $($depending on $\mathcal V$, $\f^{\star}$, $C_{\mathcal F,2}$ and $C_{\mathcal F,\infty})$ and positive constants $N_{0}$ and  $\rho^{\star}$
 $($depending on $C_{\mathcal F,2}$ and $C_{\mathcal F,\infty}$ $(${\bf Scenario A}$)$ or on $\Q^{\star}$, $C_{\mathcal F,2}$ and $C_{\mathcal F,\infty}$  $(${\bf Scenario B}$))$ such that, if 
\[
\pen (N,M) \geq \rho^{\star}\frac{M\log N}{N} 
\]
then for all $x>0$, for all $N\geq N_{0}$, for any permutation $\tau_{N}\in\mathfrak S_{K}$, 
with probability larger than $1-(e-1)^{-1}e^{-x}-\P\left(  \mathbb P_{\tau_{N}}\hat\Q\mathbb P_{\tau_{N}}^{T} \notin \mathcal V \right)$, there exists $\tau\in T_{\Q^{\star}}$ such that
\begin{align*}
\sum_{k=1}^{K}\|f^{\star}_{\tau(k)}-\hat{f}_{\tau_{N}(k)}\|_{2}^{2}\leq 
 A^{\star}&\Big[\inf_{M}\left\{\sum_{k=1}^{K}\|f^{\star}_{k}- f^{\star}_{M,k}\|_2^2+\pen(N,M)\right\}\\
& +\lVert\Q^{\star}-\mathbb P_{\tau_{N}}\hat\Q \mathbb P_{\tau_{N}}^{\top}\lVert_{F}^{2}
 +\lVert\pi^{\star}-\mathbb P_{\tau_{N}}\hat\pi\lVert_{2}^{2} +\frac{x}{N}\Big].
\end{align*}
\end{theorem}

\begin{remark}
As usual in HMM or mixture model, it is only possible to estimate the model up to label switching of the hidden states, this is the meaning of the permutation $\tau_{N}$.
\end{remark}

\begin{remark}
An important consequence of the theorem is that  a right choice of the penalty leads to a rate minimax adaptive estimator up to a $\log N$ term, see Corollary \ref{cor:adapt} below. For this purpose, one has to choose an estimator $\hat\Q$ of $\Q^{\star}$ which is, up to label switching, consistent with controlled rate. One possible choice is a spectral estimator.
\end{remark}


To apply Theorem \ref{theo:adapt} one has to choose an estimator $\hat\Q$ with controlled behavior, to be able to evaluate the probability of the event $\{\mathbb P_{\tau_{N}}\hat\Q \mathbb P_{\tau_{N}}\in\mathcal V\}$ and  the rate of convergence of $\mathbb P_{\tau_{N}}\hat\Q \mathbb P_{\tau_{N}}$ and $\mathbb P_{\tau_{N}}\hat\pi$. One possibility is to use the spectral estimator described in Section \ref{sec:spectral}.
To get the following result (proved in Section~\ref{sec:proof-coradapt}), we propose to use the spectral estimator with, for each $N$, the dimension $M_{N}$ chosen such that  $\eta_{3}(\Phi_{M_{N}})= O\big((\log N)^{1/4}\big)$, see Section~\ref{sec:spectral} for a definition of $\eta_{3}$. 
\begin{corollary}
\label{cor:adapt}
With this choice of $\hat\Q$, under the assumptions of Theorem \ref{theo:adapt} , there exists a sequence of permutations 
 $\tau_{N}\in\mathfrak S_{K}$ such that as $N$ tends to infinity, 
$$
\E\left[\sum_{k=1}^{K}\|f^{\star}_{k}-\hat{f}_{\tau_{N}(k)}\|_{2}^{2}\right]= 
O\left(\inf_{M'}\left\{\sum_{k=1}^{K}\|f^{\star}_{k}- f^{\star}_{M',k}\|_{2}^{2}+\pen(N,M')\right\}+\frac{\log N}{N}\right).
$$
\end{corollary}

Thus, choosing $\pen(N,M)= \rho M\log N /N$ for a large $\rho$ leads to the minimax asymptotic rate of convergence up to a power of $\log N$.
Indeed, 
standard results in approximation theory (see \cite{devore1993constructive} for instance) show that one can upper bound the approximation error $\lVert f^{\star}_{k}- f^{\star}_{M,k}\lVert_{2}$ by $\mathcal O(M^{-\frac sD})$ where $s>0$ denotes a regularity parameter. Then the trade-off is obtained for $M^{\frac 1D}\sim (N/\log N)^{\frac1{2s+D}}$, which leads to the quasi-optimal rate $(N/\log N)^{-\frac s{2s+D}}$ for the nonparametric estimation when the minimal smoothness of  the emission densities  is $s$.
Notice that the algorithm automatically selects the best $M$ leading to this rate. 

To implement the estimator, it remains to choose a value for $\rho$ in the penalty.  The calibration of this parameter is a classical issue and could be the subject of a full paper. In practice one can use the slope heuristic as in \cite{BMM12}.

\section{Nonparametric spectral method}
\label{sec:spectral}

This section is devoted to a short description of the nonparametric spectral method for sake of completeness: we describe the algorithm, and give the results we need to support the use of spectral estimators to initialize our algorithm. A detailed study of the nonparametric spectral method is given in \cite{YCES}.

The following procedure (see Algorithm \ref{alg:Spectral}) describes a tractable approach to estimate the transition matrix in a way that can be used for the penalized least squares estimator of the emission densities, and also for the estimation of the projections of the emission densities that may be used to initialize the least squares algorithm. The procedure is based on recent developments 
in parametric estimation of HMMs. For each fixed $M$, we estimate the projection of the emission distributions on the basis $\Phi_{M}$ using the spectral method proposed  in \cite{anandkumar2012method}. As the authors of the latter paper explain, this allows further to estimate the transition matrix (we use a modified version of their estimator), and we set the estimator of the stationary distribution as the stationary distribution of the estimator of the transition matrix. The computation of those estimators is particularly simple: it is based on one SVD, some matrix inversions and one diagonalization. One can prove, with overwhelming probability, all matrix inversions and the diagonalization can be done rightfully, see \cite{YCES}.
In the following, when $A$ is a $(p\times q)$ matrix with $p\geq q$, $A^{\top}$ denotes the transpose matrix of $A$, $A(k,l)$ its $(k,l)$th entry, $A(\ldotp,l)$ its $l$th column and $A(k,\ldotp)$ its $k$th line. When $v$ is a vector of size $p$, we denote by $\D v$ the diagonal matrix with diagonal entries $v_{i}$ and, by abuse of notation, $\D v=\D{v^{\top}}$.

{
\begin{algorithm}[t]
  \SetAlgoLined
  \KwData{An observed chain $(Y_{1},\ldots,Y_{N})$ and a number of hidden states $K$.}
  \KwResult{Spectral estimators $\hat\pi$, $\hat\Q$ and $(\hat f_{M,k})_{k\in\X}$.}
    \BlankLine
\begin{enumerate}[{\bf [Step 1]}]
\item Consider the following empirical estimators: For any $a,b,c$ in $\{1,\ldots,M\}$, $\hatl(a):=\frac1N\sum_{s=1}^{N}\varphi_{a}(Y_{1}^{(s)})$, $\hatM(a,b,c):=\frac1N\sum_{s=1}^{N}\varphi_{a}(Y_{1}^{(s)})\varphi_{b}(Y_{2}^{(s)})\varphi_{c}(Y_{3}^{(s)})$, $\hatN(a,b):=\frac1N\sum_{s=1}^{N}\varphi_{a}(Y_{1}^{(s)})\varphi_{b}( Y_{2}^{(s)})$, $\hatPj(a,c):=\frac1N\sum_{s=1}^{N}\varphi_{a}(Y_{1}^{(s)})\varphi_{c}(Y_{3}^{(s)})$.\item
Let $\hat\U$ be the $M\times K$ matrix of orthonormal right singular vectors of $\hatPj$ corresponding to its top $K$ singular values. 
\item 
Form the matrices for all $b\in\{1,\ldots,M\}$, $\hat\B(b):=(\hat\U^{\top}\hatPj\hat\U)^{-1}\hat\U^{\top}\hatM(\ldotp,b,\ldotp)\hat\U$.
\item
Set $\Theta$ a $(K\times K)$ random unitary matrix uniformly drawn and form the matrices for all $k\in\{1,\ldots,K\}$, ${\hat\C}(k):=\sum_{b=1}^{M}(\hat\U\Theta)(b,k)\hat\B(b)$.
\item
Compute $\hat\Rbf$ a $(K\times K)$ unit Euclidean norm columns matrix that diagonalizes the matrix $\hat\C(1)$: $\hat\Rbf^{-1}\hat\C(1)\hat\Rbf=\D{(\hat\Lambda(1,1),\ldots,\hat\Lambda(1,K))}$.
\item
Set for all $k,k'\in\X$, $\hat\Lambda(k,k'):=(\hat\Rbf^{-1}\hat\C(k)\hat\Rbf)(k',k')$ and $\hatO:=\hat\U\Theta\hat\Lambda$.
\item
Consider the emission laws estimator $\tilde\f:=(\tilde f_{M,k})_{k\in\X}$ defined by for all $k\in\X$, $\hat f_{M,k}:=\sum_{m=1}^{M}\hatO(m,k)\varphi_{{m}}$.
\item 
Set $\tilde\pi:=\big(\hat\U^{\top}\hatO\big)^{-1}\hat\U^{\top}\hatl$.
\item
Consider the transition matrix estimator: 
\[
\hat\Q:=\Pi_{\mathrm{TM}}\Big(\big(\hat\U^{\top}\hatO\D{\tilde\pi}\big)^{-1}\hat\U^{\top}\hatN\hat\U\big(\hatO^{\top}\hat\U\big)^{-1}\Big)\,,
\] where $\Pi_{\mathrm{TM}}$ denotes the projection onto the convex set of transition matrices, and define
$\hat\pi$ as the stationary distribution of $\hat\Q$.
\end{enumerate}
  \caption{Nonparametric spectral estimation of HMMs}
\label{alg:Spectral}
\end{algorithm}}

We now  state a  result which allows to derive the asymptotic properties of the spectral estimators. 
Let us define:
\[
\eta_{3}^{2}(\Phi_{M}):=\sup_{y,y'\in\Y^{3}}\sum_{a,b,c=1}^{M}(\varphi_{a}(y_{1})\varphi_{b}(y_{2})\varphi_{c}(y_{3})-\varphi_{a}(y'_{1})\varphi_{b}(y'_{2})\varphi_{c}(y'_{3}))^{2}.
\]
Note that in the examples {\bf (Spline)}, {\bf (Trig.)} and {\bf (Wav.)} we have:
\[
\eta_{3}(\Phi_{M})\leq C_{\eta}M^{\frac{3}{2}}
\]
where $C_{\eta}>0$ is a constant. The following theorem is proved in \cite{YCES}. Its statement concerns {\bf (Scenario B)} (same chain sampling) and the interested reader may consult \cite{YCES} for its statement under {\bf (Scenario A)}.
\begin{theorem}[Spectral estimators]
\label{thm:spectral}
Assume that {\bf [H1]-[H4]} hold. Then, there exist positive constant numbers $M_{\F^{\star}}$, $x(\Q^{\star})$, ${\mathcal C}(\Q^{\star},\F^{\star})$ and ${\mathbf N}(\Q^{\star},\F^{\star})$ such that the following holds.
For any $x\geq x(\Q^{\star})$, 
for any $\delta \in (0,1)$, for any $M\geq M_{\F^{\star}}$, there exists a permutation $\tau_{M}\in\mathfrak S_{K}$ such that the spectral method estimators $\hat f_{M,k}$, $\hat\pi$ and $\hat\Q$  satisfy: For any $N\geq {\mathbf N}(\Q^{\star},\F^{\star})\eta_3(\Phi_{M})^{2}x(-\log \delta )/\delta^{2}$, with probability greater than $1-2\delta-4 e^{-x}$, 
\begin{align*}
\lVert f^{\star}_{M,k}-\hat f_{M,\tau_{M}(k)}\lVert_{2}&\leq {\mathcal C}(\Q^{\star},\F^{\star})\frac{\sqrt{-\log \delta}}{\delta}\frac{\eta_3(\Phi_{M})}{\sqrt N}\sqrt{x}\,,\\
\lVert\pi^{\star}-\mathbb P_{\tau_{M}}\hat\pi\lVert_{2}
&\leq
 {\mathcal C}(\Q^{\star},\F^{\star})\frac{\sqrt{-\log \delta}}{\delta}\frac{\eta_3(\Phi_{M})}{\sqrt N}\sqrt{x}\,,\\
 \lVert\Q^{\star}-\mathbb P_{\tau_{M}}\hat\Q \mathbb P_{\tau_{M}}^{\top}\lVert
&\leq
 {\mathcal C}(\Q^{\star},\F^{\star})\frac{\sqrt{-\log \delta}}{\delta}\frac{\eta_3(\Phi_{M})}{\sqrt N}\sqrt{x}\,.
\end{align*}
\end{theorem}

\section{Numerical experiments}
\label{sec:simul}

\subsection{General description}
In this section we present the numerical performances of our method. We recall that the experimenter knows nothing about the underlying hidden Markov model but the number of hidden states $K$. In this set of experiments, we consider the regular histogram basis or the trigonometric basis for estimating emission laws given by beta laws from a single chain observation of length $N=5e4$. 

Our procedure is based on the computation of the empirical least squares estimators $\hat g_M$ defined as minimizers of the empirical contrast $\gamma_{N}$ on the space $\dS ({\hat \Q},{M})$ where $\hat \Q$ is an estimator of the transition matrix (for instance the spectral estimation of the transition matrix). Since the function $\gamma_{N}$ is non-convex, we use a second order approach estimating a positive definite matrix (using a covariance matrix) within an iterative procedure called CMAES for Covariance Matrix Adaptation Evolution Strategy, see \cite{hansen2006cma}. Using this latter algorithm, we search for the minimum of $\gamma_{N}$ with starting point the spectral estimation of the emission laws. 

Then, we estimate the size of the model thanks to
\begin{equation}
\label{eq:M_rho}
\hat M(\rho)\in\displaystyle\arg\min_{M=1,\ldots,M_{\max}}\left\{\gamma_{N}(\hat g_M)+\rho\,\frac{M\log N}{N}\right\}\,,
\end{equation}
where the penalty term $\rho$ has to be tuned and the maximum size of the model $M_{\max}$ can be set by the experimenter in a data-driven procedure. 

Indeed, we shall apply the slope heuristic to adjust the penalty term and to choose $M_{\max}$. As presented in \cite{BMM12}, the minimum contrast function $M\mapsto\gamma_N(\hat g_{M})$ should have a linear behavior for large values of $M$. The experimenter has to consider $M_{\max}$ large enough in order to observe this linear stabilization, as depicted in Figure~\ref{fig:Slope_Heuristic}. The slope of the linear interpolation is then $(\hat\rho/2)\log N/N$ (recall that the sample size $N$ is fixed here) where $\hat\rho$ is the slope heuristic choice on how $\rho$ should be tuned. Another procedure (theoretically equivalent) consists in plotting the function $\rho\mapsto\hat M(\rho)$ which is a non-increasing piecewise constant function. The estimated $\hat\rho$ is such that the largest drop (called ``dimension jump'') of this function occurs at point $\hat\rho/2$. We illustrate this procedure in Figure~\ref{fig:Slope_Heuristic_Jump} where one can clearly point the jump and deduce the size $\hat M$.

To summarize, our procedure reads as follows.
\begin{enumerate}
\item For all $M\leq M_{\max}$, compute the spectral estimations $(\hat \Q,\hat \pi)$ of the transition matrix and its stationary distribution and the spectral estimation $\tilde \f$ of the emission laws. This is straightforward using the procedure described by {\bf [Step1-9]} in Section \ref{sec:spectral}.
\item For all $M\leq M_{\max}$, compute a minimum $\hat g_M$ of the empirical contrast function $\gamma_N$ using ``Covariance Matrix Adaptation Evolution Strategy'', see \cite{hansen2006cma}. 
Use the estimation $\tilde \f$ of the spectral method as a starting point of CMAES.
\item Tune the penalty term using the slope heuristic procedure and select $\hat M$.
\item Return the emission laws of the solution of point $(2)$ for $M=\hat M$.
\end{enumerate}
Note that the size $M$ of the projection space for the spectral estimator has been set as the one chosen by the slope heuristic for the empirical least squares estimators. 

All the codes of the numerical experiments are available at \url{https://mycore.core-cloud.net/public.php?service=files&t=44459ccb178a3240cfb8712f27a28d75}. We shall indicate that the slope heuristic has been done using \href{http://www.math.univ-toulouse.fr/~maugis/CAPUSHE.html}{CAPUSHE}, the Matlab graphical user interface presented in \cite{BMM12}.

\subsection{Complexity}
\label{sec:Complexity}
A crucial step of our method lies in computing the empirical least squares estimators $\hat g_M$. One may struggle to compute $\hat g_M$ since the function $\gamma_{N}$ is non-convex. It follows that an acceptable procedure must start from a good approximation of $\hat g_M$. This is done by the spectral method. Observe that the key leitmotiv throughout this paper is a two steps estimation procedure that starts by the spectral estimator. This latter has rate of convergence of the order of $N^{-s/(2s+3)}$ and seems to be a good candidate to initialize an iterative scheme that will converge towards $\hat g_M$. It follows that the main consuming operations in our algorithm are the following steps.
\begin{itemize}
\item  The computation of the tensor $\hatM$ of the empirical law of three consecutive observations  where we use three loops of size $M$ and one loop of size $N$ so the complexity is $\mathcal O(NM^3)$,
\item  The singular value decomposition of $\hatPj$ in the spectral method (complexity: $\mathcal O(M^3)$),
\item The computation of the minimum of the empirical contrast function: cost of one evaluation of the empirical contrast function $\mathcal O(K^3M^3)=\mathcal O(M^3)$ times the number $f(M,K)$ of evaluations while minimizing the empirical contrast. Recall that we start from the spectral estimator solution to get the minimum so a constant number of evaluation is enough in practice, say \texttt{stopeval} =$1e4$ using CMAES.
\end{itemize}
We have to compute the minimal contrast value for all models of size $M=1,\ldots, M_{\max}$ where $M_{\max}$ has to be chosen so that one can apply the slope heuristic. We deduce that the overall complexity of our algorithm is $\mathcal O\big((f(M_{\max},K)K^3\vee N)M_{\max}^4\big)$ 
where $f(M_{\max},K)$ is the number of evaluations of $\gamma_N$ while minimizing the empirical contrast. Since we use the spectral estimator as a starting point of the minimization of the empirical contrast, we believe that $f(M_{\max},K)$ can be considered as constant, say $1e4$. Note that the upper bound $M_{\max}$ has to be large enough in order to observe a linear stabilization of $M\mapsto\hat g_M$, see \cite{BMM12} for instance. Moreover, recall that the trade-off between the approximation bias and the penalty term (accounting for the standard error of the empirical law) is obtained for $M\sim (N/\log N)^{\frac D{2s+D}}$ where $s>0$ denotes the minimal smoothness parameter of the emission laws. In order to properly apply the slope heuristic, it is enough to consider models with this order of magnitude, so that $M_{\max}=\mathcal O((N/\log N)^{\frac D{2s+D}})$. It follows that the overall complexity of our procedure can be expressed in terms of the minimal smoothness parameter $s$ of the emission laws as
\[
\mathrm{Complexity}=\mathcal O\big(N^{1+\frac {4D}{2s+D}}\big)\,,
\]
as soon as $K=\mathcal O(N^{1/3})$ which is a reasonable assumption. Nevertheless, this theoretical bound is unknown for the practitioner since it involves the unknown minimal smoothness parameter $s>0$. For chains of length $\mathcal O(1e5)$, we have witnessed that one can afford a maximal model size $M_{\max}\leq 50$ and this allows to consider problems where typical sizes of $M$ ranges between $1$ and $50$. All numerical experiments of this paper fall in this frame.

\subsection{Comparison of the variances}

\begin{figure}[t]
\centering
\includegraphics[width=0.7\textwidth]{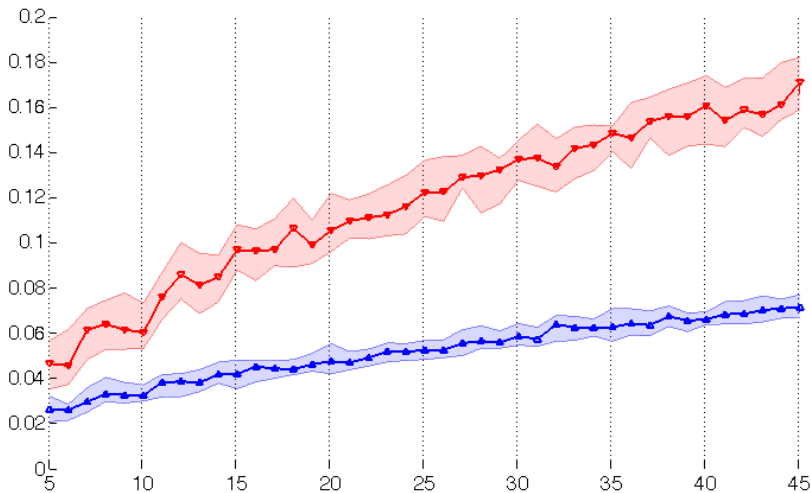}
\caption{Variance comparison of the spectral and empirical least squares estimators. The upper curve (in red) present the performance  (median value of the variance over 40 iterations) of the spectral method while the lower curve (in blue) the performance of the empirical least squares estimator. For each curve, we have plotted a shaded box plot representing the first and third quartiles.} 
\label{fig:Variance_HISTO}
\end{figure}

The quadratic loss can be expressed as a variance term and a bias term as follows
\[
\forall 1\leq k\leq K,\ \forall M\geq 0,\quad\|f^\star_k-\hat f_k\|_2^2=\|\hat f_k-f^\star_{M,k}\|_2^2+\|f^\star_k-f^\star_{M,k}\|_2^2
\]
where $f^\star_{M,k}$ is the orthogonal projection of $f_k^\star$ on $\Proj_{M}$ and $\hat f_k$ is any estimator such that $\hat f_k$ belongs to $\Proj_{M}$. Note that the bias term $\|f^\star_k-f^\star_{M,k}\|_2$ does not depend on the estimator $\hat f_k$. Hence, the variance term 
\[
\mathrm{Variance}_M(\hat f):=\min_{\tau\in\mathfrak S_{K}}\max_{1\leq k\leq K}\|\hat f_k-f^\star_{M,\tau(k)}\|_2^2\,,
\]
 accounts for the performances of the estimator $\hat f_k$.

As depicted in Figure~\ref{fig:Variance_HISTO}, we have compared, for each $M$, the variance terms obtained by the spectral method and the empirical least squares method over $40$ iterations on chains of length $N=5e4$. We have considered $K=2$ hidden states whose emission variables are distributed with respect to beta laws of parameters $(2,5)$ and $(4,2)$. This numerical experiment consolidates the idea that the least squares method significantly improves upon the spectral method. Indeed, even for small values of $M$, one may see in Figure \ref{fig:Variance_HISTO} that the variance term is divided by a constant factor. 

\subsection{Histogram basis and trigonometric basis as approximation spaces}

\begin{figure}[t]
\centering
\includegraphics[width=0.47\textwidth]{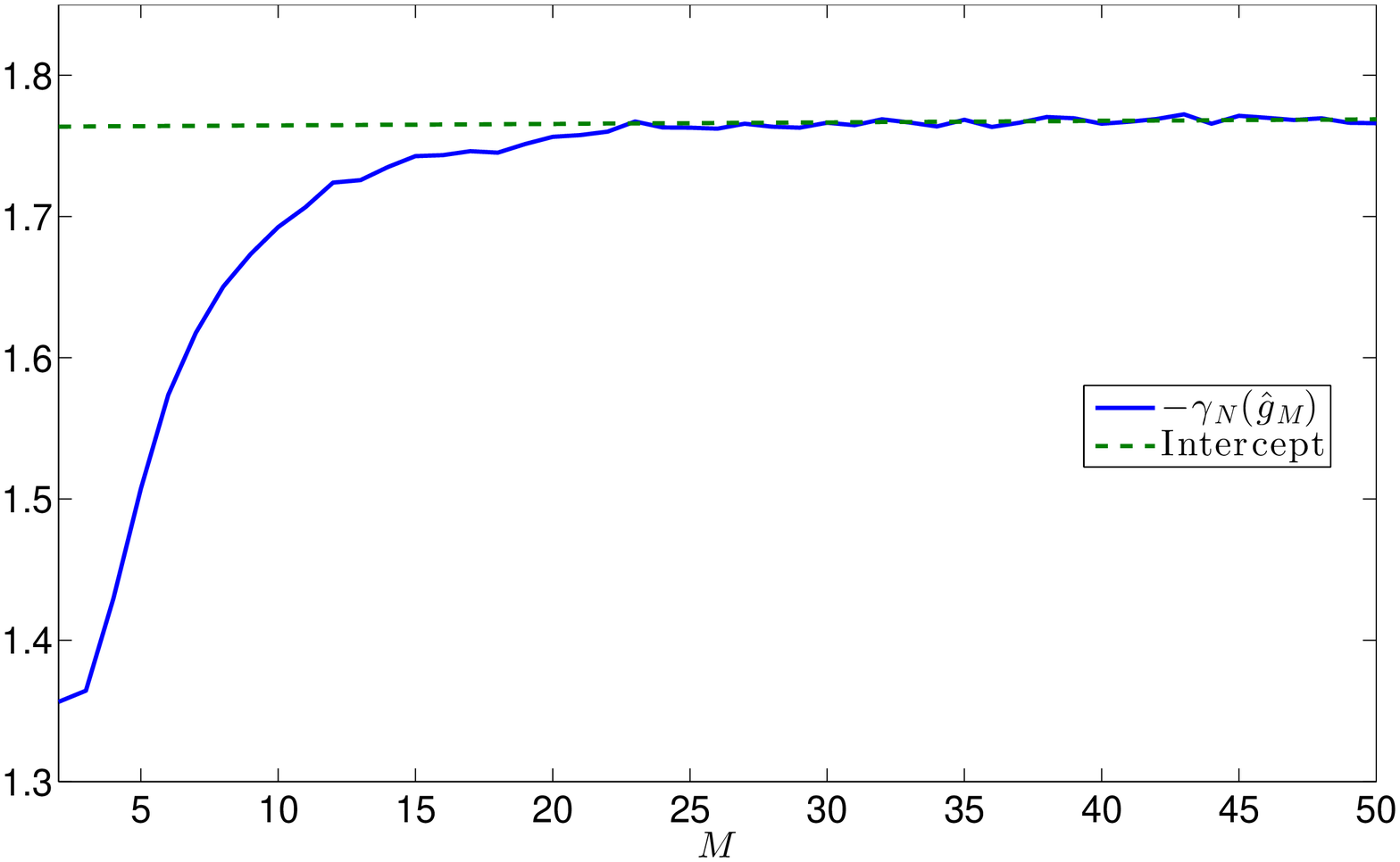}
\ \ \ 
\includegraphics[width=0.47\textwidth]{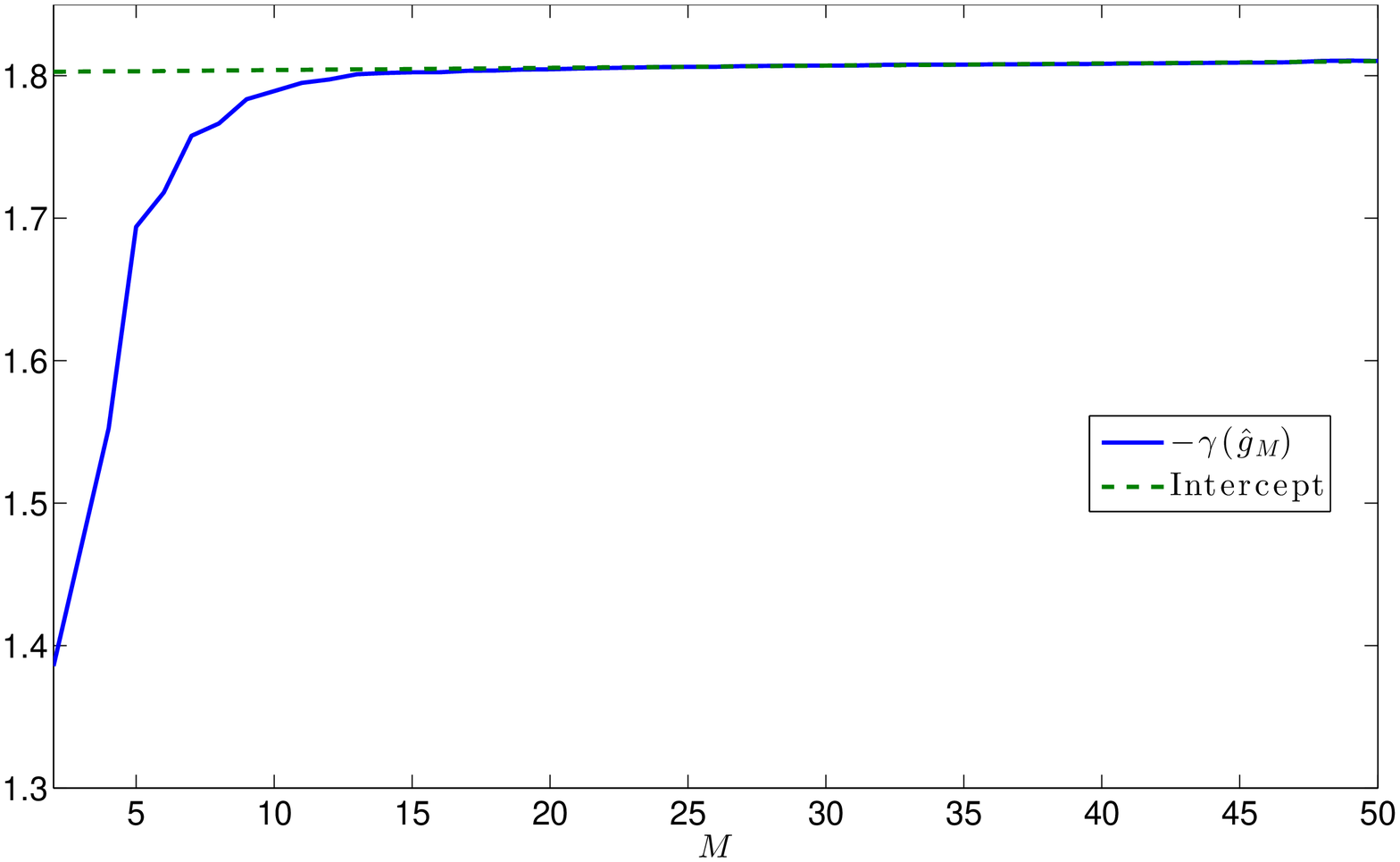}
\caption{Slope heuristic to choose $M$: the experimenter may observe a linear stabilization of the empirical contrast $\gamma_N$ for estimating beta emission laws of parameters $(2,5)$ and $(4,2)$. We have $K=2$ hidden states and $N=5e4$ samples along a single chain. On the left panel we have used the trigonometric basis as approximation space, the stabilization occurs on the points $M=30$ to $M=50$ and the interpolation of the slope leads to $\hat M=23$. On the right panel we have considered the trigonometric basis, the stabilization occurs on the points $M=20$ to $M=50$ and it leads to $\hat M=21$.}\label{fig:Slope_Heuristic}
\end{figure}

\begin{figure}[!t]
\centering
\includegraphics[height=5cm]{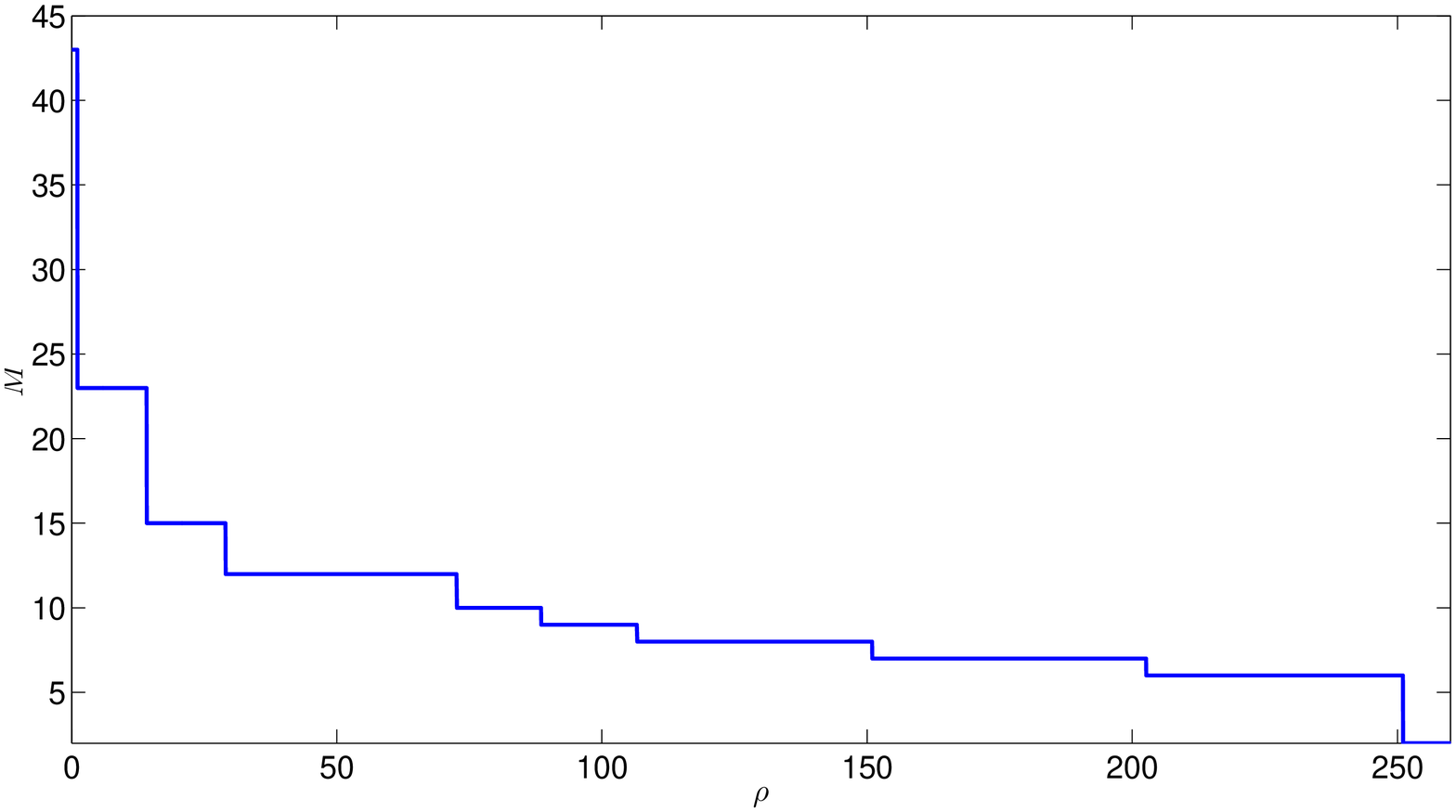}
\quad
\includegraphics[height=5cm]{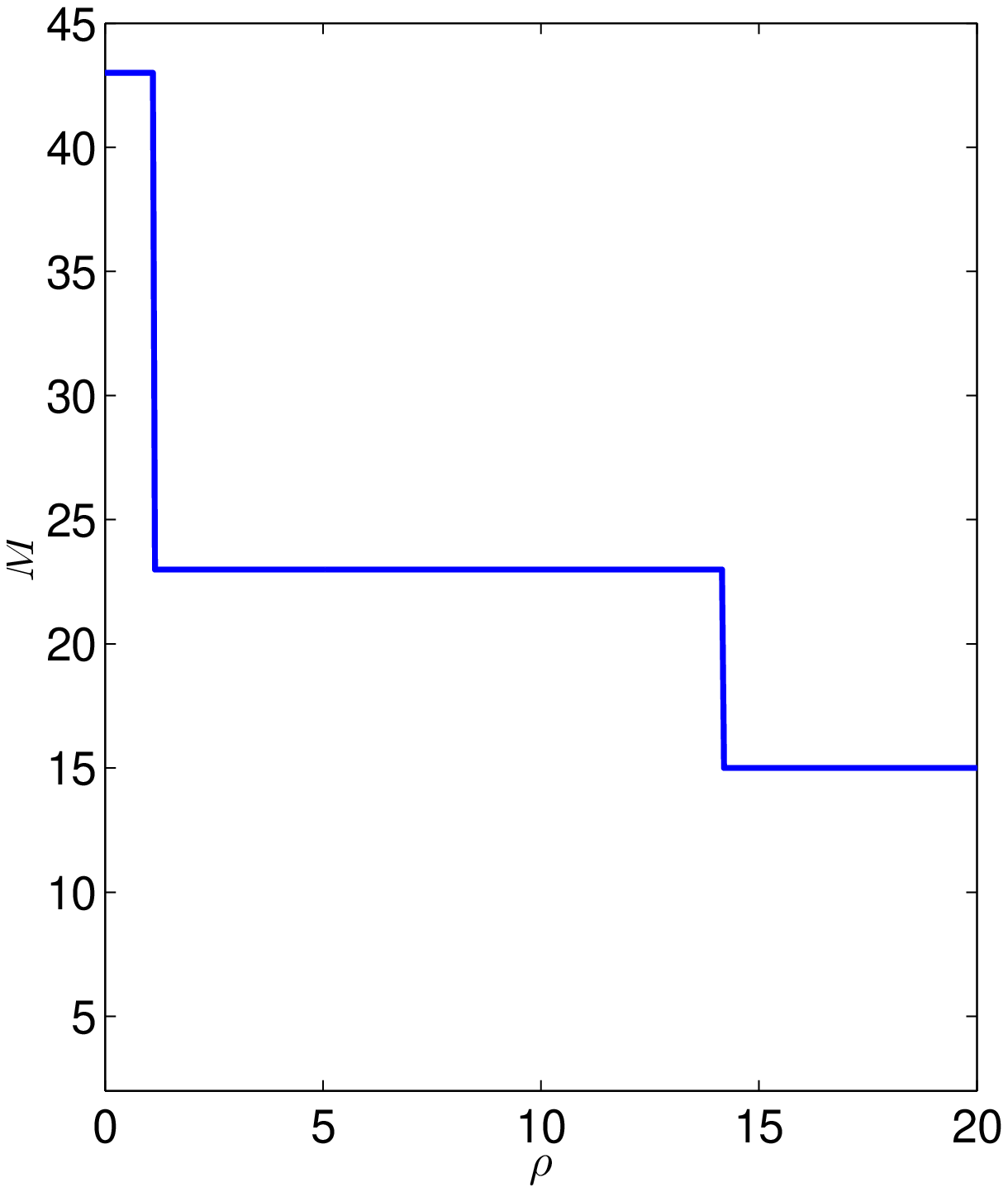}
\caption{Slope heuristic to choose $M$: the experimenter observes the largest drop of the function $\rho\mapsto \hat M(\rho)$ at $1.1$ so that $\hat\rho=2.2$ and $\hat M=23$. We have $K=2$ hidden states and a single chain of length $N=5e4$. We have used the histogram basis as approximation space.}\label{fig:Slope_Heuristic_Jump}
\end{figure}

An illustrative example of our method can be given using the histogram basis (regular basis with $M$ bins) or the trigonometric basis.  In the following experiments, we have $K=2$ hidden states and emission laws given by beta laws of parameters $(2,5)$ and $(4,2)$. Recall we observe a single chain of length $N=5e4$. 

We begin with the computation of the minimum contrast function $M\mapsto \gamma(\hat g_{M})$, as depicted in Figure \ref{fig:Slope_Heuristic}. Observe that the slope of this function unquestionably stabilizes at a critical value refer to as $\hat\rho/2$ in both the histogram and the trigonometric case. This leads to an adaptive choice of $\hat M=23$ for the histogram basis and $\hat M=21$ for the trigonometric basis, see Figures \ref{fig:Slope_Heuristic} and \ref{fig:Slope_Heuristic_Jump}.

\begin{figure}[!t]
\centering
\includegraphics[width=0.8\textwidth]{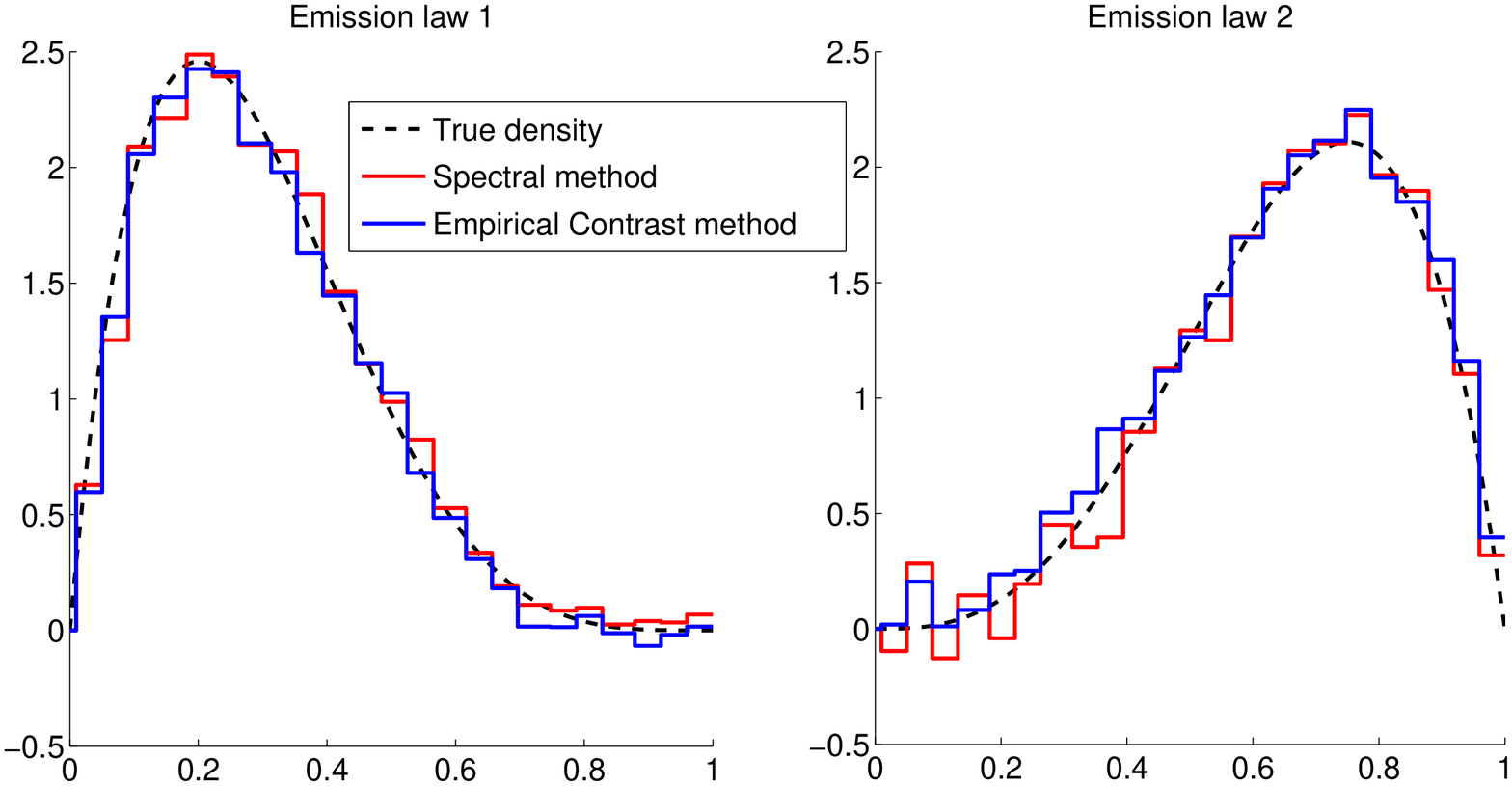}
\includegraphics[width=0.8\textwidth]{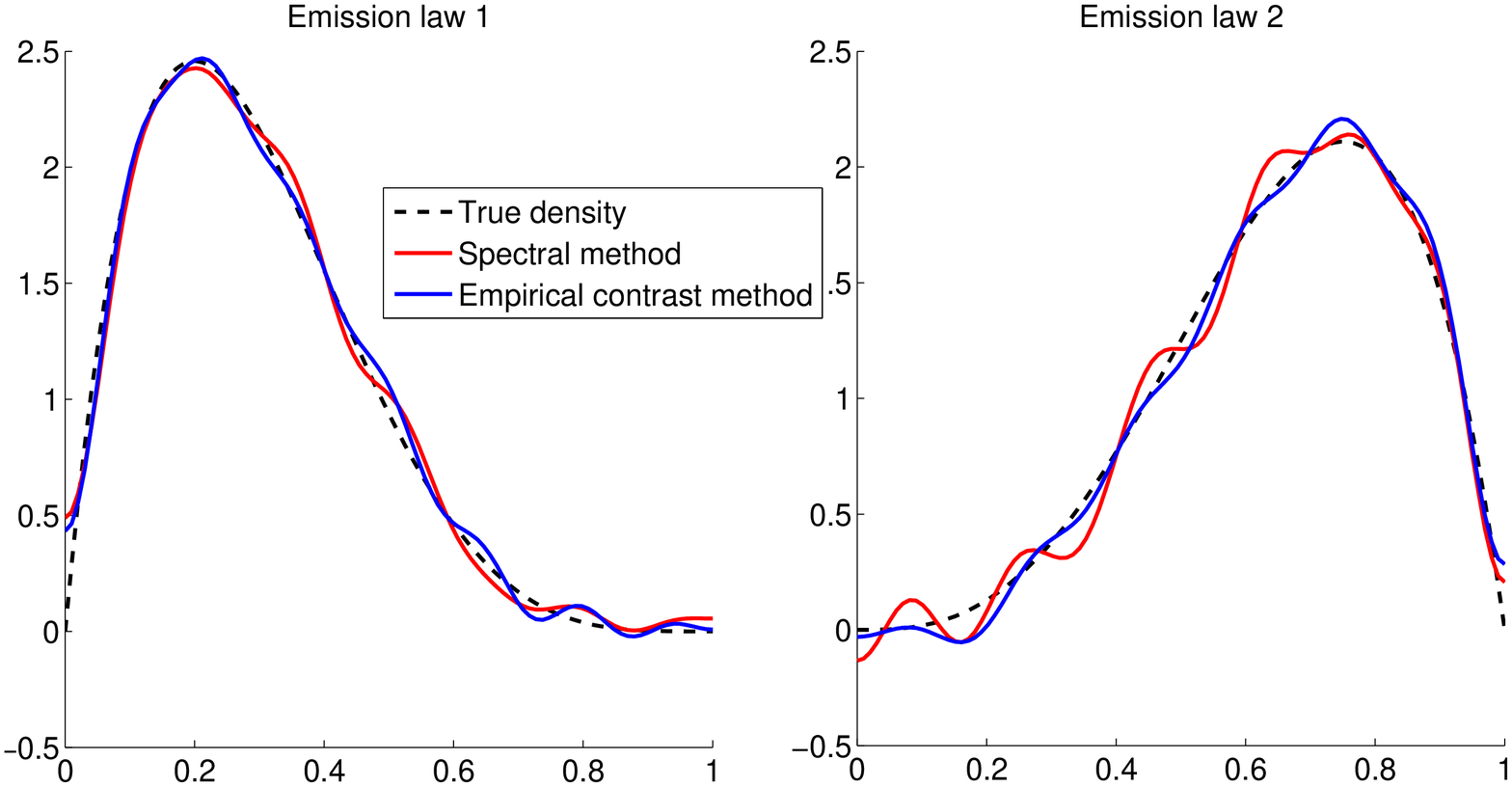}
\caption{Estimators of the emissions densities (beta laws of parameters $(2,5)$ and $(4,2)$) from the observation of a single chain of length $N=5e4$. On the top panels, we have used the histogram basis ($\hat M=23$). On the bottom panels, we have considered the trigonometric basis ($\hat M=21$).
}
\label{fig:Final_Result}
\end{figure}

Furthermore, one can see on Figure \ref{fig:Final_Result} that our method also qualitatively improves upon the spectral method in both the histogram and the trigonometric case. 
\subsection{Three states}

\begin{figure}[!t]
\centering
\includegraphics[width=0.9\textwidth]{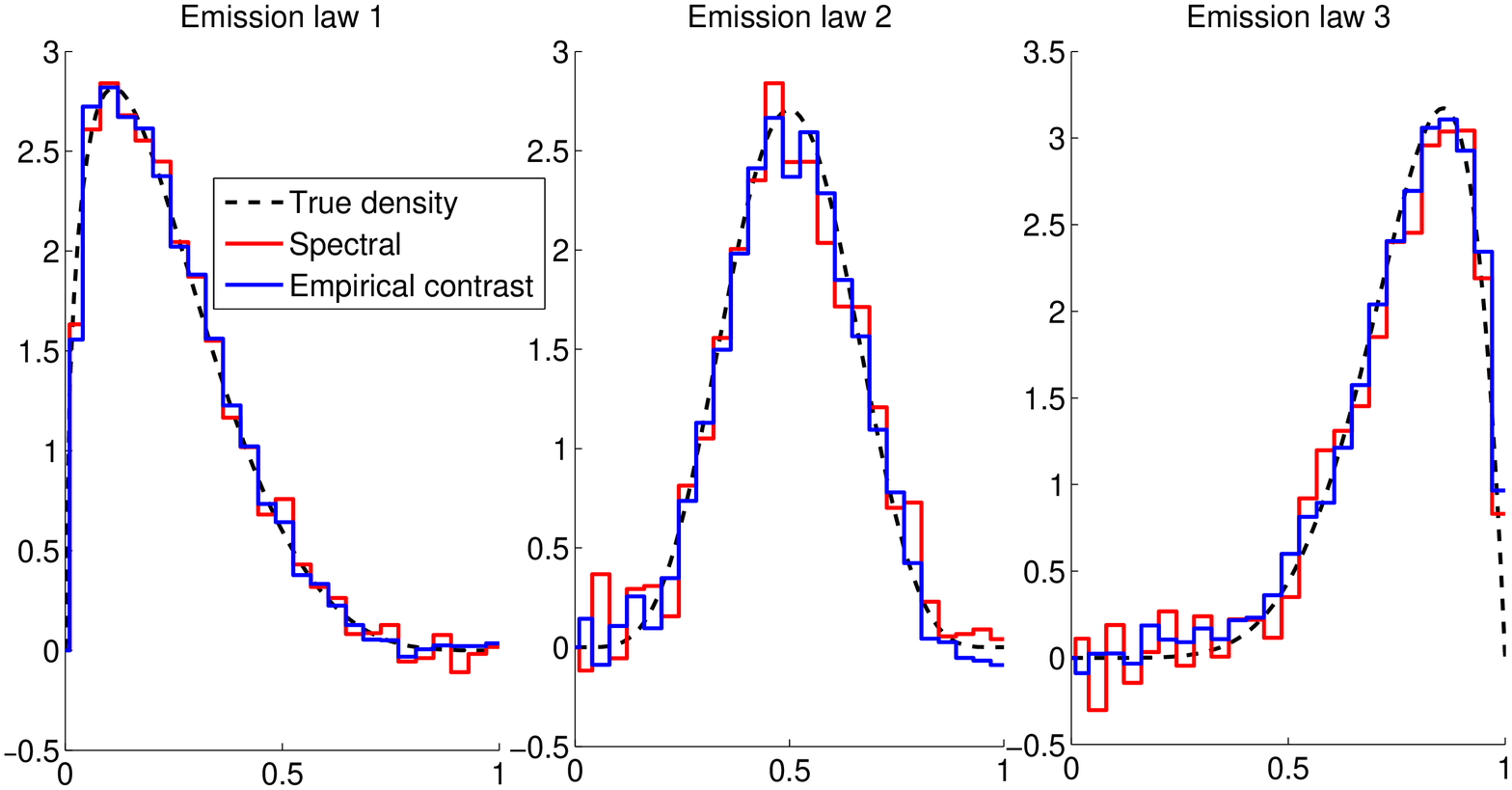}
\caption{Estimation of three densities given by beta laws of parameters $(1.5,5)$, $(6,6)$ and $(7,2)$ from a single chain of length $N=5e4$. We have used the histogram basis and we have found $\hat M=25$ using the slope heuristic.} 
\label{fig:ThreeStates}
\end{figure}

Our method can be performed for $K>2$ as illustrated in Figure \ref{fig:ThreeStates}. In this example $K=3$, the sample size is $N=5e4$ and the emission laws are three beta distributions with parameters $(1.5,5)$, $(6,6)$ and $(7,2)$. Note that the number of hidden states $K$ does not really impact on the complexity of the algorithm as we have seen in Section \ref{sec:Complexity}.

In this example, we were able to observe a linear stabilization of the minimum contrast function. The slope heuristic procedure led to an adaptive choice $\hat M=25$.

\section{Discussion}
\label{sec:discussion}
We have proposed a penalized least squares method to estimate the emission densities of the hidden chain when the transition matrix of the hidden chain is full rank and the emission probability distributions are linearly independent.
The algorithm may be initialized using spectral estimators.
The obtained estimators are adaptive rate optimal up to a $\log$ factor, where adaptivity is upon the family of emission densities.
The results hold under an assumption on the parameter that holds generically.  We have proved  that this assumption is always verified when there are two hidden states. We did not find a general argument to prove that the assumption always holds when $K>2$, 
and a natural question is to ask if, when the number of hidden states is $K>2$, this assumption is also  always verified. 
\\
 It is proved in \cite{alexandrovich2014nonparametric} that identifiability holds as soon as $f_{1}^{\star},\ldots,f_{K}^{\star}$ are distinct densities. The identifiability is obtained in that case using the marginal distribution  of dimension $2K+1$, that is the marginal distribution of $Y_{1},\ldots,Y_{2K+1}$. Thus, to get consistent estimators, one needs to use the joint distribution of $2K+1$ consecutive observations. Though linear independence is generically satisfied, one may wonder what happens when emission densities are not far to be linearly dependent. Simulations in  
\cite{Luc1} show that estimation becomes harder. 
In those practical situations where estimation becomes difficult, it is observed that the Gram matrix of $f_{1}^{\star},\ldots,f_{K}^{\star}$ has an eigenvalue close to $0$.
On the theoretical side, the proof of Theorem \ref{LEMMA:INEQ2} uses the linear independence of the emission densities by using that Gram matrices are positive.
An interesting problem  would be to investigate if  it is possible to estimate the emission densities with the classical adaptive rate for density estimation when the emission densities are linearly dependent (though all distinct). It is possible using model selection to  get the classical rate for the estimation of the density of $2K+1$ consecutive observations, but it does not seem obvious to see whether this rate can be  transferred to the estimators of the emission densities.
This is the subject of further work, see \cite{Luc2}.
\\
Another question arising from our work is whether it is possible to adapt to different smoothnesses of the emission densities.

\section{Proofs}
\label{sec:proofs}

\subsection{Proof of lemma \ref{identifiabilityL2}}
\label{sec:proof-iden}
In \cite{hsu2012spectral} it is proved that when {\bf [H1]}, {\bf [H2]}, {\bf [H3]} hold and when the rank of the matrix $\O:=(\langle \varphi_{m},f^{\star}_{k})_{1\leq m \leq M, 1\leq k\leq K}$ is $K$, 
the knowledge of the tensor $\M$ given by $\M (a,b,c)=\mathbb{E}(\varphi_{a}(Y_{1})\varphi_{b}(Y_{2})\varphi_{c}(Y_{3}))$ for all $a,b,c$ in $\{1,\ldots,M\}$ allows to recover $\O$ and $\Q$ up to relabelling of the hidden states.
 Thus, when {\bf [H1]}, {\bf [H2]}, {\bf [H3]} and {\bf [H4]} hold, the knowledge of $g^{\Q,\f^{\star}}$ is equivalent to the knowledge of the sequence $(\M)_{M}$, which allows to recover $\Q$ and the sequence $(\O)_{M}$, up to relabelling of the hidden states, which allows to recover
$\f^{\star}=(f_{1}^{\star},\ldots,f_{K}^{\star})$ up to relabelling of the hidden states, thanks to  \eqref{eq:ConvergenceL2}.
See also \cite{gassiat2013finite}. 

\subsection{Proof of Theorem \ref{prop:oracle}}
\label{sec:proof-oracle}

Throughout the proof $N$ is fixed, and we write $\gamma$ (instead of $\gamma_{N}$) for the contrast function.

\subsubsection{Beginning of the proof: algebraic manipulations}~\\
Let us fix some $M$ and some permutation $\tau$. 
Using the definitions of $\hat g_M$ and $\hat M$, we can write
\[
\gamma(\hat g_{\hat M})+\pen(N,\hat M)\leq \gamma(\hat g_{M})+\pen(N,M)\leq \gamma(g^{\hat\Q,  \f^\star_{M,\tau^{-1}}})+\pen(N,M)\,,
\]
where $ \f^\star_{M,\tau^{-1}}=( f^\star_{M,\tau^{-1}(1)}, \dots, f^\star_{M,\tau^{-1}(K)})$ (here we use that $ \f^\star_{M,\tau^{-1}}\in\mathcal F^K$).
But we can compute for all functions $t_1,t_2$,
\[
\gamma(t_1)-\gamma(t_2)=\lVert t_1-g^\star\lVert_{2}^2-\lVert t_2-g^\star\lVert_{2}^2-2\nu(t_1-t_2)\,,
\]
where $\nu$ is  the centered empirical process
\[
\nu(t)=\frac1N\sum_{s=1}^Nt(Y_1^{(s)},Y_2^{(s)},Y_3^{(s)})-\int t g^{\star}\,.
\]
This  gives
\begin{equation}\label{ee}
\lVert  \hat g_{\hat M}-g^{\star}\lVert_{2}^2\leq \lVert g^{\hat\Q,  \f^\star_{M,\tau^{-1}}}-g^{\star}\lVert_{2}^2+2\nu(\hat g_{\hat M}-
g^{\hat\Q,  \f^\star_{M,\tau^{-1}}})+\pen(N,M)-\pen(N,\hat M)
\end{equation}
Now, we denote by
$B_M=\lVert g^{\Q^\star, \f^\star_{M}}-g^{\star}\lVert_{2}^2$ a bias term
and we notice that $g^{\hat\Q,  \f^\star_{M,\tau^{-1}}}=g^{\P_{\tau}\hat\Q\mathbb P_{\tau}^{\top},  \f^\star_{M}}$. Then 
\begin{eqnarray*}
\lVert g^{\hat \Q, \f^\star_{M,\tau^{-1}}}-g^{\star}\lVert_{2}^2&\leq& 2\lVert g^{\hat \Q, \f^\star_{M,\tau^{-1}}}-g^{\Q^\star, \f^\star_{M}}\lVert_{2}^2+2\lVert g^{\Q^\star, \f^\star_{M}}-g^{\star}\lVert_{2}^2\\
&\leq& 2\lVert g^{\P_{\tau}\hat\Q\mathbb P_{\tau}^{\top},  \f^\star_{M}}-g^{\Q^\star, \f^\star_{M}}\lVert_{2}^2+2B_M.\\
\end{eqnarray*}
But, using Schwarz inequality,  $\lVert g^{\Q_1,  \f^\star_{M}}-g^{\Q_2, \f^\star_{M}}\lVert_{2}^2$ can be bounded by
\begin{align}
\label{qQ1-gQ2}
& \sum_{m_1,m_2,m_3=1}^M
 \Big|\sum_{k_1,k_2,k_3=1}^K \left(\pi_1(k_1) \Q_1(k_1,k_2)\Q_1(k_2,k_3)-\pi_2(k_1) \Q_2(k_1,k_2)\Q_2(k_2,k_3)\right)\nonumber\\
&\quad\quad\quad\langle f^\star_{k_1},\varphi_{m_1}\rangle \langle f^\star_{k_2},\varphi_{m_2}\rangle \langle f^\star_{k_3},\varphi_{m_3}\rangle
\Big|^2\nonumber\\
 &\leq\left(\sum_{k_1,k_2,k_3=1}^K \left(\pi_1(k_1) \Q_1(k_1,k_2)\Q_1(k_2,k_3)-\pi_2(k_1) \Q_2(k_1,k_2)\Q_2(k_2,k_3)\right)^2\right)
\nonumber \\
 &\sum_{m_1,m_2,m_3=1}^M\,\sum_{k_1,k_2,k_3=1}^K
 \Big| \langle f^\star_{k_1},\varphi_{m_1}\rangle \langle f^\star_{k_2},\varphi_{m_2}\rangle \langle f^\star_{k_3},\varphi_{m_3}\rangle
\Big|^2\nonumber\\
 & \leq 3 K^{3} C_{\mathcal F,2}^6\left(\lVert \pi_1-\pi_2\lVert_2^2+2\lVert \Q_1- \Q_2\lVert_F^2\right)
\end{align}
so that 
\begin{eqnarray*}
\lVert g^{\hat \Q, \f^\star_{M,\tau^{-1}}}-g^{\star}\lVert_{2}^2
&\leq& 6 K^{3} C_{\mathcal F,2}^6\left(\lVert \P_{\tau}\hat\pi-\pi^\star\lVert_2^2+2\lVert \P_{\tau}\hat\Q\mathbb P_{\tau}^{\top}-\Q^{\star}\lVert_F^2\right)
+2B_M.\\
\end{eqnarray*}
Next we set $S_{M}=\cup_{\Q}\dS ({\Q},{M})$ and 
\[ Z_{M}=\sup_{t\in S_{M}}\left[\frac{ |\nu(t-g^\star)|}{\lVert t-g^\star\lVert_{2}^2+x_{M}^2}\right]\]
for $x_{M}$ to be determined later. Then
\begin{eqnarray*}
\nu(\hat g_{\hat M}-g^{\hat\Q,  \f^\star_{M,\tau^{-1}}})&=&\nu( \hat g_{\hat M}-g^{\star})
+\nu(g^{\star}-g^{\hat\Q,  \f^\star_{M,\tau^{-1}}})\\
&\leq& Z_{\hat M}(\lVert \hat g_{\hat M}-g^\star\lVert_{2}^2+x_{\hat M}^2)+
Z_M(\lVert g^{\hat\Q,  \f^\star_{M,\tau^{-1}}}-g^\star\lVert_{2}^2+x_{M}^2).
\end{eqnarray*}
Denoting by $R_{\hat M}=\lVert \hat g_{\hat M}-g^\star\lVert_{2}^2$ the squared risk,  \eqref{ee} becomes
\begin{eqnarray*}
R_{\hat M}&\leq& 6 K^{3}C_{\mathcal F,2}^6\left(\lVert \P_{\tau}\hat\pi-\pi^\star\lVert_2^2+2\lVert \P_{\tau}\hat\Q\mathbb P_{\tau}^{\top}-\Q^{\star}\lVert_F^2\right)+2B_M
+2Z_{\hat M}(R_{\hat M}+x_{\hat M}^2)\\&&
+2Z_M\left(6 K^{3}C_{\mathcal F,2}^6\left(\lVert \P_{\tau}\hat\pi-\pi^\star\lVert_2^2+2\lVert \P_{\tau}\hat\Q\mathbb P_{\tau}^{\top}-\Q^{\star}\lVert_F^2\right)+2B_M+x_{M}^2\right)
\\&&+2\pen(N,M)-\pen(N,\hat M)-\pen(N, M)\,,\\
R_{\hat M}(1-2Z_{\hat M})&\leq &(2+4Z_{ M})B_M+2\pen(N,M)
\\&&+(1+2Z_M)6 K^{3}C_{\mathcal F,2}^6\left(\lVert \P_{\tau}\hat\pi-\pi^\star\lVert_2^2+2\lVert \P_{\tau}\hat\Q\mathbb P_{\tau}^{\top}-\Q^{\star}\lVert_F^2\right)\\
&&+2\sup_{M'}(2Z_{ M'}x_{M'}^2-\pen(N,M'))\,.
\end{eqnarray*}

\noindent
To conclude it is then sufficient to establish that, with probability larger than $1- (e-1)^{-1}e^{-x}$, it holds
\[ \sup_{M'} Z_{M'}\leq\frac14\quad\text{ and }\quad \sup_{M'}(2Z_{ M'}x_{M'}^2-\pen(N,M'))\leq A\frac{x}{N}\,,\]
with $A$ a constant depending only on $\Q^\star$ and $\f^\star$ and not on $N,M,x$.
Thus we will have, for any $M$, with probability larger than $1- (e-1)^{-1}e^{-x}$,
\begin{eqnarray*}
\frac{1}{2}R_{\hat M}&\leq &3 B_M+2\pen(N,M)+2A\frac{x}{N}\\
&&+9 C_{\mathcal F,2}^6\left(\lVert \P_{\tau}\hat\pi-\pi^\star\lVert_2^2+2\lVert \P_{\tau}\hat\Q\mathbb P_{\tau}^{\top}-\Q^{\star}\lVert_F^2\right)
\end{eqnarray*}
which is the announced result.

The heart of the proof is then the study of $Z_M$. 
We introduce $u_M$ a projection of $g^\star$ on $S_M$
and we split $Z_M$ into two terms: $Z_M\leq 4 Z_{M,1}+Z_{M,2}$ with 
\[ 
\begin{cases}
 \displaystyle  Z_{M,1}=\sup_{t\in S_{M}}\left[\frac{ |\nu(t-u_M)|}{\lVert t-u_M\lVert_{2}^2+4x_{M}^2}\right]\\
    \displaystyle  Z_{M,2}=\frac{ |\nu(u_M-g^\star)|}{\lVert u_M-g^\star\lVert_{2}^2+x_{M}^2}
  \end{cases}
\]
Indeed $u_M$ verifies: for all $t\in S_M$, 
\[\lVert u_M-g^\star\lVert_{2}\leq \lVert  t-g^\star\lVert_{2}\quad \mathrm{and} \quad\lVert u_M-t\lVert_{2}\leq 2\lVert t-g^\star\lVert_{2}\,.\]

\subsubsection{Deviation inequality for $Z_{M,2}$}
~\\
Bernstein's inequality \eqref{eq:bernstein} for HMMs (see Appendix \ref{sec:AppendixConcentration})  gives, with probability larger than $1-e^{-z}$: 
\[  |\nu(u_M-g^\star)|
\leq 2\sqrt{2c^{\star}\lVert u_M-g^\star\lVert_{2}^2\lVert g^\star\lVert_\infty\frac{z}{N}}+2\sqrt{2}c^{\star}\lVert u_M-g^\star\lVert_\infty\frac{z}{N}\,.\]
Then, using $a^2+b^2\geq 2ab$, with probability larger than $1-e^{-z}$: 
\begin{eqnarray*}
\frac{ |\nu(u_M-g^\star)|}{\lVert u_M-g^\star\lVert_{2}^2+x_{M}^2}
&\leq & 2\sqrt{2c^{\star}\lVert g^\star\lVert_\infty}\frac{1}{2x_{M}} \sqrt{\frac{z}{N}}
+2\sqrt{2}c^{\star}\frac{\lVert u_M\lVert_\infty +\lVert g^\star\lVert_\infty}{x_{M}^2}\frac{z}{N}\,.
\end{eqnarray*}
But any function $t$ in  $S_M$ can be written 
\[ t=\sum_{k_1,k_2,k_3=1}^K \pi(k_1) \Q(k_1,k_2)\Q(k_2,k_3)f_{k_1}\otimes f_{k_2}\otimes f_{k_3}\,,\]
with $f_k\in\mathcal{F}$ for $k=1,\dots, K$, so that $\sup_{t\in S_M}\lVert t\lVert_\infty\leq C_{\mathcal F,\infty}^3$.
Then, with probability larger than $1-e^{-z_M-z}$
\begin{eqnarray}\label{zm2}
Z_{M,2}\leq 
\sqrt{2c^{\star}\lVert g^\star\lVert_\infty}\sqrt{\frac{z_M+z}{x_M^2N}}+
4\sqrt{2}c^{\star}C_{\mathcal F,\infty}^3 \frac{z_M+z}{x_M^2N}\,.
\end{eqnarray}

\subsubsection{Deviation inequality for $Z_{M,1}$}
~\\
We shall first study the term $\sup_{t\in B_\sigma} |\nu(t-u_M)|$
where \[ B_\sigma=\{ t\in S_{M},\lVert t-u_M\lVert_{2}\leq \sigma\}.\]
Remark that, for all $t\in \dS ({\Q},{M})$, 
\[ \lVert t\lVert_{2}^2\leq \sum_{k_1,k_2,k_3=1}^K \pi^2(k_1) \Q^2(k_1,k_2)\Q^2(k_2,k_3) \sum_{k_1,k_2,k_3=1}^K 
C_{\mathcal F,2}^2C_{\mathcal F,2}^2C_{\mathcal F,2}^2\leq K^3
C_{\mathcal F,2}^6.\]
Then, if $t\in B_\sigma$, 
$\lVert t-u_M\lVert_{2}\leq \sigma \wedge 2K^{3/2}
C_{\mathcal F,2}^3$.
Notice also that for all  $t\in S_M$, $\lVert t-u_M\lVert_\infty\leq 2C_{\mathcal F,\infty}^3$.
Now Proposition~\ref{prop:deviation} in Appendix~\ref{sec:AppendixConcentration} (applied to a countable dense set  in $B_\sigma$) gives
that for any measurable set $A$ such that $\P(A)>0$, 
\[ \E^{A}(\sup_{t\in B_\sigma} |\nu(t-u_M)|)\leq C^\star\left[ \frac{E}{ N}+\sigma \sqrt{\frac1N\log\left(\frac{1}{\P(A)}\right)}+
\frac{2C_{\mathcal F,\infty}^3}{N}\log\left(\frac{1}{\P(A)}\right)\right],\]
and 
\[ E=\sqrt{N}\int_0^\sigma \sqrt{H(u)\wedge N}du +(2C_{\mathcal F,\infty}^3+2K^{3/2}
C_{\mathcal F,2}^3) H(\sigma) \,.\]
Here, for any integrable random variable $Z$, $E^{A}[Z]$ denotes $E[Z\one_{A}]/\P(A)$.

We shall compute $E$ later and find $\sigma_M$ and $\varphi$ such  that 
\begin{equation}\label{bornent}
\forall \sigma\geq \sigma_M\qquad
E\leq (1+2C_{\mathcal F,\infty}^3+2K^{3/2}C_{\mathcal F,2}^3)\varphi(\sigma)\sqrt{N}.
\end{equation}
(see Section~\ref{secentro}).
We then use Lemma 4.23 in \cite{MR2319879}  to write (for $x_M\geq \sigma_M$)
\[ \E^{A}\left(\sup_{t\in S_{M}}\left[\frac{ |\nu(t-u_M)|}{\lVert t-u_M\lVert_{2}^2+4x_{M}^2}\right]\right)
\leq \frac{C^\star}{x_M^2}\left[  C\frac{\varphi(2x_M)}{\sqrt N}+2x_M \sqrt{\frac1{N}\log\left(\frac{1}{\P(A)}\right)}+
\frac{2C_{\mathcal F,\infty}^3}{N}\log\left(\frac{1}{\P(A)}\right)\right]\]
Finally, Lemma 2.4 in \cite{MR2319879}  ensures that, with probability $1-e^{-z_M-z}$:
\begin{equation}\label{zm1}
Z_{M,1}=\sup_{t\in S_{M}}\left[\frac{ |\nu(t-u_M)|}{\lVert t-u_M\lVert_{2}^2+4x_{M}^2}\right]
\leq  C^\star\left[  C\frac{\varphi(2x_M)}{x_M^2\sqrt N}+ 2\sqrt{\frac{z_M+z}{x_M^2N}}+
2C_{\mathcal F,\infty}^3\frac{z_M+z}{x_M^2N}\right].
\end{equation}

\subsubsection{Computation of the entropy and function $\varphi$}~\\
\label{secentro}%
The definition of $H$ given in Proposition~\ref{prop:deviation} shows that 
$H(\delta)$ is bounded by the classical bracketing entropy for $\mathbf{L}^{2}$ distance at point $\delta/C_{\mathcal F,\infty}^3$ 
(where $C_{\mathcal F,\infty}^3$ bounds the sup norm of $g^\star$):
$H(\delta)\leq H(\delta/C_{\mathcal F,\infty}^3,S_M, \mathbf{L}^{2})$. We denote by $N(u,S, \mathbf{L}^{2})=e^{H(u,S,\mathbf{L}^{2})}$
the minimal number of brackets of radius $u$ to cover $S$. Recall that when $t_{1}$ and $t_{2}$ are real valued functions, the bracket $[t_{1},t_{2}]$ is the set of real valued functions $t$ such that $t_{1}(\cdot)\leq t(\cdot) \leq t_{2}(\cdot)$, and the radius of the bracket is $\|t_{2}-t_{1}\|_{2}$.
Now, observe that $S_M=\cup_{\Q} \dS ({ \Q},{M})$ is a set of mixtures of parametric functions. Denoting $\mathbf k=(k_1,k_2,k_3)$, 
$S_M$ is included in 
\begin{multline*}
\left\{\sum_{\mathbf k\in \{1,\dots,K\}^3}\mu(\mathbf k)f_{k_1}\otimes f_{k_2}\otimes f_{k_3},\
 \mu\geq 0,\
 \sum_{\mathbf k\in \{1,\dots,K\}^3} \mu(\mathbf k)=1,\right.\\
\left.  f_{k_i}\in {\mathcal F}\cap {\rm Span}(\varphi_1,\dots, \varphi_M),\;i=1,2,3\right\}\,.
\end{multline*}
Set 
$$
\mathcal{A}
=\{f_{1}\otimes f_{2}\otimes f_{3}, f_{i}\in {\mathcal F}\cap {\rm Span}(\varphi_1,\dots, \varphi_M), \;i=1,2,3\}.
$$
Then following the proof in Appendix A of \cite{BonTou2013}, we can prove 
\begin{equation}
\label{eq:bracketM}
 N(\varepsilon,S_M, \mathbf{L}^{2})\leq \left(\frac{C_1}{\varepsilon}\right)^{K^3-1}
\left[N\left(\frac{\varepsilon}{3},\mathcal{A}, \mathbf{L}^{2}\right)\right]^{K^{3}}\,.
\end{equation}
where $C_1$ depends on $K$ and $C_{\mathcal F,2}$. 
Denote $\mathcal{B}
= {\mathcal F}\cap {\rm Span}(\varphi_1,\dots, \varphi_M)
$.
Let $a=(a_{m})_{1\leq m \leq M}\in\R^{M}$ and $b=(b_{m})_{1\leq m \leq M}\in\R^{M}$ such that $a_{m}<b_{m}$, $m=1,\ldots,M$.  For each $m=1,\ldots,M$ and $y\in\Y$, let
\begin{align*}
u_{m}(y) &=
	\begin{cases}
	a_{m} \qquad \text{if } \varphi_{m}(y) \geq 0 \\
	b_{m} \qquad \text{otherwise}
	\end{cases} \\
v_{m}(y) &= a_{m} + b_{m} - u_{m}(y).
\end{align*}
Then, if $(c_{m})_{1\leq m \leq M}\in\R^{M}$ is such that for all $m=1,\ldots,M$, $a_{m}\leq c_{m} \leq b_{m}$, then
\begin{equation*}
U_{a,b}^{1}(y): = \sum_{m=1}^M u_{m}(y) \varphi_{m}(y)
\leq \sum_{m=1}^M c_{m} \varphi_{m}(y)
\leq \sum_{m=1}^M v_{m}(y) \varphi_{m}(y) = U^{2}_{a,b}(y).
\end{equation*}
Moreover,
\begin{align*}
\|U^{2}_{a,b} - U^{1}_{a,b}\|_{2}^{2}
	&=\| \sum_{m=1}^M |b_{m} - a_{m}|. |\varphi_{m}| \|_{2}^{2}\\
	&\leq M \|b-a\|_{2}^{2}
\end{align*}
using Cauchy-Schwarz inequality. Thus, one may cover $\mathcal{B}$ with brackets of form  $[U^{1}_{a,b},U^{2}_{a,b}]$.
Also, for $i=1,2$,
\begin{align*}
\|U^{i}_{a,b}\|_{2}^{2}
	&\leq\| \sum_{m_i=1}^M |b_{m} + a_{m}|. |\varphi_{m}| \|_{2}^{2}\\
	&\leq 2 M (\|a\|_{2}^{2}+\|b\|_{2}^{2}).
\end{align*}
If now for some $a^{i}$, $b^{i}$ in $\R^{M}$, $f_{i}\in [U^{1}_{a^{i},b^{i}},U^{2}_{a^{i},b^{i}}]$,  $i=1,2,3$, then
$$
f_{1}\otimes f_{2}\otimes f_{3} \in \left[ V,W  \right]
$$
with
$$
V= \min \{U^{i_{1}}_{a^{1},b^{1}}U^{i_{2}}_{a^{2},b^{2}}U^{i_{3}}_{a^{3},b^{3}},\;i_{1},i_{2},i_{3} \in \{1,2\}   \}
$$
and
$$
W= \max \{U^{i_{1}}_{a^{1},b^{1}}U^{i_{2}}_{a^{2},b^{2}}U^{i_{3}}_{a^{3},b^{3}},\;i_{1},i_{2},i_{3} \in \{1,2\}   \},
$$
pointwise.
Moreover, one can see that
\begin{eqnarray*}
\left\vert W-V\right\vert & \leq& \left \vert U^{2}_{a^{1},b^{1}}-U^{1}_{a^{1},b^{1}}\right\vert \max_{j_{1},j_{2}\in \{1,2\}}  \left \vert U^{j_{1}}_{a^{2},b^{2}}\right\vert .\left \vert U^{j_{2}}_{a^{3},b^{3}}\right\vert\\
&&+\left \vert U^{2}_{a^{2},b^{2}}-U^{1}_{a^{2},b^{2}}\right\vert \max_{j_{1},j_{2}\in \{1,2\}}  \left \vert U^{j_{1}}_{a^{1},b^{1}}\right\vert .\left \vert U^{j_{2}}_{a^{3},b^{3}}\right\vert\\
&&
+\left \vert U^{2}_{a^{3},b^{3}}-U^{1}_{a^{3},b^{3}}\right\vert \max_{j_{1},j_{2}\in \{1,2\}}  \left \vert U^{j_{1}}_{a^{1},b^{1}}\right\vert .\left \vert U^{j_{2}}_{a^{2},b^{2}}\right\vert\\
&\leq &\sum_{i=1}^{3}\left\vert U^{2}_{a^{i},b^{i}}-U^{1}_{a^{i},b^{i}}\right\vert \prod_{j\neq i}\left( \left \vert U^{1}_{a^{j},b^{j}}\right\vert+\left \vert U^{2}_{a^{j},b^{j}}\right\vert\right)
\end{eqnarray*}
so that
\begin{eqnarray*}
\left\| W-V\right\|_{2}^{2} & \leq& 12\sum_{i=1}^{3}\left\| U^{2}_{a^{i},b^{i}}-U^{1}_{a^{i},b^{i}}\right\|_{2}^{2} \prod_{j\neq i}\left( \left \| U^{1}_{a^{j},b^{j}}\right\|_{2}^{2}+\left \|U^{2}_{a^{j},b^{j}}\right\|_{2}^{2}\right)\\
&\leq & 48 M^{3} \sum_{i=1}^{3} \|b^{i}-a^{i}\|_{2}^{2}\prod_{j\neq i}\left( \|a^{j}\|_{2}^{2}+\|b^{j}\|_{2}^{2}\right)\\
&\leq & 192 M^{3} C_{\mathcal{F},2}^{4}\sum_{i=1}^{3} \|b^{i}-a^{i}\|_{2}^{2}.
\end{eqnarray*}
Thus one may cover $\mathcal A$ by covering the ball of radius $C_{\mathcal{F},2}$ in $\R^{M}$ with hypercubes $[a,b]$,  for which $\|a\|_{2}$, $\|b\|_{2}$ are less than $C_{\mathcal{F},2}$.
To get a bracket with radius $u$, it is enough that $\|b^{i}-a^{i}\|_{2}^{2} \leq u^{2} / (576 M^{3} C_{\mathcal{F},2}^{4})$, $i=1,2,3$.
We finally obtain that
\begin{equation}
\label{eq:bracketA}
N\left(u,\mathcal{A}, \mathbf{L}^{2}\right)\leq \left(\frac{48 \sqrt{3} M^{3/2} C_{\mathcal{F},2}^3}{u}
\right)^{3M}.
\end{equation}
We deduce from (\ref{eq:bracketM}) and (\ref{eq:bracketA}) that 
\[ N(u, S_M, \mathbf{L}^{2})\leq \left(\frac{C_1}{u}\right)^{K^3-1}
\left(\frac{48 \sqrt{3} M^{3/2} C_{\mathcal{F},2}^3}{u}
\right)^{3MK^3}\,,\]
and then
\[ H(u, S_M,\mathbf{L}^{2})\leq (K^3-1)\log(\frac{C_1}u)+3MK^3\log\left(\frac{C_2M^{3/2} }{u}\right)\,,\]
with $C_2$ depending on $K$ and $C_{\mathcal F,2}$.
To conclude we use that
$\int_0^{\sigma}\sqrt{\log\left(\frac{1}{x}\right)}dx\leq \sigma(\sqrt{\pi}+\sqrt{\log\left(\frac{1}{\sigma}\right)})$, see \cite{BMM12}.
Finally we can write for $\sigma \leq M^{3/2}$:
\[ \int_0^\sigma \sqrt{H(u)}du\leq
C_3\sqrt{M}\sigma\left(1+\sqrt{\log\left(\frac{M^{3/2}}{\sigma}\right)}\right),\]
where $C_3$ depends on $K$, $C_{\mathcal F,2}$ and $C_{\mathcal F,\infty}$.
Set \[ \varphi(x)= C_3\sqrt{M}x\left(1+\sqrt{\log\left(\frac{M^{3/2}}{x}\right)}\right)\]
The function
$\varphi$ is increasing on $]0,M^{3/2}]$, and $\varphi(x)/x$ is decreasing. 
Moreover $\varphi(\sigma)\geq \int_0^\sigma \sqrt{H(u)}du$ and  $\varphi^2(\sigma)\geq \sigma^2 H(\sigma)$. \\

\subsubsection{End of the proof, choice of parameters}~\\
As soon as $N\geq C_{3}^{2}/M^{2}:=N_{0}$, we may define $\sigma_M$ as the solution of equation $\varphi(x)=\sqrt{N} x^2$. 
Then, for all $\sigma\geq \sigma_M$, \[ H(\sigma)\leq \frac{\varphi(\sigma)^2}{\sigma^2}
\leq \frac{\varphi(\sigma)}{\sigma}\sigma \sqrt{N}.\]
This yields, for all $\sigma\geq \sigma_M,$
\[ E\leq (1+2C_{\mathcal F,\infty}^3+2K^{3/2}C_{\mathcal F,2}^3)\varphi(\sigma)\sqrt{N},\] 
which was required in \eqref{bornent}.

Moreover $\frac{\varphi(2x_M)}{x_M\sqrt N}\leq 2\sigma_M$ as soon as $x_M\geq \sigma_M$.
Combining \eqref{zm1} and \eqref{zm2}, we obtain,
with probability $1-e^{-z_M-z}$:
\[ Z_{M}\leq C^{\star\star} \left[\frac{\sigma_M}{x_M}+\sqrt{\frac{z_M+z}{x_M^2N}}+
\frac{z_M+z}{x_M^2N}\right]\,,\]
where $C^{\star\star} $ depends on $K$,  $C_{\mathcal F,2}$, $C_{\mathcal F,\infty}$, $\Q^\star$. 
Now let us choose
$x_M=\theta^{-1}\sqrt{\sigma_M^2+\frac{z_M+z}{N}}$ with $\theta$ such that 
$2\theta+\theta^2\leq(C^{\star\star})^{-1}/4$.
This choice entails:
$x_M\geq \theta^{-1}\sigma_M$ and $x_M^2\geq\theta^{-2}\frac{z_M+z}{N}.$
Then with probability $1-e^{-z_M-z}$:
\[ Z_{M}\leq C^{\star\star}( \theta+\theta+
\theta^2). \]
We now choose $z_M=M $ which implies $\sum_{M\geq 1} e^{-z_M}=(e-1)^{-1}$. 
Then, with probability $1-(e-1)^{-1} e^{-z}$,
\[ \forall M \qquad Z_{M}\leq C^{\star\star}(2 \theta+\theta^2)\leq \frac{1}{4}\,,\]
and for all $M$,
\begin{eqnarray*}
 Z_{M}x_M^2&\leq &C^{\star\star} \left[{\sigma_M}{x_M}+x_M\sqrt{\frac{z_M+z}{N}}+
\frac{z_M+z}{N}\right]\\
&\leq &C^{\star\star}\theta^{-1} \left(\sigma_M+\sqrt{\frac{z_M+z}{N}}\right)^2+
C^{\star\star}\frac{z_M+z}{N}\,.
\end{eqnarray*}
Then, with probability $1-(e-1)^{-1} e^{-z}$, for all $M$,
\[ 
Z_{M}x_M^2-C^{\star\star}\left(2\theta^{-1} \sigma_M^2+(2\theta^{-1}+1)\frac{M}{N}\right)
\leq
C^{\star\star}(2\theta^{-1}+1)\frac{z}{N}\,.
\]
Then the result is proved as soon as 
\begin{equation}
\label{eq:penlower}
\pen(N,M)\geq 2C^{\star\star}\left(2\theta^{-1} \sigma_M^2+(2\theta^{-1}+1)\frac{M}{N}\right). 
\end{equation}
It remains to get an upper bound for $\sigma_M$.
Recall that $\sigma_M$ is defined as the solution of equation $C_3\sqrt{M}x(1+\sqrt{\log\left(\frac{M^{3}}{x}\right)})=\sqrt{N} x^2$.
Then we obtain that for some $C_{4}$ 
\[ \sigma_M\leq C_4\sqrt{\frac{M}{N}}(1+\sqrt{\log(N)})\,,\]
and (\ref{eq:penlower}) holds as soon as
$$\pen(N,M)\geq \rho^{\star}\frac{M\log(N)}{N}
$$
for some constant $\rho^{\star}$ depending on $C_{\mathcal F,2}$ and $C_{\mathcal F,\infty}$ ({\bf Scenario A}) or on $\Q^{\star}$, $C_{\mathcal F,2}$ and $C_{\mathcal F,\infty}$  ({\bf Scenario B}).

\subsection{Proof of Theorem \ref{LEMMA:INEQ2}}
\label{sec:proof-ineq2}

\noindent
For any $\h \in {\mathcal K}^{K}$ and $\Q\in  {\mathcal V}$, denote $N(\Q,\h)=\lVert g^{\Q,\f^{\star}+\h}-g^{\Q,\f^{\star}}\lVert_{2}^{2}$. What we want to prove is that
\[ 
c:=c({\mathcal K}, {\mathcal V}, \F^{\star})^{2}:=\inf_{\Q\in {\mathcal V},\h \in {\mathcal K}^{K},\lVert {\h}\lVert_{2}\neq 0}\frac{N({\Q},{\h})}{\lVert {\h}\lVert_{{2}}^{2}}
>0.
\]
One can compute:
\begin{multline*}
N(\Q,\h)=\sum_{k_{1},k_{2},k_{3},k'_{1},k'_{2},k'_{3}=1}^{K}\pi(k_1) \Q(k_1,k_2)\Q(k_{2},k_3)\pi(k'_1) \Q(k'_1,k'_2)\Q(k'_{2},k'_3)\\
\left(\prod_{i=1}^{3}\langle f^{\star}_{k_{i}}+h_{k_{i}},f^{\star}_{k'_{i}}+h_{k'_{i}}\rangle +\prod_{i=1}^{3}\langle f^{\star}_{k_{i}},f^{\star}_{k'_{i}}\rangle
-\prod_{i=1}^{3}\langle f^{\star}_{k_{i}}+h_{k_{i}},f^{\star}_{k'_{i}}\rangle-\prod_{i=1}^{3}\langle f^{\star}_{k_{i}},f^{\star}_{k'_{i}}+h_{k'_{i}}\rangle \right).
\end{multline*}
Let ${\mathbf{u}}=(u_{1},\ldots , u_{K})$ be such that $u_{i}$, $i=1,\ldots,K$, is the orthogonal projection of $h_{i}$ on the subspace of $\mathbf{L}^{2}(\Y,\L)$ spanned by $f_{1}^{\star}\ldots,f_{K}^{\star}$. Then
\begin{equation}
\label{eq:decompo}
N(\Q,\h)=N(\Q,\mathbf{u})+M(\Q,\mathbf{u},\h-\mathbf{u})
\end{equation}
where, for any $\mathbf{a}=(a_{1},\ldots,a_{K})\in \mathbf{L}^{2}(\Y,\L)^{K}$,
\begin{multline*}
M(\Q,\mathbf{u},\mathbf{a})=\sum_{k_{1},k_{2},k_{3},k'_{1},k'_{2},k'_{3}=1}^{K}
\pi(k_1) \Q(k_1,k_2)\Q(k_{2},k_3)\pi(k'_1) \Q(k'_1,k'_2)\Q(k'_{2},k'_3)\\
\left(\prod_{i=1}^{3}\langle a_{k_{i}},a_{k'_{i}}\rangle +\sum_{i=1}^{3}\langle a_{k_{i}},a_{k'_{i}}\rangle\prod_{j\neq i}\langle f^{\star}_{k_{j}}+u_{k_{j}},f^{\star}_{k'_{j}}+u_{k'_{j}}\rangle
+\sum_{i=1}^{3}\langle f^{\star}_{k_{i}}+u_{k_{i}},f^{\star}_{k'_{i}}+u_{k'_{i}}\rangle\prod_{j\neq i}\langle a_{k_{j}},a_{k'_{j}}\rangle \right).
\end{multline*}
Let $A=\D\pi$ with $\pi$ the stationary distribution of $\Q$. Then $M(\Q,\mathbf{u},\mathbf{a})$ may be rewritten as:
\begin{eqnarray*}
M(\Q,\mathbf{u},\mathbf{a})&=\sum_{i,j=1}^{K}& \langle (\Q^{T}A \mathbf{a})_{i}, ( \Q^{T}A \mathbf{a})_{j}\rangle \langle {a}_{i},  {a}_{j}\rangle
\langle (\Q  \mathbf{a})_{i},  (\Q  \mathbf{a})_{j}\rangle\\
&&+\langle (\Q^{T}A \mathbf{a})_{i},  (\Q^{T}A \mathbf{a})_{j}\rangle \langle (\f^{\star}+\mathbf{u}_{i}),  (\f^{\star}+\mathbf{u})_{j}\rangle
\langle (\Q  (\f^{\star}+\mathbf{u}))_{i},  (\Q  (\f^{\star}+\mathbf{u}))_{j}\rangle
\\
&&+ \langle (\Q^{T}A (\f^{\star}+\mathbf{u}))_{i}, ( \Q^{T}A (\f^{\star}+\mathbf{u}))_{j}\rangle \langle {a}_{i},  {a}_{j}\rangle
\langle (\Q  (\f^{\star}+\mathbf{u}))_{i},  (\Q (\f^{\star}+\mathbf{u}))_{j}\rangle
\\
&&+\langle (\Q^{T}A (\f^{\star}+\mathbf{u}))_{i},  (\Q^{T}A(\f^{\star}+\mathbf{u}))_{j}\rangle \langle (\f^{\star}+\mathbf{u})_{i},  (\f^{\star}+\mathbf{u})_{j}\rangle
\langle (\Q  \mathbf{a})_{i},  (\Q  \mathbf{a})_{j}\rangle\\
&&+ \langle (\Q^{T}A (\f^{\star}+\mathbf{u}))_{i}, ( \Q^{T}A (\f^{\star}+\mathbf{u}))_{j}\rangle \langle {a}_{i},  {a}_{j}\rangle
\langle (\Q  \mathbf{a})_{i},  (\Q  \mathbf{a})_{j}\rangle
\\
&&+\langle (\Q^{T}A \mathbf{a})_{i},  (\Q^{T}A \mathbf{a})_{j}\rangle \langle (\f^{\star}+\mathbf{u})_{i}, (\f^{\star}+\mathbf{u})_{j}\rangle
\langle (\Q  \mathbf{a})_{i},  (\Q  \mathbf{a})_{j}\rangle
\\
&&+\langle (\Q^{T}A \mathbf{a})_{i}, ( \Q^{T}A \mathbf{a})_{j}\rangle \langle {a}_{i},  {a}_{j}\rangle
\langle (\Q  (\f^{\star}+\mathbf{u}))_{i},  (\Q (\f^{\star}+\mathbf{u}))_{j}\rangle.
\end{eqnarray*}
All terms in this sum are non negative. Let us prove it for the first one, the argument for the others is similar. Define $V$ the $K\times K$ matrix given by
$$
V_{i,j}=\langle (\Q^{T}A \mathbf{a})_{i}, ( \Q^{T}A \mathbf{a})_{j}\rangle \langle (\Q  \mathbf{a})_{i},  (\Q  \mathbf{a})_{j}\rangle,\quad i,j=1,\ldots,K.
$$
$V$ is the Hadamard product of two Gram matrices which are non negative, thus $V$ is itself non negative by the Schur product Theorem, see \cite{Schur}, and
$$
\sum_{i,j=1}^{K} V_{i,j} \langle{a}_{i},  {a}_{j}\rangle=\int \mathbf{a}(y)^{T}V\mathbf{a}(y) dy \geq 0.
$$
Thus we have that $M(\Q,\mathbf{u},\mathbf{a})$ is lower bounded by one term of the sum so that
$$
M(\Q,\mathbf{u},\mathbf{a}) \geq \sum_{i,j=1}^{K}\langle (\Q^{T}A (\f^{\star}+\mathbf{u}))_{i}, ( \Q^{T}A (\f^{\star}+\mathbf{u}))_{j}\rangle \langle {a}_{i},  {a}_{j}\rangle
\langle (\Q  (\f^{\star}+\mathbf{u}))_{i},  (\Q (\f^{\star}+\mathbf{u}))_{j}\rangle.
$$
The minimal eigenvalue of the Hadamard product of two non negative matrices is lower bounded by the product of the minimal eignevalues of each matrix, and we get
\begin{equation}
\label{eq:Schur}
M(\Q,\mathbf{u},\mathbf{a}) \geq \left(\min_{i=1,\ldots,K}\lambda_{i} (\Q^{T}A (\f^{\star}+\mathbf{u}))\right)\left(\min_{i=1,\ldots,K}\lambda_{i} (\Q(\f^{\star}+\mathbf{u}))\right)
\|\mathbf{a}\|_{2}^{2}
\end{equation} 
where $\|\mathbf{a}\|_{2}^{2}=\sum_{k=1}^{k}\|a_{k}\|_{2}^{2}$, and where, if $\h \in \mathbf{L}^{2}(\Y,\L)^{K}$, $\lambda_{1}(\h)$, \ldots, $\lambda_{K}(\h)$ are the (non negative) eigenvalues of the Gram matrix of $h_{1},\ldots,h_{K}$.
\\

\noindent
 Let now $(\Q_{n},\h_{n})_{n}$ be a sequence in ${\mathcal V}\times {\mathcal K}$ such that $c=\lim_{n}\frac{N(\Q_{n},\mathbf{h}_{n})}{\lVert \mathbf{h}_{n}\lVert_{\Q^{\star}}^{2}}$. 
 Let ${\mathbf{u}}_{n}$ be the vector of the orthogonal projections of the coordinate functions of $\h_{n}$ on the subspace of $\mathbf{L}^{2}(\Y,\L)$ spanned by $f_{1}^{\star}\ldots,f_{K}^{\star}$. Notice that 
 $$
 \|\h_{n}\|_{\Q^{\star}}^{2}=\|{\mathbf{u}}_{n}\|_{\Q^{\star}}^{2}+\|\h_{n}-{\mathbf{u}}_{n}\|_{2}^{2}.
 $$
   Let $C_{{\cal K},2}$ be the upper bound of the norm of elements of $\cal K$. We have,  for any $n\geq 1$,
 $$
  \|\h_{n}\|_{\Q^{\star}}^{2} \leq K(C_{{\cal K},2}+2C_{{\cal F},2})^{2}
 $$
 so that for any $n\geq 1$,
 $$
 \|{\mathbf{u}}_{n}\|_{\Q^{\star}}^{2} \leq K(C_{{\cal K},2}+2C_{{\cal F},2})^{2}.
 $$
  Since  $(\Q_{n},{\mathbf{u}}_{n})_{n}$  is a bounded sequence in a finite dimensional space it has a limit point $(\Q,{\mathbf{u}})$.
Now,  using (\ref{eq:decompo}) and the non negativity of $M(\Q_{n},\mathbf{u}_{n},\h_{n}-\mathbf{u}_{n})$, we get on the converging subsequence
$$
c \geq \lim_{n\rightarrow +\infty}\frac{N(\Q_{n}, \mathbf{u}_{n})}{K(C_{{\cal K},2}+2C_{{\cal F},2})^{2}}= \frac{N(\Q,{\mathbf{u}})}{K(C_{{\cal K},2}+2C_{{\cal F},2})^{2}}.
$$
Since $\Q\in{\mathcal V}$, $T_{\Q}\subset T_{\Q^{\star}}$ so that $\lVert \mathbf{u}\lVert_{\Q}\geq \lVert \mathbf{u}\lVert_{\Q^{\star}}$. Thus if  $\lVert \mathbf{u}\lVert_{\Q^{\star}}\neq 0$, $\lVert \mathbf{u}\lVert_{\Q}\neq 0$, and using Lemma \ref{identifiabilityL2}, $N(\Q,{\mathbf{u}})\neq 0$ so that $c>0$ in this case.
 \\
 
\noindent
Consider now the situation where  $\lVert \mathbf{u}\lVert_{\Q^{\star}} = 0$. 
Since $\lim_{n\rightarrow +\infty}\lVert \mathbf{u}_{n}\lVert_{\Q^{\star}}=0$, there exists $n_{1}$ and $\tau\in T_{\Q^{\star}}$ such that for all $n\geq n_{1}$, one has $\lVert \mathbf{u}_{n}\lVert_{\Q^{\star}}^{2} =\sum_{k=1}^{K}\lVert u_{n,k}+f^{\star}_{k}-f^{\star}_{\tau(k)}\lVert^{2}_{2}$,
and it is possible to exchange the states in the transition matrix using $\tau$ so that we just have to consider the situation where $\lVert \mathbf{u}_{n}\lVert_{\Q^{\star}}^{2} = \lVert \mathbf{u}_{n}\lVert_{2}^{2}$ for large enough $n$.\\
Eigenvalues of Gram matrices of functions are continuous in the functions so that using  (\ref{eq:decompo}) and (\ref{eq:Schur}) we get 
\begin{multline*}
c \geq \lim_{n\rightarrow +\infty} \frac{N(\Q_{n}, \mathbf{u}_{n})}{\lVert \mathbf{u}_{n}\lVert_{2}^{2} +\|\h_{n}- \mathbf{u}_{n}\|_{2}^{2} }\\+
\left(\min_{i=1,\ldots,K}\lambda_{i} (\Q^{T}A \f^{\star})\right)\left(\min_{i=1,\ldots,K}\lambda_{i} (\Q\f^{\star})\right)\liminf_{n\rightarrow +\infty}\frac{\|\h_{n}- \mathbf{u}_{n}\|_{2}^{2}}{\lVert \mathbf{u}_{n}\lVert_{2}^{2} +\|\h_{n}- \mathbf{u}_{n}\|_{2}^{2}}.
\end{multline*}
Under assumptions {\bf[H1]} and {\bf[H4]}, $\Q^{T}A \f^{\star}$ is a vector of linearly independent functions and $\Q\f^{\star}$ also, so that
$$
\left(\min_{i=1,\ldots,K}\lambda_{i} (\Q^{T}A \f^{\star})\right)\left(\min_{i=1,\ldots,K}\lambda_{i} (\Q\f^{\star})\right) >0.
$$
Thus, if
$
\liminf_{n\rightarrow +\infty}\frac{\|\h_{n}-\mathbf{u}_{n}\|_{2}^{2}}{\lVert \mathbf{u}_{n}\lVert_{2}^{2} +\|\h_{n}-\mathbf{u}_{n}\|_{2}^{2}}>0
$
we obtain $c>0$.\\
If now $
\liminf_{n\rightarrow +\infty}\frac{\|\h_{n}-\mathbf{u}_{n}\|_{2}^{2}}{\lVert \mathbf{u}_{n}\lVert_{2}^{2} +\|\h_{n}-\mathbf{u}_{n}\|_{2}^{2}}=0$, 
we have on a subsequence
\begin{equation}
\label{eq:presque}
c \geq \lim_{n\rightarrow +\infty} \frac{N(\Q_{n}, \mathbf{u}_{n})}{\lVert \mathbf{u}_{n}\lVert_{2}^{2} }\frac{\lVert \mathbf{u}_{n}\lVert_{2}^{2}}{\lVert \mathbf{u}_{n}\lVert_{2}^{2} +\|\h_{n}-\mathbf{u}_{n}\|_{2}^{2}}= \lim_{n\rightarrow +\infty} \frac{N(\Q_{n}, \mathbf{u}_{n})}{\lVert \mathbf{u}_{n}\lVert_{2}^{2} }
\end{equation}
with $(\mathbf{u}_{n})_{n}$ a sequence of vectors of functions in the finite dimensional space spanned by $f_{1}^{\star},\ldots,f_{K}^{\star}$ and such that for $i=1,\ldots,K$,
\begin{equation}
\label{eq:zero}
\lim_{n\rightarrow +\infty}\int \frac{(\mathbf{u}_{n})_{i}(y)}{\lVert \mathbf{u}_{n}\lVert_{2} }dy=0
\end{equation}
since for all $n$ and all $i=1,\ldots,K$, $\int \frac{(\h_{n}-\mathbf{u}_{n})_{i}+(\mathbf{u}_{n})_{i}(y)}{\|\h_{n}\lVert_{2} }dy=0$.\\
Let us return to general considerations on the function $N(\cdot,\cdot)$.
As it may be seen from its formula, $N({\tilde{\Q}},{\tilde{\h}})$ is polynomial in the variables $\tilde{\Q}_{i,j}$,  $\langle f_{i}^{\star},f_{j}^{\star}\rangle$, $\langle \tilde{h}_{i},f_{j}^{\star}\rangle$, $\langle \tilde{h}_{i},\tilde{h}_{j}\rangle$, $i,j=1,\ldots,K$. 
Let $D({\tilde{\Q}},{\tilde{\h}})$ denote the part of $N({\tilde{\Q}},{\tilde{\h}})$ which is homogeneous of degree $2$ with respect to the variable ${\tilde{\h}}$, that is
\begin{eqnarray}
D({\tilde{\Q}},{\tilde{\h}})&=&\sum_{k_{1},k_{2},k_{3},k'_{1},k'_{2},k'_{3}=1}^{K}\tilde{\pi}(k_1){ \tilde{\Q}}(k_1,k_2){\tilde{\Q}}(k_{2},k_3)\tilde{\pi}(k'_1) {\tilde{\Q}}(k'_1,k'_2){\tilde{\Q}}(k'_{2},k'_3)\nonumber\\
&&\left(\sum_{i=1}^{3}\langle {\tilde{h}}_{k_{i}},{\tilde{h}}_{k'_{i}}\rangle \prod_{j\neq i}\langle f^{\star}_{k_{j}},f^{\star}_{k'_{j}}\rangle
+2\sum_{i=1}^{3}\langle f^{\star}_{k_{i}},f^{\star}_{k'_{i}}\rangle
\prod_{j\neq j' \neq i}\langle f^{\star}_{k_{j}},{\tilde{h}}_{k'_{j}}\rangle\langle {\tilde{h}}_{k_{j'}},f^{\star}_{k'_{j'}}\rangle
 \right). \label{eq:D}
\end{eqnarray}
One gets
\[
N({\tilde{\Q}},{\tilde{\h}})= D({\tilde{\Q}},{\tilde{\h}})+O\left(\lVert {\tilde{\h}}\lVert_{2}^{3}\right)
\]
where the $O(\cdot)$ depends only on $\f^{\star}$.
Let us first notice that $D(\cdot,\cdot)$ is always non negative. Indeed,
since for all ${\tilde{\Q}}\in {\mathcal V}$ and all ${\tilde{\h}}\in (\mathbf{L}^{2}(\Y,\L))^{K}$ one has $N({\tilde{\Q}},{\tilde{\h}})\geq 0$, it holds
\[
\forall {\tilde{\Q}}\in {\mathcal V}, \;\forall {\tilde{\h}}\in (\mathbf{L}^{2}(\Y,\L))^{K},\;\frac{D({\tilde{\Q}},{\tilde{\h}})}{\lVert {\tilde{\h}\lVert_{2}^{2}}}+O(\lVert {\tilde{\h}}\lVert_{2})\geq 0, 
\]
so that, since for all $\lambda\in\R$, $D({\tilde{\Q}},\lambda{\tilde{\h}})=\lambda^{2}D({\tilde{\Q}},{\tilde{\h}})$,
\begin{equation}
\label{Dpositif}
\forall{\tilde{\Q}}\in {\mathcal V}, \;\forall {\tilde{\h}}\in (\mathbf{L}^{2}(\Y,\L))^{K},\;D({\tilde{\Q}},{\tilde{\h}})\geq 0.
\end{equation}
Then we obtain from (\ref{eq:presque})
\[ 
c\geq \liminf_{n\rightarrow +\infty} D\big(\Q_{n},\frac{\mathbf{u_{n}}}{\lVert \mathbf{u_{n}}\lVert_{2}}\big).
\]
Let  ${\mathbf{b}}=(b_{1},\ldots, b_{K})$ be a limit point of the sequence $(\frac{\mathbf{u_{n}}}{\lVert \mathbf{u_{n}}\lVert_{2}})_{n}$. 
We then have
$$
c\geq D\left(\Q,{\mathbf{b}}\right).
$$
Now,
using (\ref{eq:zero}) we get that
$$
\int b_{k}d\L=0,\;k=1,\ldots,K.
$$

\noindent
Thus there exists a $K\times K$ matrix $U$ such that $\mathbf{b}^{T}=U(\f^{\star})^{T}$ and $U{\mathbf{1}}_{k}=0$, and equation (\ref{eq:D}) leads to
\begin{multline*}
D(\Q,\mathbf{b})=\\
 \sum_{i,j}\left\{ \left( \Q^{T}A UG^{\star} U^{T} A\Q\right)_{i,j}\left( G^{\star} \right)_{i,j}\left( \Q G^{\star}\Q^{T}\right)_{i,j}+
  \left( \Q^{T}A G^{\star} A\Q\right)_{i,j}\left(  UG^{\star}U^{T} \right)_{i,j}\left( \Q G^{\star}\Q^{T}\right)_{i,j}\right.\\
\left.+ \left( \Q^{T}AG^{\star}A\Q\right)_{i,j}\left( G^{\star} \right)_{i,j}\left( \Q UG^{\star}U^{T} \Q^{T}\right)_{i,j}\right\}
+ 2\sum_{i,j}\left\{ \left( \Q^{T}A UG^{\star} A\Q\right)_{i,j}\left( UG^{\star}\right)_{j,i}\left( \Q G^{\star}\Q^{T}\right)_{i,j}\right.\\+
  \left( \Q^{T}A UG^{\star} A\Q\right)_{i,j}\left(  \Q UG^{\star}\Q^{T} \right)_{j,i}\left( G^{\star}\right)_{i,j}
\left.+ \left( UG^{\star}\right)_{i,j}\left( \Q UG^{\star} \Q^{T}\right)_{j,i}\left( \Q^{T}AG^{\star}A\Q\right)_{i,j}\right\}.
\end{multline*}
with $G^{\star}$ the  $K\times K$ Gram matrix such that $(G^{\star})_{i,j}=\langle f^{\star}_{i},f^{\star}_{j}\rangle$, $i=1,\ldots,K$.
\\
This is the quadratic form  ${\mathcal D}$  in $U_{i,j}$, $i=1,\ldots,K$, $j=1,\ldots,K-1$ defined in Section \ref{sec:MAIN}. This quadratic form is non negative, and as soon as it is positive, we get that $c>0$. But the quadratic form ${\mathcal D}$ is positive as soon as its determinant is non zero, that is if and only  if $H(\Q,G(\f^{\star}))\neq 0$.

\subsection{Proof of Lemma \ref{LEMMA:INEQ3}}
\label{sec:bruteforce}
Here we specialize to the situation where $K=2$. In such a case, $\f^{\star}=(f_{1}^{\star},f_{2}^{\star})$, and 
\[
\Q^{\star}=\left(\begin{array}{cc} 1-p^{\star} & p^{\star} \\ q^{\star} & 1-q^{\star}   \end{array}\right)
\]
for some $p^{\star}, q^{\star}$ in $[0,1]$ for which $0<p^{\star}<1,\;0<q^{\star}<1,\;
p^{\star}\neq 1-q^{\star}.
$
Now 
\[
U=\left(\begin{array}{cc} \alpha& -\alpha \\ \beta& -\beta   \end{array}\right)
\] 
for some real numbers $\alpha$ and $\beta$, and
brute force computation gives $D\big(\Q,\mathbf{b}\big)=D_{1,1}\alpha^{2} + 2 D_{1,2}\alpha \beta + D_{2,2}\beta^{2}$ with, denoting $p=\Q(1,2)$ and $q=\Q(2,1)$:
\begin{align*}
\frac{(p+q)^{2}D_{1,1}}{q^{2} }=& 2(1-p)^{2}\lVert f^{\star}_{1}-f^{\star}_{2}\lVert ^{2}\lVert f^{\star}_{1}\lVert ^{2}\lVert (1-p)f^{\star}_{1}+pf^{\star}_{2}\lVert ^{2}+\lVert (1-p)f^{\star}_{1}+pf^{\star}_{2}\lVert ^{4}\lVert f^{\star}_{1}-f^{\star}_{2}\lVert ^{2}\\
&+4p(1-p)\left(\langle(1-p) f^{\star}_{1} +pf^{\star}_{2}, q f^{\star}_{1}+(1-q)f^{\star}_{2}\rangle\right)\langle f^{\star}_{1},f^{\star}_{2}\rangle \lVert f^{\star}_{1}-f^{\star}_{2}\lVert ^{2}
\\
&+2p^{2}\lVert f^{\star}_{1}-f^{\star}_{2}\lVert ^{2}\lVert f^{\star}_{2}\lVert ^{2}\lVert qf^{\star}_{1}+(1-q)f^{\star}_{2}\lVert ^{2}\\
&+2(1-p)^{2} \left(\langle (1-p)f^{\star}_{1}+pf^{\star}_{2},f^{\star}_{1}-f^{\star}_{2}\rangle\right)^{2}  \lVert  f^{\star}_{1}\lVert ^{2}\\
&+2p^{2}\left( \langle q f^{\star}_{1}+(1-q) f^{\star}_{2},f^{\star}_{1}-f^{\star}_{2}\rangle\right)^{2} \lVert f^{\star}_{2}\lVert ^{2}\\
&+4 p(1-p) \langle q f^{\star}_{1}+(1-q) f^{\star}_{2},f^{\star}_{1}-f^{\star}_{2}\rangle \langle (1-p)f^{\star}_{1}+pf^{\star}_{2},f^{\star}_{1}-f^{\star}_{2}\rangle \langle f^{\star}_{1},f^{\star}_{2}\rangle\\
&+4 (1-p)\langle (1-p)f^{\star}_{1}+pf^{\star}_{2},f^{\star}_{1}-f^{\star}_{2}\rangle \langle f^{\star}_{1},f^{\star}_{1}-f^{\star}_{2}\rangle \lVert  (1-p)f^{\star}_{1}+pf^{\star}_{2}\lVert ^{2} \\
&+4 p\langle (1-p)f^{\star}_{1}+pf^{\star}_{2},f^{\star}_{1}-f^{\star}_{2}\rangle  \langle f^{\star}_{2},f^{\star}_{1}-f^{\star}_{2}\rangle \langle(1-p) f^{\star}_{1} +pf^{\star}_{2}, q f^{\star}_{1}+(1-q)f^{\star}_{2}\rangle,\\
\end{align*}
\begin{align*}
\frac{(p+q)^{2}D_{2,2}}{p^{2}}=&2q^{2}\lVert f^{\star}_{1}-f^{\star}_{2}\lVert ^{2}\lVert f^{\star}_{1}\lVert ^{2}\lVert (1-p)f^{\star}_{1}+pf^{\star}_{2}\lVert ^{2}+\lVert qf^{\star}_{1}+(1-q)f^{\star}_{2}\lVert ^{4}\lVert f^{\star}_{1}-f^{\star}_{2}\lVert ^{2}\\
&+4(1-q)q\left(\langle(1-p) f^{\star}_{1} +pf^{\star}_{2}, q f^{\star}_{1}+(1-q)f^{\star}_{2}\rangle\right)\langle f^{\star}_{1},f^{\star}_{2}\rangle \lVert f^{\star}_{1}-f^{\star}_{2}\lVert ^{2}\\
&+2(1-q)^{2}\lVert f^{\star}_{1}-f^{\star}_{2}\lVert ^{2}\lVert f^{\star}_{2}\lVert ^{2}\lVert qf^{\star}_{1}+(1-q)f^{\star}_{2}\lVert ^{2}\\
&+2q^{2}  \left(\langle (1-p)f^{\star}_{1}+pf^{\star}_{2},f^{\star}_{1}-f^{\star}_{2}\rangle\right)^{2}  \lVert  f^{\star}_{1}\lVert ^{2}\\
&+2(1-q)^{2}\left( \langle q f^{\star}_{1}+(1-q) f^{\star}_{2},f^{\star}_{1}-f^{\star}_{2}\rangle\right)^{2} \lVert f_{2}\lVert ^{2}\\
&+4 q(1-q)\langle q f^{\star}_{1}+(1-q) f^{\star}_{2},f^{\star}_{1}-f^{\star}_{2}\rangle \langle (1-p)f^{\star}_{1}+pf^{\star}_{2},f^{\star}_{1}-f^{\star}_{2}\rangle \langle f^{\star}_{1},f^{\star}_{2}\rangle\\
&+4 q  \langle q f^{\star}_{1}+(1-q) f^{\star}_{2},f^{\star}_{1}-f^{\star}_{2}\rangle  \langle f^{\star}_{1},f^{\star}_{1}-f^{\star}_{2}\rangle \langle(1-p) f^{\star}_{1} +pf^{\star}_{2}, q f^{\star}_{1}+(1-q)f^{\star}_{2}\rangle \\
&+4 (1-q)  \langle q f^{\star}_{1}+(1-q) f^{\star}_{2},f^{\star}_{1}-f^{\star}_{2}\rangle  \langle f^{\star}_{2},f^{\star}_{1}-f^{\star}_{2}\rangle \lVert qf^{\star}_{1}+(1-q) f^{\star}_{2}\lVert ^{2}\,,
\end{align*}
and:
\begin{align*}
\frac{(p+q)^{2}D_{1,2}}{pq}=&2(1-p)q\lVert f^{\star}_{1}-f^{\star}_{2}\lVert ^{2}\lVert f^{\star}_{1}\lVert ^{2}\lVert (1-p)f^{\star}_{1}+pf^{\star}_{2}\lVert ^{2}\\
&+2[pq+(1-p)(1-q)]\left(\langle(1-p) f^{\star}_{1} +pf^{\star}_{2}, q f^{\star}_{1}+(1-q)f^{\star}_{2}\rangle\right)\langle f^{\star}_{1},f^{\star}_{2}\rangle \lVert f^{\star}_{1}-f^{\star}_{2}\lVert ^{2}\\
&+\left(\langle(1-p) f^{\star}_{1} +pf^{\star}_{2}, q f^{\star}_{1}+(1-q)f^{\star}_{2}\rangle\right)^{2} \lVert f^{\star}_{1}-f^{\star}_{2}\lVert ^{2}\\
&+2p(1-q)\lVert f^{\star}_{1}-f^{\star}_{2}\lVert ^{2}\lVert f^{\star}_{2}\lVert ^{2}\lVert qf^{\star}_{1}+(1-q)f^{\star}_{2}\lVert ^{2}\\
&+2q(1-p) \left(\langle (1-p)f^{\star}_{1}+pf^{\star}_{2},f^{\star}_{1}-f^{\star}_{2}\rangle\right)^{2}  \lVert  f^{\star}_{1}\lVert ^{2}\\
&+2p(1-q)\left( \langle q f^{\star}_{1}+(1-q) f^{\star}_{2},f^{\star}_{1}-f^{\star}_{2}\rangle\right)^{2} \lVert f^{\star}_{2}\lVert ^{2}\\
&+ 2pq \langle q f^{\star}_{1}+(1-q) f^{\star}_{2},f^{\star}_{1}-f^{\star}_{2}\rangle \langle (1-p)f^{\star}_{1}+pf^{\star}_{2},f_{1}-f^{\star}_{2}\rangle \langle f^{\star}_{1},f^{\star}_{2}\rangle\\
&+ 2(1-p)  (1-q)\langle q f^{\star}_{1}+(1-q) f^{\star}_{2},f_{1}-f^{\star}_{2}\rangle \langle (1-p)f^{\star}_{1}+pf^{\star}_{2},f^{\star}_{1}-f^{\star}_{2}\rangle \langle f^{\star}_{1},f^{\star}_{2}\rangle\\
&+q\langle (1-p)f^{\star}_{1}+pf^{\star}_{2},f^{\star}_{1}-f^{\star}_{2}\rangle \langle f^{\star}_{1},f^{\star}_{1}-f^{\star}_{2}\rangle \lVert  (1-p)f^{\star}_{1}+pf^{\star}_{2}\lVert ^{2} \\
&+2 (1-p)  \langle q f^{\star}_{1}+(1-q) f^{\star}_{2},f^{\star}_{1}-f^{\star}_{2}\rangle  \langle f^{\star}_{1},f^{\star}_{1}-f^{\star}_{2}\rangle \langle(1-p) f^{\star}_{1} +pf^{\star}_{2}, q f^{\star}_{1}+(1-q)f^{\star}_{2}\rangle \\
&+2  (1-q)\langle (1-p)f^{\star}_{1}+pf^{\star}_{2},f^{\star}_{1}-f^{\star}_{2}\rangle  \langle f^{\star}_{2},f^{\star}_{1}-f^{\star}_{2}\rangle \langle(1-p) f^{\star}_{1} +pf^{\star}_{2}, q f^{\star}_{1}+(1-q)f^{\star}_{2}\rangle\\
&+ 2p\langle q f^{\star}_{1}+(1-q) f^{\star}_{2},f^{\star}_{1}-f^{\star}_{2}\rangle  \langle f^{\star}_{2},f^{\star}_{1}-f^{\star}_{2}\rangle \lVert qf^{\star}_{1}+(1-q) f^{\star}_{2}\lVert ^{2}.
\end{align*}
We have:
\[
H(\Q,G(\f^{\star}))=D_{1,1}D_{2,2}-D_{1,2}^{2}.
\,
\]
We shall now write  $H(\Q,G(\f^{\star}))$ using
\[ n_{1}=\lVert f^{\star}_{1}\lVert_{2},\;n_{2}=\lVert f^{\star}_{2}\lVert_{2},\;a=\frac{\langle f^{\star}_{1} ,f^{\star}_{2} \rangle}{\lVert f^{\star}_{1}\lVert_{2}\lVert f^{\star}_{2}\lVert_{2}},
\]
for which the range is $[1,\infty[^{2}\times [0,1[$. Doing so, we obtain a polynomial $P_{1}$ in the variables $n_{1}$, $n_{2}$, $a$, $p$ and $q$.

\noindent
First observe that, by symmetry,
\[
P_{1}\left(n_{1},n_{2},a,p,q\right)=P_{1}\left(n_{2},n_{1},a,q,p\right),\,
\]
so that it is sufficient to prove that the polynomial $P_{1}$ is positive on the domain 
\eq
\label{eq:symElisabeth}
1\leq n_{2}\leq n_{1}\,,
\qe 
and $0\leq a<1$ and $0<p\neq q<1$.

\noindent
Furthermore, consider the change of variable 
\[
q=1-p+d
\]
then we have a polynomial $P_{2}$ in the variables $n_{1}$, $n_{2}$, $a$, $p$ and $d$ which factorizes with 
\[
\frac{p^2(1 - a^2) d^2 n_{1}^2 n_{2}^2 (1 + d - p)^2}{(1 + d)^4}\,.
\] 
Dividing by this factor, one gets a polynomial $P_{3}$ which is homogeneous of degree $8$ in $n_{1}$ and $n_{2}$, so that one may set $n_{1}=1$ and keep $b=n_{2}\in ]0,1]$ (observe that we have used \eqref{eq:symElisabeth} to reduce the problem to the domain $n_{2}/n_{1}\leq 1$) and obtain a polynomial $P_{4}$ in the variables $b$, $a$, $p$ and $d$. It remains to prove that $P_{4}$ is positive on $\mathcal D_{4}=\{b\in]0,1],a\in[0,1[, p\in]0,1[, d\in]p-1,0[\cup]0,p[\}$.

\noindent
Consider now the following change of variables
\[ 
b =\frac1{1 + x^2}\,,\quad  a =\frac{y^2}{1 + y^2}\,,\quad  p=\frac{z^2}{1 + z^2}\,,\quad\mathrm{and}\quad d =\frac{(tz)^2 - 1}{(1 + t^2)(1 + z^2)}\,,
\]
mapping $(x,y,z,t)\in\R^{4}$ onto $(b,a,p,d)\in\mathcal D_{5}=\{b\in]0,1],a\in[0,1[, p\in[0,1[, d\in]p-1,p[\}$ which contains $\mathcal D_{4}$. This change of variables maps $P_{4}$ onto a rational fraction with positive denominator, namely 
\[
(1 + t^2)^4 (1 + y^2)^4 (1 + z^2)^4 (1 + x^2)^8
\]
So it remains to prove that its numerator $P_{5}$, which is polynomial, is positive on $\R^{4}$. An expression of $P_{5}$ can be found in Appendix \ref{sec:P5}.
\noindent
Observe that $P_{5}$ is polynomial in $x^{2},y^{2},z^{2}$ and $t^{2}$ and there are only three monomials with negative coefficients. These monomials can be expressed as sum of squares using others monomials, namely:
\begin{itemize}
\item $ - 18x^{12}t^2 +27x^{12}+ 1979x^{12}t^4 = 18x^{12} + 
   9(x^6 - x^6t^2)^2 + 1970x^{12}t^4$,
\item $ - 108x^{10}t^2 + 1970x^{12}t^4+495x^8= 439x^8 + 56(x^4 - x^6t^2)^2 + 1914x^{12}t^4 + 4 t^2 x^{10}$,
\item and $- 114x^8t^2 +972x^4 + 1914x^{12}t^4 =915x^4 + 
   57(x^2 - x^6t^2)^2 + 1857x^{12}t^4$.
\end{itemize}
Thus $P_{5}$ is equal to $144$ more a sum of squares, hence it is positive. This proves that $H(\Q,G(\f^{\star}))$ is always positive.

\subsection{Proof of Theorem~\ref{theo:adapt}}
\label{sec:proof-theoadapt}

Let ${\mathcal K}=\{\h=\f-\f^{\star},\f\in{\mathcal F}^{K}\}$. 
Using  Theorem \ref{prop:oracle} 
we get that for all $x>0$, for all $N\geq N_{0}$,
with probability $1-(e-1)^{-1}e^{-x}$, one has for any permutation ${\tau}_{N}$, 
\begin{eqnarray}
\label{eq:oracle}
\|\hat{g}-g^{\star}\|_{2}^2
&\leq& 6 \inf_{M}\left\{\|g^{\star}-g^{\Q^{\star},{\f}^{\star}_{M}}\|_{2}^2+\pen(N,M)\right\}+ A^{\star}_{1}\frac{x}{N}\\
\nonumber &&+18 C_{\mathcal F,2}^6\left(\|\Q^\star-\mathbb P_{{\tau}_{N}}\hat\Q_N\mathbb P_{{\tau}_{N}}^{\top}\|_{F}^2+\|\pi^\star-\mathbb P_{{\tau}_{N}}\hat\pi\|_{2}^2\right).
\end{eqnarray}
Notice that writing
\begin{align*}
\hat{g}\left(y_{1},y_{2},y_{3}\right)=\sum_{k_1,k_{2},k_3=1}^K 
&(\mathbb P_{{\tau}_{N}}\hat\pi)(k_1) (\mathbb P_{{\tau}_{N}}\hat\Q \mathbb P_{{\tau}_{N}}^{\top})(k_1,k_2)(\mathbb P_{{\tau}_{N}}\hat\Q \mathbb P_{{\tau}_{N}}^{\top})(k_{2},k_3)\\
&\times\hat f_{{\tau}_{N}(k_1)}(y_1) \hat f_{{\tau}_{N}(k_2)}(y_2) \hat f_{{\tau}_{N}(k_3)}(y_3)\,,
\end{align*}
and applying Theorem \ref{LEMMA:INEQ2} we get that, on the event  $\mathbb P_{{\tau}_{N}}\hat\Q \mathbb P_{{\tau}_{N}}^{\top}\in \mathcal V$, there exists $\tau\in T_{\Q^{\star}}$ such that
\begin{equation}
\label{eq:step1}
\sum_{k=1}^{K}\|f^{\star}_{\tau(k)}-\hat{f}_{\tau_{N}(k)}\|_{2}^{2} 
\leq  \frac{1}{c({\mathcal K}, {\mathcal V}, \F^{\star})^{2}}\|\hat{g} - g^{\mathbb P_{\tau_{N}}\hat\Q \mathbb P_{\tau_{N}}^{\top},\f^{\star}}\|_{2}^{2}.
\end{equation}
Now by the triangular inequality
\begin{equation}
\label{eq:step3}
\|\hat{g} - g^{\mathbb P_{\tau_{N}}\hat\Q \mathbb P_{\tau_{N}}^{\top},\f^{\star}}\|_{2}\leq \|\hat{g}-g^{\star}\|_{2} + \|g^{\Q^{\star},\f^{\star}} - g^{\mathbb P_{\tau_{N}}\hat\Q \mathbb P_{\tau_{N}}^{\top},\f^{\star}}\|_{2}.
\end{equation}
Similarly to \eqref{qQ1-gQ2}, we have

\begin{equation}
\label{eq:step4}
\|g^{\Q^{\star},\f^{\star}} - g^{\mathbb P_{\tau_{N}}\hat\Q \mathbb P_{\tau_{N}}^{\top},\f^{\star}}\|_{2}^2
\leq 3K^{3}C_{\mathcal F,2}^{6}\left[\|\pi^\star-\mathbb P_{\tau_{N}}\hat\pi\|_{2}^2+2\|\Q^\star-\mathbb P_{\tau_{N}}\hat\Q\mathbb P_{\tau_{N}}^{\top}\|_{F}^2 \right].
\end{equation}
In the same way,
\begin{multline*}
\left(g^{\star}-g^{\Q^{\star},\f^{\star}_{M}}\right)\left(y_{1},y_{2},y_{3}\right)=\\
\sum_{k_1,k_{2},k_3=1}^K
\pi^{\star}(k_1) \Q^{\star}(k_1,k_2)\Q^{\star}(k_{2},k_3)\left(f^{\star}_{k_1}(y_1) f^{\star}_{k_2}(y_2) f^{\star}_{k_3}(y_3)-f^{\star}_{M,k_1}(y_1) f^{\star}_{M,k_2}(y_2) f^{\star}_{M,k_3}(y_3)
\right)
\end{multline*}
so that
\[
\|g^{\star}-g^{\Q^{\star},\f^{\star}_{M}}\|_{2}^2\leq {3}K^{3}C_{\mathcal F,2}^{4}\max\{\|f^{\star}_{k}- f^{\star}_{M,k}\|_{2}^2,\;k=1,\ldots,K\}.
\]
Thus collecting \eqref{eq:oracle},
 \eqref{eq:step1}, 
 \eqref{eq:step3}, \eqref{eq:step4} and with an appropriate choice of $A^{\star}$ we get Theorem~\ref{theo:adapt}.

\subsection{Proof of Corollary~\ref{cor:adapt}}
\label{sec:proof-coradapt}
We shall apply Theorem \ref{thm:spectral} where, for each $N$, we define $\delta_{N}$ such that $(-\log \delta_{N})/\delta_{N}^{2}:=(\log N)^{1/2}$. Notice first that $\delta_{N}$ goes to $0$ and that $M_{N}$ tends to infinity as $N$ tends to infinity, so that for large enough $N$, $M_{N}\geq M_{\F^{\star}}$.  By denoting $\tau_{N}$ the $\tau_{M_{N}}$ given by Theorem \ref{thm:spectral} we get that  for all $x\geq x(\Q^{\star})$, for all $N\geq{\mathbf N}(\Q^{\star},\F^{\star})x \log N$, with probability $1-[4+(e-1)^{-1}]e^{-x}-2\delta_{N}$,
$$\lVert\pi^{\star}-\mathbb P_{\tau_{N}}\hat\pi\lVert_{2}
\leq
 {\mathcal C}(\Q^{\star},\F^{\star})\sqrt{\frac{\log N}{N}}\sqrt{x}
 $$
and
$$
\lVert\Q^{\star}-\mathbb P_{\tau_{N}}\hat\Q \mathbb P_{\tau_{N}}^{\top}\lVert
\leq
 {\mathcal C}(\Q^{\star},\F^{\star})\sqrt{\frac{\log N}{N}}\sqrt{x}\,.
$$
We first obtain that
\begin{multline*}
\limsup_{N\rightarrow +\infty}\E\left[ \frac{ N}{\log N}\lVert \Q^{\star}-\mathbb P_{\tau_{N}}\hat\Q\mathbb P_{\tau_{N}}^{T}\lVert^{2}\right] \leq \\
{\mathcal C}(\Q^{\star},\F^{\star})^{2}\int_{0}^{+\infty}\limsup_{N\rightarrow +\infty}\P\left( \frac{\sqrt N}{{\mathcal C}(\Q^{\star},\F^{\star})\sqrt{\log N} }\lVert \Q^{\star}-\mathbb P_{\tau_{N}}\hat\Q\mathbb P_{\tau_{N}}^{T}\lVert\geq  \sqrt{x}\right)dx 
\leq 
\\
{\mathcal C}(\Q^{\star},\F^{\star})^{2} x(\Q^{\star})+{\mathcal C}(\Q^{\star},\F^{\star})^{2}\int_{x(\Q^{\star})}^{+\infty}[4+(e-1)^{-1}]e^{-x}dx <+\infty
\end{multline*}
so that
\[
\E\left[\lVert \Q^{\star}-\mathbb P_{\tau_{N}}\hat\Q\mathbb P_{\tau_{N}}^{T}\lVert^{2}\right]=O\left(\frac{\log N}{N}\right).
\]
Similarly, one has $\E\left[ \lVert \pi^{\star}-\mathbb P_{\tau_{N}}\hat\pi\lVert^{2}\right]=O\left(\frac{\log N}{N}\right)$. We also obtain, by taking $x=N/(\log N)^{1/4}$, that
$$
\limsup_{N\rightarrow +\infty}\P\left(  \mathbb P_{\tau_{N}}\hat\Q\mathbb P_{\tau_{N}}^{T} \notin \mathcal V\right)  =0,
$$
so that, using Theorem \ref{theo:adapt}, we get for some $\tau\in T_{\Q^{\star}}$,
 \begin{multline*}
\limsup_{N\rightarrow +\infty}\P\left(
\frac{N}{A^{\star}}\left[\sum_{k=1}^{K}\|f^{\star}_{k}-\hat{f}_{\tau^{-1}\circ\tau_{N}(k)}\|_{2}^{2}-\inf_{M}\left\{\sum_{k=1}^{K}\|f^{\star}_{k}- f^{\star}_{M,k}\|_2^2+\pen(N,M)\right\}\right.\right.\\
 \left.\left.-\lVert\Q^{\star}-\mathbb P_{\tau_{N}}\hat\Q \mathbb P_{\tau_{N}}^{\top}\lVert_{F}^{2}
 -\lVert\pi^{\star}-\mathbb P_{\tau_{N}}\hat\pi\lVert_{2}^{2} +\frac{x}{N}\right]
\geq x\right) \leq (e-1)^{-1}e^{-x}.
\end{multline*}
Thus,  by integration and the previous results, Corollary \ref{cor:adapt} follows. 



\acks{We would like to thank N. Curien and D. Henrion for their help in the computation of $P_{5}$. We are deeply grateful to V. Magron for pointing us that $P_{5}$ can be explicitly expressed as a Sum-Of-Squares. Thanks are due to Luc Lehericy for a careful reading of the manuscript. We are grateful to anonymous referees for their careful reading of this work.
 }


\newpage

\appendix
\section{Concentration inequalities}
\label{sec:AppendixConcentration}



We first recall  results that hold both for  {\bf (Scenario A)} (where we  consider $N$ i.i.d. samples $(Y_{1}^{(s)},Y_{2}^{(s)},Y_{3}^{(s)})_{s=1}^{N}$ of three consecutive observations) and for {\bf (Scenario B)} (where we consider consecutive observations of the same chain).

The following proposition is the classical Bernstein's inequality for  {\bf (Scenario A)} and is proved in  \cite{paulin2012concentrationv2}, Theorem 2.4, 
for  {\bf (Scenario B)}. 
\begin{proposition} \label{prop:bernstein}
Let $t$ be a real valued and measurable bounded function on $\Y^{3}$.  Let $V=\E [t^{2}(Z_{1})]$. There exists a positive constant $c^{\star}$ depending only on $\Q^{\star}$ such that for all $0\leq \lambda \leq 1/(2\sqrt{2}c^{\star}\lVert t\lVert_{\infty})$ :
\begin{equation}
\label{eq:loglaplace}
\log \E \exp \left[ \lambda \sum_{s=1}^{N}\left(t(Z_{s})-\E t(Z_{s})\right)  \right] \leq \frac{2Nc^{\star}V\lambda^{2}}{1-2\sqrt{2}c^{\star}\lVert t\lVert_{\infty}\lambda}
\end{equation}
so that for all $x\geq 0$,
\begin{equation}
\label{eq:bernstein}
\P\left(\sum_{s=1}^{N}\left(t(Z_{s})-\E t(Z_{s})\right)\geq 2\sqrt{2Nc^{\star}Vx} +2\sqrt{2}c^{\star}\lVert t\lVert_{\infty}x \right)\leq e^{-x}.
\end{equation}
\end{proposition}
We now state a deviation inequality, which comes from \cite{MR2319879}
Theorem 6.8 and Corollary 6.9 for {\bf (Scenario A)}. For {\bf (Scenario B)}
the proof of the following proposition follows \textit{mutatis mutandis} from the proof of Theorem 6.8 (and then Corollary 6.9) in \cite{MR2319879} the early first step being equation (\ref{eq:loglaplace}).
Recall that when $t_{1}$ and $t_{2}$ are real valued functions, the bracket $[t_{1},t_{2}]$ is the set of real valued functions $t$ such that $t_{1}(\cdot)\leq t(\cdot) \leq t_{2}(\cdot)$. For any measurable set $A$ such that $\P(A)>0$, and any integrable random variable $Z$, denote $E^{A}[Z]=E[Z\one_{A}]/\P(A)$.
\begin{proposition}
\label{prop:deviation}
 Let $\mathcal{T}$ be some countable class of real valued and measurable
functions on $\Y^{3}$. Assume that there exists some positive numbers $\sigma$ and $b$
such that for all $t\in \mathcal{T}$, $\lVert t\lVert_{\infty}\leq b$ and $\E [t^{2}(Z_{1})]\leq \sigma^{2}$.\\
Assume furthermore that for any positive number $\delta$, there exists some finite
set $B_\delta$ of brackets covering $\mathcal{F}$ such that for any bracket $[t_1, t_2]\in
B_\delta$, $\lVert t_{1}-t_{2}\lVert_{\infty}\leq b$ and $\E[(t_{1}-t_{2})^{2}(Z_{1})]\leq \delta^{2}$.
Let $e^{H(\delta)}$
denote the minimal cardinality of such a covering. Then, there exists a positive constant $C^{\star}$ depending only on $\Q^{\star}$ such that:
for any measurable set $A$, 
\[ \E^{A}\left(\sup_{t\in \mathcal{T}}\sum_{s=1}^{N}\left(t(Z_{s})-\E t(Z_{s}\right))\right)\leq C^{\star}\left[
E+\sigma\sqrt{N\log\left(\frac{1}{\P(A)}\right)}+b\log\left(\frac{1}{\P(A)}\right)\right]\]
and for  all positive number $x$
\[ \P\left(\sup_{t\in \mathcal{T}}\sum_{s=1}^{N}\left(t(Z_{s})-\E t(Z_{s}\right))\geq C^{\star}[E+\sigma\sqrt{Nx}+bx]\right) \leq\exp(-x),\]
where 
 \[ E=\sqrt{N}\int_0^{\sigma}\sqrt{H(u)\wedge N}du + (b+\sigma)H(\sigma).\]
\end{proposition}

\section{Expression of polynomial $P_{5}$}
\label{sec:P5}
\noindent
Computer assisted computations (available at \url{https://mycore.core-cloud.net/public.php?service=files&t=db7b8c1a2bcbcca157dcda5ecab84374}) give that:

\medskip

{
\small
\noindent
$P_{5}=$
\begin{verbatim}
 144 - 114 t^2 x^8 - 108 t^2 x^10  - 18 t^2 x^12 +
 192 t^2 + 128 t^4 + 256 t^6 + 176 t^8 + 576 x^2 + 624 t^2 x^2 + 
 672 t^4 x^2 + 1776 t^6 x^2 + 1152 t^8 x^2 + 972 x^4 + 720 t^2 x^4 + 
 1884 t^4 x^4 + 5496 t^6 x^4 + 3360 t^8 x^4 + 900 x^6 + 264 t^2 x^6 + 
 3556 t^4 x^6 + 9920 t^6 x^6 + 5728 t^8 x^6 + 495 x^8 + 
 4551 t^4 x^8 + 11424 t^6 x^8 + 6264 t^8 x^8 + 162 x^10 + 
 3810 t^4 x^10 + 8592 t^6 x^10 + 4512 t^8 x^10 + 
 27 x^12 + 1979 t^4 x^12 + 4120 t^6 x^12 + 
 2096 t^8 x^12 + 576 t^4 x^14 + 1152 t^6 x^14 + 576 t^8 x^14 + 
 72 t^4 x^16 + 144 t^6 x^16 + 72 t^8 x^16 + 144 y^2 + 480 t^2 y^2 + 
 784 t^4 y^2 + 704 t^6 y^2 + 256 t^8 y^2 + 576 x^2 y^2 + 
 2064 t^2 x^2 y^2 + 4192 t^4 x^2 y^2 + 4496 t^6 x^2 y^2 + 
 1792 t^8 x^2 y^2 + 1080 x^4 y^2 + 4104 t^2 x^4 y^2 + 
 10760 t^4 x^4 y^2 + 13528 t^6 x^4 y^2 + 5792 t^8 x^4 y^2 + 
 1224 x^6 y^2 + 5016 t^2 x^6 y^2 + 17592 t^4 x^6 y^2 + 
 25032 t^6 x^6 y^2 + 11232 t^8 x^6 y^2 + 900 x^8 y^2 + 
 4224 t^2 x^8 y^2 + 19924 t^4 x^8 y^2 + 30776 t^6 x^8 y^2 + 
 14176 t^8 x^8 y^2 + 432 x^10 y^2 + 2520 t^2 x^10 y^2 + 
 15584 t^4 x^10 y^2 + 25336 t^6 x^10 y^2 + 11840 t^8 x^10 y^2 + 
 108 x^12 y^2 + 936 t^2 x^12 y^2 + 7916 t^4 x^12 y^2 + 
 13456 t^6 x^12 y^2 + 6368 t^8 x^12 y^2 + 144 t^2 x^14 y^2 + 
 2304 t^4 x^14 y^2 + 4176 t^6 x^14 y^2 + 2016 t^8 x^14 y^2 + 
 288 t^4 x^16 y^2 + 576 t^6 x^16 y^2 + 288 t^8 x^16 y^2 + 144 y^4 + 
 480 t^2 y^4 + 624 t^4 y^4 + 384 t^6 y^4 + 96 t^8 y^4 + 576 x^2 y^4 + 
 2208 t^2 x^2 y^4 + 3392 t^4 x^2 y^4 + 2464 t^6 x^2 y^4 + 
 704 t^8 x^2 y^4 + 1188 x^4 y^4 + 5256 t^2 x^4 y^4 + 
 9636 t^4 x^4 y^4 + 8256 t^6 x^4 y^4 + 2688 t^8 x^4 y^4 + 
 1548 x^6 y^4 + 8112 t^2 x^6 y^4 + 18076 t^4 x^6 y^4 + 
 18008 t^6 x^6 y^4 + 6496 t^8 x^6 y^4 + 1359 x^8 y^4 + 
 8598 t^2 x^8 y^4 + 23375 t^4 x^8 y^4 + 26392 t^6 x^8 y^4 + 
 10256 t^8 x^8 y^4 + 810 x^10 y^4 + 6156 t^2 x^10 y^4 + 
 20442 t^4 x^10 y^4 + 25656 t^6 x^10 y^4 + 10560 t^8 x^10 y^4 + 
 243 x^12 y^4 + 2574 t^2 x^12 y^4 + 11299 t^4 x^12 y^4 + 
 15848 t^6 x^12 y^4 + 6880 t^8 x^12 y^4 + 432 t^2 x^14 y^4 + 
 3456 t^4 x^14 y^4 + 5616 t^6 x^14 y^4 + 2592 t^8 x^14 y^4 + 
 432 t^4 x^16 y^4 + 864 t^6 x^16 y^4 + 432 t^8 x^16 y^4 + 
 216 x^4 y^6 + 720 t^2 x^4 y^6 + 952 t^4 x^4 y^6 + 608 t^6 x^4 y^6 + 
 160 t^8 x^4 y^6 + 648 x^6 y^6 + 2592 t^2 x^6 y^6 + 
 4168 t^4 x^6 y^6 + 3152 t^6 x^6 y^6 + 928 t^8 x^6 y^6 + 
 918 x^8 y^6 + 4428 t^2 x^8 y^6 + 8502 t^4 x^8 y^6 + 
 7392 t^6 x^8 y^6 + 2400 t^8 x^8 y^6 + 756 x^10 y^6 + 
 4392 t^2 x^10 y^6 + 10036 t^4 x^10 y^6 + 9920 t^6 x^10 y^6 + 
 3520 t^8 x^10 y^6 + 270 x^12 y^6 + 2268 t^2 x^12 y^6 + 
 6766 t^4 x^12 y^6 + 7808 t^6 x^12 y^6 + 3040 t^8 x^12 y^6 + 
 432 t^2 x^14 y^6 + 2304 t^4 x^14 y^6 + 3312 t^6 x^14 y^6 + 
 1440 t^8 x^14 y^6 + 288 t^4 x^16 y^6 + 576 t^6 x^16 y^6 + 
 288 t^8 x^16 y^6 + 108 x^8 y^8 + 360 t^2 x^8 y^8 + 468 t^4 x^8 y^8 + 
 288 t^6 x^8 y^8 + 72 t^8 x^8 y^8 + 216 x^10 y^8 + 864 t^2 x^10 y^8 + 
 1368 t^4 x^10 y^8 + 1008 t^6 x^10 y^8 + 288 t^8 x^10 y^8 + 
 108 x^12 y^8 + 648 t^2 x^12 y^8 + 1404 t^4 x^12 y^8 + 
 1296 t^6 x^12 y^8 + 432 t^8 x^12 y^8 + 144 t^2 x^14 y^8 + 
 576 t^4 x^14 y^8 + 720 t^6 x^14 y^8 + 288 t^8 x^14 y^8 + 
 72 t^4 x^16 y^8 + 144 t^6 x^16 y^8 + 72 t^8 x^16 y^8 + 192 z^2 + 
 416 t^2 z^2 + 288 t^4 z^2 + 320 t^6 z^2 + 256 t^8 z^2 + 
 912 x^2 z^2 + 1664 t^2 x^2 z^2 + 1248 t^4 x^2 z^2 + 
 2304 t^6 x^2 z^2 + 1808 t^8 x^2 z^2 + 1728 x^4 z^2 + 
 2520 t^2 x^4 z^2 + 2776 t^4 x^4 z^2 + 7624 t^6 x^4 z^2 + 
 5640 t^8 x^4 z^2 + 1704 x^6 z^2 + 1736 t^2 x^6 z^2 + 
 4664 t^4 x^6 z^2 + 14808 t^6 x^6 z^2 + 10176 t^8 x^6 z^2 + 
 966 x^8 z^2 + 494 t^2 x^8 z^2 + 6098 t^4 x^8 z^2 + 
 18218 t^6 x^8 z^2 + 11648 t^8 x^8 z^2 + 324 x^10 z^2 + 
 36 t^2 x^10 z^2 + 5468 t^4 x^10 z^2 + 14444 t^6 x^10 z^2 + 
 8688 t^8 x^10 z^2 + 54 x^12 z^2 + 6 t^2 x^12 z^2 + 
 3002 t^4 x^12 z^2 + 7186 t^6 x^12 z^2 + 4136 t^8 x^12 z^2 + 
 896 t^4 x^14 z^2 + 2048 t^6 x^14 z^2 + 1152 t^8 x^14 z^2 + 
 112 t^4 x^16 z^2 + 256 t^6 x^16 z^2 + 144 t^8 x^16 z^2 + 
 480 y^2 z^2 + 1312 t^2 y^2 z^2 + 1888 t^4 y^2 z^2 + 
 1760 t^6 y^2 z^2 + 704 t^8 y^2 z^2 + 1776 x^2 y^2 z^2 + 
 5248 t^2 x^2 y^2 z^2 + 9504 t^4 x^2 y^2 z^2 + 
 10624 t^6 x^2 y^2 z^2 + 4592 t^8 x^2 y^2 z^2 + 3096 x^4 y^2 z^2 + 
 9904 t^2 x^4 y^2 z^2 + 23104 t^4 x^4 y^2 z^2 + 
 30288 t^6 x^4 y^2 z^2 + 13992 t^8 x^4 y^2 z^2 + 3144 x^6 y^2 z^2 + 
 11344 t^2 x^6 y^2 z^2 + 35712 t^4 x^6 y^2 z^2 + 
 53424 t^6 x^6 y^2 z^2 + 25912 t^8 x^6 y^2 z^2 + 2064 x^8 y^2 z^2 + 
 9016 t^2 x^8 y^2 z^2 + 38552 t^4 x^8 y^2 z^2 + 
 63192 t^6 x^8 y^2 z^2 + 31592 t^8 x^8 y^2 z^2 + 936 x^10 y^2 z^2 + 
 5248 t^2 x^10 y^2 z^2 + 29072 t^4 x^10 y^2 z^2 + 
 50464 t^6 x^10 y^2 z^2 + 25704 t^8 x^10 y^2 z^2 + 216 x^12 y^2 z^2 + 
 1872 t^2 x^12 y^2 z^2 + 14192 t^4 x^12 y^2 z^2 + 
 26056 t^6 x^12 y^2 z^2 + 13520 t^8 x^12 y^2 z^2 + 
 264 t^2 x^14 y^2 z^2 + 3896 t^4 x^14 y^2 z^2 + 
 7808 t^6 x^14 y^2 z^2 + 4176 t^8 x^14 y^2 z^2 + 
 448 t^4 x^16 y^2 z^2 + 1024 t^6 x^16 y^2 z^2 + 
 576 t^8 x^16 y^2 z^2 + 480 y^4 z^2 + 1632 t^2 y^4 z^2 + 
 2208 t^4 y^4 z^2 + 1440 t^6 y^4 z^2 + 384 t^8 y^4 z^2 + 
 1632 x^2 y^4 z^2 + 6528 t^2 x^2 y^4 z^2 + 10688 t^4 x^2 y^4 z^2 + 
 8320 t^6 x^2 y^4 z^2 + 2528 t^8 x^2 y^4 z^2 + 3240 x^4 y^4 z^2 + 
 14280 t^2 x^4 y^4 z^2 + 27448 t^4 x^4 y^4 z^2 + 
 25048 t^6 x^4 y^4 z^2 + 8640 t^8 x^4 y^4 z^2 + 3936 x^6 y^4 z^2 + 
 19992 t^2 x^6 y^4 z^2 + 46552 t^4 x^6 y^4 z^2 + 
 49352 t^6 x^6 y^4 z^2 + 18856 t^8 x^6 y^4 z^2 + 3198 x^8 y^4 z^2 + 
 19518 t^2 x^8 y^4 z^2 + 55218 t^4 x^8 y^4 z^2 + 
 66170 t^6 x^8 y^4 z^2 + 27272 t^8 x^8 y^4 z^2 + 1836 x^10 y^4 z^2 + 
 13332 t^2 x^10 y^4 z^2 + 44988 t^4 x^10 y^4 z^2 + 
 59580 t^6 x^10 y^4 z^2 + 26088 t^8 x^10 y^4 z^2 + 486 x^12 y^4 z^2 + 
 5214 t^2 x^12 y^4 z^2 + 22994 t^4 x^12 y^4 z^2 + 
 34194 t^6 x^12 y^4 z^2 + 15928 t^8 x^12 y^4 z^2 + 
 792 t^2 x^14 y^4 z^2 + 6312 t^4 x^14 y^4 z^2 + 
 11136 t^6 x^14 y^4 z^2 + 5616 t^8 x^14 y^4 z^2 + 
 672 t^4 x^16 y^4 z^2 + 1536 t^6 x^16 y^4 z^2 + 
 864 t^8 x^16 y^4 z^2 + 720 x^4 y^6 z^2 + 2480 t^2 x^4 y^6 z^2 + 
 3472 t^4 x^4 y^6 z^2 + 2384 t^6 x^4 y^6 z^2 + 672 t^8 x^4 y^6 z^2 + 
 1728 x^6 y^6 z^2 + 7440 t^2 x^6 y^6 z^2 + 13072 t^4 x^6 y^6 z^2 + 
 10736 t^6 x^6 y^6 z^2 + 3376 t^8 x^6 y^6 z^2 + 2268 x^8 y^6 z^2 + 
 11484 t^2 x^8 y^6 z^2 + 23812 t^4 x^8 y^6 z^2 + 
 22276 t^6 x^8 y^6 z^2 + 7680 t^8 x^8 y^6 z^2 + 1800 x^10 y^6 z^2 + 
 10568 t^2 x^10 y^6 z^2 + 25560 t^4 x^10 y^6 z^2 + 
 26872 t^6 x^10 y^6 z^2 + 10080 t^8 x^10 y^6 z^2 + 540 x^12 y^6 z^2 + 
 4836 t^2 x^12 y^6 z^2 + 15420 t^4 x^12 y^6 z^2 + 
 18964 t^6 x^12 y^6 z^2 + 7840 t^8 x^12 y^6 z^2 + 
 792 t^2 x^14 y^6 z^2 + 4520 t^4 x^14 y^6 z^2 + 
 7040 t^6 x^14 y^6 z^2 + 3312 t^8 x^14 y^6 z^2 + 
 448 t^4 x^16 y^6 z^2 + 1024 t^6 x^16 y^6 z^2 + 
 576 t^8 x^16 y^6 z^2 + 360 x^8 y^8 z^2 + 1224 t^2 x^8 y^8 z^2 + 
 1656 t^4 x^8 y^8 z^2 + 1080 t^6 x^8 y^8 z^2 + 288 t^8 x^8 y^8 z^2 + 
 576 x^10 y^8 z^2 + 2448 t^2 x^10 y^8 z^2 + 4176 t^4 x^10 y^8 z^2 + 
 3312 t^6 x^10 y^8 z^2 + 1008 t^8 x^10 y^8 z^2 + 216 x^12 y^8 z^2 + 
 1488 t^2 x^12 y^8 z^2 + 3616 t^4 x^12 y^8 z^2 + 
 3640 t^6 x^12 y^8 z^2 + 1296 t^8 x^12 y^8 z^2 + 
 264 t^2 x^14 y^8 z^2 + 1208 t^4 x^14 y^8 z^2 + 
 1664 t^6 x^14 y^8 z^2 + 720 t^8 x^14 y^8 z^2 + 
 112 t^4 x^16 y^8 z^2 + 256 t^6 x^16 y^8 z^2 + 144 t^8 x^16 y^8 z^2 + 
 128 z^4 + 288 t^2 z^4 + 352 t^4 z^4 + 384 t^6 z^4 + 256 t^8 z^4 + 
 352 x^2 z^4 + 1056 t^2 x^2 z^4 + 1408 t^4 x^2 z^4 + 
 1952 t^6 x^2 z^4 + 1504 t^8 x^2 z^4 + 764 x^4 z^4 + 
 2104 t^2 x^4 z^4 + 2616 t^4 x^4 z^4 + 5016 t^6 x^4 z^4 + 
 4252 t^8 x^4 z^4 + 804 x^6 z^4 + 1912 t^2 x^6 z^4 + 
 2920 t^4 x^6 z^4 + 8536 t^6 x^6 z^4 + 7364 t^8 x^6 z^4 + 
 471 x^8 z^4 + 898 t^2 x^8 z^4 + 2694 t^4 x^8 z^4 + 
 10058 t^6 x^8 z^4 + 8335 t^8 x^8 z^4 + 162 x^10 z^4 + 
 252 t^2 x^10 z^4 + 2164 t^4 x^10 z^4 + 7980 t^6 x^10 z^4 + 
 6226 t^8 x^10 z^4 + 27 x^12 z^4 + 42 t^2 x^12 z^4 + 
 1182 t^4 x^12 z^4 + 4018 t^6 x^12 z^4 + 2979 t^8 x^12 z^4 + 
 352 t^4 x^14 z^4 + 1152 t^6 x^14 z^4 + 832 t^8 x^14 z^4 + 
 44 t^4 x^16 z^4 + 144 t^6 x^16 z^4 + 104 t^8 x^16 z^4 + 
 784 y^2 z^4 + 1888 t^2 y^2 z^4 + 2208 t^4 y^2 z^4 + 
 1888 t^6 y^2 z^4 + 784 t^8 y^2 z^4 + 2080 x^2 y^2 z^4 + 
 5600 t^2 x^2 y^2 z^4 + 8832 t^4 x^2 y^2 z^4 + 9952 t^6 x^2 y^2 z^4 + 
 4640 t^8 x^2 y^2 z^4 + 3368 x^4 y^2 z^4 + 9440 t^2 x^4 y^2 z^4 + 
 18928 t^4 x^4 y^2 z^4 + 25952 t^6 x^4 y^2 z^4 + 
 13224 t^8 x^4 y^2 z^4 + 2840 x^6 y^2 z^4 + 9056 t^2 x^6 y^2 z^4 + 
 25872 t^4 x^6 y^2 z^4 + 42464 t^6 x^6 y^2 z^4 + 
 23192 t^8 x^6 y^2 z^4 + 1524 x^8 y^2 z^4 + 6072 t^2 x^8 y^2 z^4 + 
 25016 t^4 x^8 y^2 z^4 + 46792 t^6 x^8 y^2 z^4 + 
 26900 t^8 x^8 y^2 z^4 + 576 x^10 y^2 z^4 + 3184 t^2 x^10 y^2 z^4 + 
 17216 t^4 x^10 y^2 z^4 + 35024 t^6 x^10 y^2 z^4 + 
 20928 t^8 x^10 y^2 z^4 + 108 x^12 y^2 z^4 + 1008 t^2 x^12 y^2 z^4 + 
 7584 t^4 x^12 y^2 z^4 + 16968 t^6 x^12 y^2 z^4 + 
 10572 t^8 x^12 y^2 z^4 + 120 t^2 x^14 y^2 z^4 + 
 1816 t^4 x^14 y^2 z^4 + 4736 t^6 x^14 y^2 z^4 + 
 3136 t^8 x^14 y^2 z^4 + 176 t^4 x^16 y^2 z^4 + 
 576 t^6 x^16 y^2 z^4 + 416 t^8 x^16 y^2 z^4 + 624 y^4 z^4 + 
 2208 t^2 y^4 z^4 + 3168 t^4 y^4 z^4 + 2208 t^6 y^4 z^4 + 
 624 t^8 y^4 z^4 + 1600 x^2 y^4 z^4 + 6976 t^2 x^2 y^4 z^4 + 
 12672 t^4 x^2 y^4 z^4 + 10816 t^6 x^2 y^4 z^4 + 
 3520 t^8 x^2 y^4 z^4 + 3364 x^4 y^4 z^4 + 14456 t^2 x^4 y^4 z^4 + 
 29416 t^4 x^4 y^4 z^4 + 29016 t^6 x^4 y^4 z^4 + 
 10692 t^8 x^4 y^4 z^4 + 3452 x^6 y^4 z^4 + 17336 t^2 x^6 y^4 z^4 + 
 43896 t^4 x^6 y^4 z^4 + 51032 t^6 x^6 y^4 z^4 + 
 21020 t^8 x^6 y^4 z^4 + 2495 x^8 y^4 z^4 + 14658 t^2 x^8 y^4 z^4 + 
 45814 t^4 x^8 y^4 z^4 + 61162 t^6 x^8 y^4 z^4 + 
 27607 t^8 x^8 y^4 z^4 + 1242 x^10 y^4 z^4 + 8892 t^2 x^10 y^4 z^4 + 
 33252 t^4 x^10 y^4 z^4 + 49644 t^6 x^10 y^4 z^4 + 
 24234 t^8 x^10 y^4 z^4 + 243 x^12 y^4 z^4 + 2914 t^2 x^12 y^4 z^4 + 
 14758 t^4 x^12 y^4 z^4 + 25538 t^6 x^12 y^4 z^4 + 
 13643 t^8 x^12 y^4 z^4 + 360 t^2 x^14 y^4 z^4 + 
 3336 t^4 x^14 y^4 z^4 + 7296 t^6 x^14 y^4 z^4 + 
 4416 t^8 x^14 y^4 z^4 + 264 t^4 x^16 y^4 z^4 + 
 864 t^6 x^16 y^4 z^4 + 624 t^8 x^16 y^4 z^4 + 952 x^4 y^6 z^4 + 
 3472 t^2 x^4 y^6 z^4 + 5232 t^4 x^4 y^6 z^4 + 3856 t^6 x^4 y^6 z^4 + 
 1144 t^8 x^4 y^6 z^4 + 1544 x^6 y^6 z^4 + 7760 t^2 x^6 y^6 z^4 + 
 15696 t^4 x^6 y^6 z^4 + 14288 t^6 x^6 y^6 z^4 + 
 4808 t^8 x^6 y^6 z^4 + 1942 x^8 y^6 z^4 + 10532 t^2 x^8 y^6 z^4 + 
 24556 t^4 x^8 y^6 z^4 + 25380 t^6 x^8 y^6 z^4 + 
 9414 t^8 x^8 y^6 z^4 + 1332 x^10 y^6 z^4 + 8408 t^2 x^10 y^6 z^4 + 
 22952 t^4 x^10 y^6 z^4 + 26776 t^6 x^10 y^6 z^4 + 
 10900 t^8 x^10 y^6 z^4 + 270 x^12 y^6 z^4 + 2972 t^2 x^12 y^6 z^4 + 
 11492 t^4 x^12 y^6 z^4 + 16244 t^6 x^12 y^6 z^4 + 
 7486 t^8 x^12 y^6 z^4 + 360 t^2 x^14 y^6 z^4 + 
 2632 t^4 x^14 y^6 z^4 + 4992 t^6 x^14 y^6 z^4 + 
 2752 t^8 x^14 y^6 z^4 + 176 t^4 x^16 y^6 z^4 + 
 576 t^6 x^16 y^6 z^4 + 416 t^8 x^16 y^6 z^4 + 468 x^8 y^8 z^4 + 
 1656 t^2 x^8 y^8 z^4 + 2376 t^4 x^8 y^8 z^4 + 1656 t^6 x^8 y^8 z^4 + 
 468 t^8 x^8 y^8 z^4 + 504 x^10 y^8 z^4 + 2448 t^2 x^10 y^8 z^4 + 
 4752 t^4 x^10 y^8 z^4 + 4176 t^6 x^10 y^8 z^4 + 
 1368 t^8 x^10 y^8 z^4 + 108 x^12 y^8 z^4 + 1024 t^2 x^12 y^8 z^4 + 
 3136 t^4 x^12 y^8 z^4 + 3656 t^6 x^12 y^8 z^4 + 
 1436 t^8 x^12 y^8 z^4 + 120 t^2 x^14 y^8 z^4 + 
 760 t^4 x^14 y^8 z^4 + 1280 t^6 x^14 y^8 z^4 + 
 640 t^8 x^14 y^8 z^4 + 44 t^4 x^16 y^8 z^4 + 144 t^6 x^16 y^8 z^4 + 
 104 t^8 x^16 y^8 z^4 + 256 z^6 + 320 t^2 z^6 + 384 t^4 z^6 + 
 352 t^6 z^6 + 160 t^8 z^6 + 272 x^2 z^6 + 256 t^2 x^2 z^6 + 
 1120 t^4 x^2 z^6 + 1408 t^6 x^2 z^6 + 784 t^8 x^2 z^6 + 
 232 x^4 z^6 + 456 t^2 x^4 z^6 + 2104 t^4 x^4 z^6 + 
 2712 t^6 x^4 z^6 + 1856 t^8 x^4 z^6 + 96 x^6 z^6 + 472 t^2 x^6 z^6 + 
 2072 t^4 x^6 z^6 + 3208 t^6 x^6 z^6 + 2792 t^8 x^6 z^6 + 
 24 x^8 z^6 + 298 t^2 x^8 z^6 + 1178 t^4 x^8 z^6 + 2686 t^6 x^8 z^6 + 
 2870 t^8 x^8 z^6 + 108 t^2 x^10 z^6 + 396 t^4 x^10 z^6 + 
 1668 t^6 x^10 z^6 + 2020 t^8 x^10 z^6 + 18 t^2 x^12 z^6 + 
 66 t^4 x^12 z^6 + 726 t^6 x^12 z^6 + 934 t^8 x^12 z^6 + 
 192 t^6 x^14 z^6 + 256 t^8 x^14 z^6 + 24 t^6 x^16 z^6 + 
 32 t^8 x^16 z^6 + 704 y^2 z^6 + 1760 t^2 y^2 z^6 + 
 1888 t^4 y^2 z^6 + 1312 t^6 y^2 z^6 + 480 t^8 y^2 z^6 + 
 1136 x^2 y^2 z^6 + 3456 t^2 x^2 y^2 z^6 + 5152 t^4 x^2 y^2 z^6 + 
 5248 t^6 x^2 y^2 z^6 + 2416 t^8 x^2 y^2 z^6 + 1768 x^4 y^2 z^6 + 
 5200 t^2 x^4 y^2 z^6 + 9152 t^4 x^4 y^2 z^6 + 
 11696 t^6 x^4 y^2 z^6 + 6232 t^8 x^4 y^2 z^6 + 1144 x^6 y^2 z^6 + 
 3760 t^2 x^6 y^2 z^6 + 9984 t^4 x^6 y^2 z^6 + 
 16720 t^6 x^6 y^2 z^6 + 10120 t^8 x^6 y^2 z^6 + 456 x^8 y^2 z^6 + 
 1752 t^2 x^8 y^2 z^6 + 7592 t^4 x^8 y^2 z^6 + 
 16024 t^6 x^8 y^2 z^6 + 10880 t^8 x^8 y^2 z^6 + 72 x^10 y^2 z^6 + 
 544 t^2 x^10 y^2 z^6 + 3952 t^4 x^10 y^2 z^6 + 
 10304 t^6 x^10 y^2 z^6 + 7848 t^8 x^10 y^2 z^6 + 
 72 t^2 x^12 y^2 z^6 + 1160 t^4 x^12 y^2 z^6 + 
 4192 t^6 x^12 y^2 z^6 + 3680 t^8 x^12 y^2 z^6 + 
 128 t^4 x^14 y^2 z^6 + 952 t^6 x^14 y^2 z^6 + 
 1016 t^8 x^14 y^2 z^6 + 96 t^6 x^16 y^2 z^6 + 128 t^8 x^16 y^2 z^6 + 
 384 y^4 z^6 + 1440 t^2 y^4 z^6 + 2208 t^4 y^4 z^6 + 
 1632 t^6 y^4 z^6 + 480 t^8 y^4 z^6 + 608 x^2 y^4 z^6 + 
 3200 t^2 x^2 y^4 z^6 + 6848 t^4 x^2 y^4 z^6 + 6528 t^6 x^2 y^4 z^6 + 
 2272 t^8 x^2 y^4 z^6 + 1760 x^4 y^4 z^6 + 7128 t^2 x^4 y^4 z^6 + 
 15128 t^4 x^4 y^4 z^6 + 16008 t^6 x^4 y^4 z^6 + 
 6248 t^8 x^4 y^4 z^6 + 1288 x^6 y^4 z^6 + 6856 t^2 x^6 y^4 z^6 + 
 19576 t^4 x^6 y^4 z^6 + 25176 t^6 x^6 y^4 z^6 + 
 11168 t^8 x^6 y^4 z^6 + 832 x^8 y^4 z^6 + 4730 t^2 x^8 y^4 z^6 + 
 17242 t^4 x^8 y^4 z^6 + 26382 t^6 x^8 y^4 z^6 + 
 13230 t^8 x^8 y^4 z^6 + 216 x^10 y^4 z^6 + 1980 t^2 x^10 y^4 z^6 + 
 10092 t^4 x^10 y^4 z^6 + 18420 t^6 x^10 y^4 z^6 + 
 10476 t^8 x^10 y^4 z^6 + 274 t^2 x^12 y^4 z^6 + 
 3186 t^4 x^12 y^4 z^6 + 7806 t^6 x^12 y^4 z^6 + 
 5278 t^8 x^12 y^4 z^6 + 384 t^4 x^14 y^4 z^6 + 
 1704 t^6 x^14 y^4 z^6 + 1512 t^8 x^14 y^4 z^6 + 
 144 t^6 x^16 y^4 z^6 + 192 t^8 x^16 y^4 z^6 + 608 x^4 y^6 z^6 + 
 2384 t^2 x^4 y^6 z^6 + 3856 t^4 x^4 y^6 z^6 + 2992 t^6 x^4 y^6 z^6 + 
 912 t^8 x^4 y^6 z^6 + 496 x^6 y^6 z^6 + 3568 t^2 x^6 y^6 z^6 + 
 8848 t^4 x^6 y^6 z^6 + 8976 t^6 x^6 y^6 z^6 + 3200 t^8 x^6 y^6 z^6 + 
 752 x^8 y^6 z^6 + 4356 t^2 x^8 y^6 z^6 + 11780 t^4 x^8 y^6 z^6 + 
 13596 t^6 x^8 y^6 z^6 + 5420 t^8 x^8 y^6 z^6 + 288 x^10 y^6 z^6 + 
 2552 t^2 x^10 y^6 z^6 + 8984 t^4 x^10 y^6 z^6 + 
 12232 t^6 x^10 y^6 z^6 + 5512 t^8 x^10 y^6 z^6 + 
 404 t^2 x^12 y^6 z^6 + 3156 t^4 x^12 y^6 z^6 + 
 5940 t^6 x^12 y^6 z^6 + 3252 t^8 x^12 y^6 z^6 + 
 384 t^4 x^14 y^6 z^6 + 1320 t^6 x^14 y^6 z^6 + 
 1000 t^8 x^14 y^6 z^6 + 96 t^6 x^16 y^6 z^6 + 128 t^8 x^16 y^6 z^6 + 
 288 x^8 y^8 z^6 + 1080 t^2 x^8 y^8 z^6 + 1656 t^4 x^8 y^8 z^6 + 
 1224 t^6 x^8 y^8 z^6 + 360 t^8 x^8 y^8 z^6 + 144 x^10 y^8 z^6 + 
 1008 t^2 x^10 y^8 z^6 + 2448 t^4 x^10 y^8 z^6 + 
 2448 t^6 x^10 y^8 z^6 + 864 t^8 x^10 y^8 z^6 + 
 184 t^2 x^12 y^8 z^6 + 1064 t^4 x^12 y^8 z^6 + 
 1600 t^6 x^12 y^8 z^6 + 720 t^8 x^12 y^8 z^6 + 
 128 t^4 x^14 y^8 z^6 + 376 t^6 x^14 y^8 z^6 + 248 t^8 x^14 y^8 z^6 + 
 24 t^6 x^16 y^8 z^6 + 32 t^8 x^16 y^8 z^6 + 176 z^8 + 256 t^2 z^8 + 
 256 t^4 z^8 + 160 t^6 z^8 + 48 t^8 z^8 + 256 x^2 z^8 + 
 240 t^2 x^2 z^8 + 544 t^4 x^2 z^8 + 496 t^6 x^2 z^8 + 
 192 t^8 x^2 z^8 + 224 x^4 z^8 + 152 t^2 x^4 z^8 + 892 t^4 x^4 z^8 + 
 848 t^6 x^4 z^8 + 396 t^8 x^4 z^8 + 96 x^6 z^8 + 32 t^2 x^6 z^8 + 
 900 t^4 x^6 z^8 + 840 t^6 x^6 z^8 + 516 t^8 x^6 z^8 + 24 x^8 z^8 + 
 8 t^2 x^8 z^8 + 575 t^4 x^8 z^8 + 510 t^6 x^8 z^8 + 
 463 t^8 x^8 z^8 + 210 t^4 x^10 z^8 + 180 t^6 x^10 z^8 + 
 290 t^8 x^10 z^8 + 35 t^4 x^12 z^8 + 30 t^6 x^12 z^8 + 
 123 t^8 x^12 z^8 + 32 t^8 x^14 z^8 + 4 t^8 x^16 z^8 + 256 y^2 z^8 + 
 704 t^2 y^2 z^8 + 784 t^4 y^2 z^8 + 480 t^6 y^2 z^8 + 
 144 t^8 y^2 z^8 + 256 x^2 y^2 z^8 + 1040 t^2 x^2 y^2 z^8 + 
 1632 t^4 x^2 y^2 z^8 + 1424 t^6 x^2 y^2 z^8 + 576 t^8 x^2 y^2 z^8 + 
 416 x^4 y^2 z^8 + 1560 t^2 x^4 y^2 z^8 + 2696 t^4 x^4 y^2 z^8 + 
 2760 t^6 x^4 y^2 z^8 + 1336 t^8 x^4 y^2 z^8 + 224 x^6 y^2 z^8 + 
 1032 t^2 x^6 y^2 z^8 + 2616 t^4 x^6 y^2 z^8 + 3416 t^6 x^6 y^2 z^8 + 
 1992 t^8 x^6 y^2 z^8 + 96 x^8 y^2 z^8 + 472 t^2 x^8 y^2 z^8 + 
 1780 t^4 x^8 y^2 z^8 + 2800 t^6 x^8 y^2 z^8 + 1972 t^8 x^8 y^2 z^8 + 
 88 t^2 x^10 y^2 z^8 + 736 t^4 x^10 y^2 z^8 + 1432 t^6 x^10 y^2 z^8 + 
 1296 t^8 x^10 y^2 z^8 + 140 t^4 x^12 y^2 z^8 + 
 400 t^6 x^12 y^2 z^8 + 548 t^8 x^12 y^2 z^8 + 40 t^6 x^14 y^2 z^8 + 
 136 t^8 x^14 y^2 z^8 + 16 t^8 x^16 y^2 z^8 + 96 y^4 z^8 + 
 384 t^2 y^4 z^8 + 624 t^4 y^4 z^8 + 480 t^6 y^4 z^8 + 
 144 t^8 y^4 z^8 + 64 x^2 y^4 z^8 + 544 t^2 x^2 y^4 z^8 + 
 1472 t^4 x^2 y^4 z^8 + 1568 t^6 x^2 y^4 z^8 + 576 t^8 x^2 y^4 z^8 + 
 448 x^4 y^4 z^8 + 1696 t^2 x^4 y^4 z^8 + 3524 t^4 x^4 y^4 z^8 + 
 3784 t^6 x^4 y^4 z^8 + 1508 t^8 x^4 y^4 z^8 + 224 x^6 y^4 z^8 + 
 1400 t^2 x^6 y^4 z^8 + 4156 t^4 x^6 y^4 z^8 + 5488 t^6 x^6 y^4 z^8 + 
 2508 t^8 x^6 y^4 z^8 + 176 x^8 y^4 z^8 + 992 t^2 x^8 y^4 z^8 + 
 3367 t^4 x^8 y^4 z^8 + 5190 t^6 x^8 y^4 z^8 + 2735 t^8 x^8 y^4 z^8 + 
 264 t^2 x^10 y^4 z^8 + 1578 t^4 x^10 y^4 z^8 + 
 3084 t^6 x^10 y^4 z^8 + 1962 t^8 x^10 y^4 z^8 + 
 315 t^4 x^12 y^4 z^8 + 998 t^6 x^12 y^4 z^8 + 875 t^8 x^12 y^4 z^8 + 
 120 t^6 x^14 y^4 z^8 + 216 t^8 x^14 y^4 z^8 + 24 t^8 x^16 y^4 z^8 + 
 160 x^4 y^6 z^8 + 672 t^2 x^4 y^6 z^8 + 1144 t^4 x^4 y^6 z^8 + 
 912 t^6 x^4 y^6 z^8 + 280 t^8 x^4 y^6 z^8 + 32 x^6 y^6 z^8 + 
 656 t^2 x^6 y^6 z^8 + 2056 t^4 x^6 y^6 z^8 + 2272 t^6 x^6 y^6 z^8 + 
 840 t^8 x^6 y^6 z^8 + 160 x^8 y^6 z^8 + 880 t^2 x^8 y^6 z^8 + 
 2534 t^4 x^8 y^6 z^8 + 3100 t^6 x^8 y^6 z^8 + 1286 t^8 x^8 y^6 z^8 + 
 320 t^2 x^10 y^6 z^8 + 1556 t^4 x^10 y^6 z^8 + 
 2408 t^6 x^10 y^6 z^8 + 1172 t^8 x^10 y^6 z^8 + 
 350 t^4 x^12 y^6 z^8 + 916 t^6 x^12 y^6 z^8 + 598 t^8 x^12 y^6 z^8 + 
 120 t^6 x^14 y^6 z^8 + 152 t^8 x^14 y^6 z^8 + 16 t^8 x^16 y^6 z^8 + 
 72 x^8 y^8 z^8 + 288 t^2 x^8 y^8 z^8 + 468 t^4 x^8 y^8 z^8 + 
 360 t^6 x^8 y^8 z^8 + 108 t^8 x^8 y^8 z^8 + 144 t^2 x^10 y^8 z^8 + 
 504 t^4 x^10 y^8 z^8 + 576 t^6 x^10 y^8 z^8 + 216 t^8 x^10 y^8 z^8 + 
 140 t^4 x^12 y^8 z^8 + 288 t^6 x^12 y^8 z^8 + 148 t^8 x^12 y^8 z^8 + 
 40 t^6 x^14 y^8 z^8 + 40 t^8 x^14 y^8 z^8 + 4 t^8 x^16 y^8 z^8
\end{verbatim}
}

\vskip 0.2in
\bibliography{biblio}

\end{document}